\documentclass[journal]{IEEEtran}
\usepackage{indentfirst}
\usepackage{mathrsfs}
\usepackage{graphicx}
\usepackage{epstopdf}
\usepackage{amsmath, amssymb, amsthm, bm}
\usepackage{leftidx}
\usepackage{extarrows}
\usepackage{stfloats}
\usepackage{picinpar}
\usepackage{enumerate}
\usepackage{algorithm}
\usepackage{algorithmicx}
\usepackage{algpseudocode}
\usepackage{epsfig}
\usepackage{latexsym}
\usepackage{amsfonts}
\usepackage{enumerate}
\usepackage{graphics}
\usepackage{float}
\usepackage{pict2e}
\usepackage{tikz}
\usepackage{lineno}
\usepackage[pdftex,colorlinks,linkcolor=blue,anchorcolor=blue,citecolor=blue]{hyperref}
\usepackage{moreverb}
\usepackage{color}
\usepackage{makecell}
\usepackage{multirow}
\usepackage{float}
\usepackage{subfigure}
\usepackage{bbm}
\usepackage{bbding}

\usepackage[nosort]{cite}

\usepackage{cuted}

\usetikzlibrary{arrows,automata}

\allowdisplaybreaks[4]

\newtheorem{thm}{Theorem}
\newtheorem{ass}{Assumption}

\newtheorem{lem}{Lemma}

\newtheorem{rmk}{Remark}
\newtheorem{defn}{Definition}

\newtheorem{prop}{Proposition}

\newcommand{\tr}{^{\sf T}}
\newcommand{\M}[1]{{\bf{#1}}}

\newcommand{\m}[1]{{\mathrm{#1}}}

\ifCLASSINFOpdf

\else

\fi

\begin{document}
%
\title{Decentralized Inexact Proximal Gradient Method With Network-Independent Stepsizes for Convex Composite Optimization}
%
%
%

\author{Luyao~Guo,
        Xinli~Shi,~\IEEEmembership{Senior Member,~IEEE,}
        Jinde~Cao,~\IEEEmembership{Fellow,~IEEE,} Zihao Wang.
\thanks{This work was supported by the National Natural Science Foundation of China under Grant Nos. 61833005 and 62003084, and the Natural Science Foundation of Jiangsu Province of China under Grant No. BK20200355. (Corresponding author: Jinde Cao.)}

\thanks{L. Guo is with the School of Mathematics and the Frontiers Science Center for Mobile Information Communication and Security, Southeast University, Nanjing 210096, China (e-mail: guo\_luyaoo@163.com).}
\thanks{X. Shi and Z. Wang are with School of Cyber Science \& Engineering, Southeast University, Nanjing 210096, China (xinli\_shi@seu.edu.cn; wwzh0328@seu.edu.cn).}
\thanks{J. Cao is with School of Mathematics, Frontiers Science Center for Mobile Information Communication and Security, Southeast University, Nanjing 210096, China, with Purple Mountain Laboratories, Nanjing 211111, China, and also with Yonsei Frontier Lab, Yonsei University, Seoul 03722, South Korea (e-mail: jdcao@seu.edu.cn).}
}

\maketitle

\begin{abstract}
This paper proposes a novel CTA (Combine-Then-Adapt)-based decentralized algorithm for solving convex composite optimization problems over undirected and connected networks. The local loss function in these problems contains both smooth and nonsmooth terms. The proposed algorithm uses uncoordinated network-independent constant stepsizes and only needs to approximately solve a sequence of proximal mappings, which is advantageous for solving decentralized composite optimization problems where the proximal mappings of the nonsmooth loss functions may not have analytical solutions. For the general convex case, we prove an $O(1/k)$ convergence rate of the proposed algorithm, which can be improved to $o(1/k)$ if the proximal mappings are solved exactly. Furthermore, with metric subregularity, we establish a linear convergence rate for the proposed algorithm. Numerical experiments demonstrate the efficiency of the algorithm.
\end{abstract}


%
\IEEEpeerreviewmaketitle

\section{Introduction}
This paper considers the following decentralized composite optimization over networks with $N$ computational agents:
\begin{align}\label{Problem}
\m{x}^*\in\mathop{\arg\min}\limits_{\m{x}\in \mathbb{R}^n}~\big\{\sum_{i=1}^N (f_i(\m{x})+g_i(\m{x}))\big\},
\end{align}
where the local loss function $f_i:\mathbb{R}^n\rightarrow\mathbb{R}$ and $g_i:\mathbb{R}^n\rightarrow\mathbb{R}\cup\{\infty\}$ are convex and can only be accessed by
corresponding agents. We assume that $f_i$ is $L_i$-smooth, and $g_i$ is proper, closed, but not necessarily smooth. This setting is general and is encountered in various application fields, including multi-agent control, wireless communication, and decentralized machine learning \cite{Nedic20151}.

Arguably, the straightforward approach to address the decentralized settings is to employ the classical gradient descent method combined with the network structure\cite{Nedic2009,Yuan2016,Jako2014,ZhangJ2019}. On one hand, in \cite{Nedic2009,Yuan2016,Jako2014}, the mixing matrix is used to realize the consensus of the local state through the weighted average of neighbor states. These algorithms can be regarded as the quadratic penalty methods for equality consensus constraints, which leads to an inexact convergence when the stepsizes are fixed \cite{Yuan2016}. To obtain the exact optimal solution, the diminishing stepsizes are adopted \cite{Jako2014}. On the other hand, by using the sign of the relative state rather than averaging over the local iterates via the mixing matrix, \cite{ZhangJ2019} proposes an algorithm based on the exact penalty method. However, these primal methods can not achieve the exact convergence under the constant stepsize.

A recent line of work in the primal-dual domain provide a number of alternative methods that can exactly converge to the optimum with the fixed step sizes. When $g_i=0$ and $f_i$ is strongly convex, several linear convergence algorithms have been proposed including D-ADMM\cite{DADMM}, DLM\cite{DLM}, EXTRA\cite{EXTRA}, DIGing\cite{DIGing}, Harnessing\cite{Harnessing}, AugDGM \cite{AugDGM1,AugDGM2}, and Exact Diffusion\cite{ExactDiffusion}. When $g_i\neq0$ and the nonsmooth terms are identical among all the agents, under the assumption that $f_i$ is strongly convex, \cite{Sulaiman2021,Xu2021,Lorenzo2016,Sun2019} present several linear convergence algorithms for problem \eqref{Problem}. When $g_i\neq0$ and the nonsmooth terms may be distinct among these agents, IC-ADMM\cite{Chang2015}, PG-ADMM\cite{Aybat2017}, PG-EXTRA\cite{PGEXTR}, NIDS\cite{NDIS} and PAD\cite{Zhang2021} have been provided. In \cite{Chang2015}, it proves the convergence of IC-ADMM under the strong convexity of $f_i$. In \cite{Aybat2017}, PG-ADMM is shown to have the $O({1}/{k})$ ergodic convergence rate. Based on EXTRA\cite{EXTRA}, PG-EXTRA and NIDS have been provided in \cite{PGEXTR} and \cite{NDIS}, respectively. For PG-EXTRA, it is shown that the algorithm has the $O({1}/{k})$ convergence rate and has the $o({1}/{k})$ convergence rate when the smooth part vanishes. For NIDS, which has a network independent step size, it is proved to have the $o({1}/{k})$ convergence rate. In \cite{Zhang2021}, by the quadratic penalty method and linearized-ADMM, an inexact convergence algorithm called PAD is proposed. For PAD, it proves that for a fixed penalty parameter ${1}/{\tilde{\epsilon}}>0$, the distance between the solution of problem \eqref{Problem} and the limit point of PAD is of order $O(\tilde{\epsilon}^{\frac{1}{2}})$. Moreover, it proves the $O({1}/{k})$ ergodic convergence rate under general convexity, and the linear convergence rate under the strongly convex of $f_i$.

\begin{table*}[!htb]
\renewcommand\arraystretch{1.5}
\begin{center}
\caption{Convergence Properties of Decentralized Algorithms For Convex Composite Optimization On Networks}
\label{Table1}
\scalebox{1}{
\begin{tabular}{|c|c|c|c|c|c|}
\hline
\textbf{Algorithm}&$f_i$&$g_i$&\textbf{Stepsize condition}&\textbf{Rate}&\textbf{Exactness}\\
\hline
\hline
PUDA\cite{Sulaiman2021}&\multirow{3}{*}{strongly convex}&\multirow{3}{*}{\makecell[c]{$g_1=\cdots=g_N$; \\proximal mapping of $g_i$ \\is efficiently computable}}& $\tau<\frac{2-\sigma_M(I-W)}{L_F}$&\multirow{3}{*}{linear}&\multirow{3}{*}{exact}\\
\cline{1-1}\cline{4-4}
ABC-Algorithm\cite{Xu2021}&~&~&$\tau<\frac{1}{L_F+\mu_F}$&~&~\\
\cline{1-1}\cline{4-4}
NEXT\cite{Lorenzo2016}/SONATA\cite{Sun2019}&~&~& $\tau<\min\big\{\frac{\tilde{\mu}_{m}}{\sum_{i=1}^{N}L_i},\varpi\big\}$&~&~\\
\hline
\hline
IC-ADMM\cite{Chang2015}&strongly convex&\multirow{4}{*}{\makecell[c]{proximal mapping of $g_i$ \\is efficiently computable}}&$\tau<\frac{\beta\sigma_m(\mathcal{L}_+)}{L_F^2/\mu_F^2}$&converges&\multirow{4}{*}{exact}\\
\cline{1-2}\cline{4-5}
PG-ADMM\cite{Aybat2017}&\multirow{3}{*}{\makecell[c]{~\\convex}}&~&$\tau_i<\frac{1}{L_i+\beta d_i}$&ergodic: $O(\frac{1}{k})$&~\\
\cline{1-1}\cline{4-5}
PG-EXTRA\cite{PGEXTR}&~&~&$\tau<\frac{1+\sigma_m(W)}{L_F}$& non-ergodic: $O(\frac{1}{k})$&~\\
\cline{1-1}\cline{4-5}
NIDS\cite{NDIS}&~&~&$\tau_i<\frac{2}{L_i}$, $0<\beta\leq0.5/\max_i \tau_i$&non-ergodic: $o(\frac{1}{k})$&~\\
\hline
\hline
\multirow{2}{*}{PAD\cite{Zhang2021}}&convex&\multirow{2}{*}{\makecell[c]{proximal mapping of $g_i$ \\is efficiently computable}}&$\tau<\frac{1}{L_F+\beta\sigma_M(I-W)}$&ergodic: $O(\frac{1}{k})$&\multirow{2}{*}{inexact}\\
\cline{2-2}\cline{4-5}
~&strongly convex&~&$\tau<\frac{1}{L_F^2/\mu_F+\beta\sigma_M(I-W)}$&linear&~\\
\hline
\hline
\multirow{2}{*}{\textbf{D-iPGM}}&\multicolumn{2}{c|}{convex}&\multirow{2}{*}{$\tau_i<\frac{2}{L_i}$, $0<\beta\leq0.5/\max_i \tau_i$}&\makecell[c]{ergodic: $O(\frac{1}{k})$\\non-ergodic: $o(\frac{1}{k})$}&\multirow{2}{*}{exact}\\
\cline{2-3}\cline{5-5}
&\multicolumn{2}{c|}{metrically subregular}&~&linear&\\
\hline
\hline
\multicolumn{6}{|c|}{\makecell[c]{$\tau$ and $\beta$ are the stepsizes of primal update and dual update, respectively. $\tau_i$ is the primal stepsize of local agent.\\
$L_i$ is the Lipschitz constant of the gradients $\nabla f_i$ and $L_F=\max_i\{L_i\}$. $\mu_F$ is the strong convexity constant of the functions $f_i$.\\
$\tilde{\mu}_m$ is the strong convexity parameter of the successive convex approximation surrogate of $f_i$. $\varpi>0$ is a given constant. \\
$\mathcal{L}_+$ is the signless Laplacian matrix. The definition of $W$, $\sigma_m(\cdot)$, and $\sigma_M(\cdot)$ can be found in Assumption \ref{ass2} and Notations.}}\\
\hline
\end{tabular}}
\end{center}
\end{table*}

{These algorithms can be divided into ATC (Adapt-Then-Combine)- and CTA (Combine-Then-Adapt)-based decentralized algorithms according to their characteristics \cite{Xu2021}. For ATC-based decentralized algorithms, there are several algorithms with network-independent constant stepsizes such as AugDGM \cite{AugDGM1,AugDGM2}, Exact Diffusion\cite{ExactDiffusion}, ABC-Algorithm \cite{Xu2021}, NEXT\cite{Lorenzo2016}/SONATA\cite{Sun2019}, and NIDS \cite{NDIS} (see \cite{Xu2021} for more details). However, to
the best of our knowledge, the choice of the stepsize of the existing CTA-based decentralized algorithms are related to the network topology. }Moreover, in the algorithms mentioned above for composite optimization, the proximal mapping of $g_i$ is required to have analytic solutions or can be efficiently computed. However, in many applications, there is no analytic solution to the proximal mapping of $g_i$. Such problems include generalized LASSO regularizer \cite{Example3}, fused LASSO regularizer \cite{Tibshirani2005}, the octagonal selection and clustering algorithm for regression (OSCAR) regularizer \cite{Howard2008}, and the group LASSO regularizer \cite{Jacob2009}. Therefore, it prompts us to introduce an inexact inner-solver to compute an approximation of the proximal mapping. There are several works focusing on inexact versions of ALM \cite{MP2021}, primal-dual hybrid gradient method \cite{Han2021}, ADMM \cite{YaoW2017} and decentralized ALM \cite{IDEAL2020,Hong2016}. To our knowledge, the corresponding results for problem \eqref{Problem} are relatively scarce. Furthermore, to the problem \eqref{Problem}, when the nonsmooth parts $g_i$ are distinct, the linear convergence rate of corresponding algorithms under certain assumptions have not been established. Although \cite{Sulaiman2021} shows that when $f_i$ is a strongly convex smooth function and $g_i$ is nonsmooth convex with closed-form solution to proximal mappings, then the dimension independent exact global linear convergence cannot be attained for the composite optimization problem \eqref{Problem}, it does not mean that the linear convergence rate related to dimension is impossible under some conditions. Therefore, it naturally leads to the following three problems.

\textbf{Question 1.} Can one provide a CTA-based decentralized algorithm with network-independent constant stepsizes?

\textbf{Question 2.} Can one provide a decentralized algorithm when the proximal mapping of $g_i$ has no analytic solution?

\textbf{Question 3.} Can one establish the linear convergence of the proposed decentralized algorithm for the problem \eqref{Problem}?

This paper aims at addressing Questions 1, 2 and 3 for the problem \eqref{Problem} over an undirected and connected network. We develop a novel CTA-based decentralized inexact proximal gradient method (D-iPGM) with network-independent constant stepsizes. In the implementation of D-iPGM, computable approximations of the proximal mapping is allowed, which implies that the proposed algorithm can reduce the computational cost when proximal mapping of $g_i$ has no closed-form solution. To convergence analysis, we prove the $O(1/k)$ convergence rate under general convexity. When the proximal mapping of $g_i$ is solved exactly, we establish the $o(1/k)$ non-ergodic convergence of D-iPGM. We further prove the linear convergence rate of D-iPGM under the metric subregularity which is weaker than strong convexity. Table \ref{Table1} shows a comparison of D-iPGM with the existing decentralized algorithms to convex composite optimization problems over networks.

This paper is organized as follows. Section \ref{Sec3} casts the problem \eqref{Problem} into a constrained form, and based on the reformulation, the D-iPGM is proposed. Additionally, some relations between D-iPGM and PG-EXTRA \cite{PGEXTR} and NIDS \cite{NDIS} are discussed. Then, the convergence properties under general convexity and metric subregularity are investigated in Section \ref{Sec4}. Finally, the numerical experiments are implemented in Section \ref{Sec5} and conclusions are given in Section \ref{Sec6}.

\emph{Notations:} $\mathbb{R}^n$ denotes the $n$-dimensional vector space with inner-product $\langle \cdot,\cdot\rangle$. The $l_1$-norm, $l_2$-norm and $l_{\infty}$-norm are denoted as $\|\cdot\|_1$, $\|\cdot\|$ and $\|\cdot\|_{\infty}$, respectively. $\mathbf{0}$, $\mathbf{I}$, and $\mathbf{1}$ denote the null matrix, the identity matrix, and the all-ones matrix of appropriate dimensions, respectively. For $\m{x}_i\in \mathbb{R}^n,i=1\ldots,m$, let $\mathrm{col}\{\m{x}_1,\cdots,\m{x}_N\}=[\m{x}_1\tr,\cdots,\m{x}_N\tr]\tr$. For a matrix $A\in \mathbb{R}^{m\times n}$, $\sigma_m(A)$ and $\sigma_M(A)$ is the smallest and largest nonzero singular value of $A$, respectively. If $A$ is symmetric, $A\succ 0$ means that $A$ is positive definite. Let $\mathcal{H}\succ0$, and $\mathcal{D}$ be a non-empty closed convex subset of $\mathbb{R}^n$. For any $\m{x}\in\mathbb{R}^n$, we define $\|\m{x}\|^2_{\mathcal{H}}:=\langle \m{x},\mathcal{H}\m{x}\rangle$, $\mathrm{dist}_{\mathcal{H}}(\m{x},\mathcal{D}):=\min_{\m{x}'\in\mathcal{D}}\|\m{x}-\m{x}'\|_{\mathcal{H}}$, and  $\mathcal{P}_{\mathcal{D}}^{\mathcal{H}}(\m{x}):=\arg\min_{\m{x}'\in \mathcal{D}}\|\m{x}-\m{x}'\|_{\mathcal{H}}$. For a given closed proper convex function $\theta$, the proximal mapping of $\theta$ relative to $\|\cdot\|_{\mathcal{H}}$ is defined as
$\mathrm{prox}_\theta^{\mathcal{H}}(\m{x}):=\arg\min_\m{y}\{\theta(\m{y})+\frac{1}{2}\|\m{y}-\m{x}\|_{\mathcal{H}}^2\}$. When $\mathcal{H}=\M{I}$, we just omit it from these notations. $\mathcal{B}_{r}(\bar{\m{x}}):=\{\m{x}:\|\m{x}-\bar{\m{x}}\|<r\}$ denotes the open Euclidean norm ball around $\m{x}$ with modulus $r>0$.

\section{Decentralized Inexact Proximal Gradient Method}\label{Sec3}
Consider an undirected and connected network $\mathcal{G}(\mathcal{V},\mathcal{E})$, where $\mathcal{V}=\{v_1,\ldots,v_N\}$ denotes the vertex set, and the edge set $\mathcal{E}\subseteq \mathcal{V}\times \mathcal{V}$ specifies the connectivity in the network, i.e., a communication link between agents $i$ and $j$ exists iff $(i,j)\in \mathcal{E}$. Let $\mathrm{x}_i\in \mathbb{R}^n$ be the decision variable held by the agent $i$. To cast the problem \eqref{Problem} into a consensus-constrained form and make sure that $\mathrm{x}_1=\cdots=\mathrm{x}_N$, we introduce the global mixing matrix $\mathbf{W}=W\otimes \mathbf{I}_n$, where $W=[W_{ij}]\in \mathbb{R}^{N\times N}$, and make the following standard assumption in decentralized optimization.
\begin{ass}\label{ass2}
The mixing matrix $W$ is symmetric and doubly stochastic, and the null space of $\mathbf{I}-\mathbf{W}$ is $\mathrm{Span}(\mathbf{1})$. Moreover, if $(i, j) \notin \mathcal{E}$ and $i \neq j$, $W_{i j}=W_{j i}=0$; otherwise, $W_{ij}>0$.
\end{ass}
In addition, with this mixing matrix, similar as \cite{PGEXTR}, we introduce another mixing matrix $\tilde{\M{W}}$ defined as $\tilde{\M{W}}:=\frac{1}{2}(\M{I+W})$. By Assumption \ref{ass2}, the constraint $\mathrm{x}_{1}=\cdots=\mathrm{x}_{N}$ is equivalent to the identity $(\mathbf{I}-\mathbf{W})\mathbf{x}=\mathbf{0}$, which is again equivalent to $\sqrt{\mathbf{I}-\mathbf{W}} \mathbf{x}=\mathbf{0}$. Then, we have the following equivalent formulation of problem \eqref{Problem}.
\begin{align}\label{Problem2}
\min_{\M{x}}~\underbrace{\sum_{i=1}^N f_i(\m{x}_i)}_{:=F(\M{x})}+\underbrace{\sum_{i=1}^N g_i(\m{x}_i)}_{:=G(\M{x})},~\mathrm{s.t.~} \sqrt{\mathbf{I}-\mathbf{W}}\mathbf{x}=\mathbf{0},
\end{align}
where $\mathbf{x}=[\mathrm{x}_1\tr,\cdots,\mathrm{x}_N\tr]\tr\in\mathbb{R}^{Nn}$. To guarantee the wellposedness of \eqref{Problem}, it is necessary to suppose that the optimal solution of \eqref{Problem} exists. Throughout the paper, we give the following assumption.
\begin{ass}\label{ass1}
The local function $f_i$ is convex and $L_i$-smooth.
\end{ass}
According to \cite[Theorem 5.8]{Beck2017} and Assumption \ref{ass1}, for any $\M{x},\M{y},\M{z}\in\mathbb{R}^{Nn}$, we have the following commonly used inequalities
\begin{align}
&\langle \nabla F(\M{x})-\nabla F(\M{y}), \M{x-y} \rangle\geq \| \nabla F(\M{x})-\nabla F(\M{y})\|^2_{\M{L}_F^{-1}},\label{SM3}\\
&\langle \M{x}-\M{y},\nabla F(\M{z})\rangle\leq F(\M{x})-F(\M{y})+\frac{1}{2}\|\M{y}-\M{z}\|_{\M{L}_F}^2,\label{SM2}\\
&\langle \M{x}-\M{y},\nabla F(\M{z})-\nabla F(\M{x})\rangle\leq \frac{1}{4} \|\M{y}-\M{z}\|_{\M{L}_F}^2,\label{SM1}
\end{align}
where $\M{L}_F=\mathrm{diag}\{L_1\M{I}_n,\cdots,L_N\M{I}_n\}$.

\subsection{Algorithm Development}
Recall $\m{ABC}$-unified framework proposed in \cite{Xu2021}
\begin{subequations}\label{ABC-1}
\begin{align}
\M{s}^{k}&=\m{A}\M{x}^k-\tau\m{B}\nabla F(\M{x}^k)-\bm{\lambda}^k,\\
\bm{\lambda}^{k+1}&=\bm{\lambda}^k+\m{C}\M{s}^{k},\\
\M{x}^{k+1}&=\mathrm{prox}_{\tau G}(\M{s}^{k}),
\end{align}
\end{subequations}
where $\m{A,B,C}$ are suitably chosen weight matrices, $\tau>0$ is the primal stepsize, and the nonsmooth function $g_i$ is common to all agents. Let $\m{C}=\frac{1}{2}(\M{I-W})$. According to this unified framework, when $\m{A}=\m{B}\neq\M{I}$, by eliminating the dual variable $\bm{\lambda}$, it holds that
\begin{align*}
\mathbf{x}^{k+1}&=\mathrm{prox}_{\tau G}(\mathbf{s}^{k}),\\
\M{s}^{k+1}&=\tilde{\M{W}}\M{s}^{k}+\underbrace{\m{A}(\underbrace{\M{x}^{k+1}-\M{x}^k+\tau\nabla F(\mathbf{x}^{k})-\tau\nabla F(\M{x}^{k+1})}_{\text{First: Adapt}})}_{\text{Then: Combine}},
\end{align*}
and these algorithms are the ATC-based decentralized algorithms. Similarly, when $\m{B}=\M{I}$, it deduces that
\begin{align*}
\mathbf{x}^{k+1}&=\mathrm{prox}_{\tau G}(\mathbf{s}^{k}),\\
\M{s}^{k+1}&=\tilde{\M{W}}\M{s}^{k}+\underbrace{\underbrace{\m{A}(\M{x}^{k+1}-\M{x}^k)}_{\text{First: Combine}}+\tau\nabla F(\M{x}^{k})-\tau\nabla F(\M{x}^{k+1})}_{\text{Then: Adapt}},
\end{align*}
and these algorithms are the CTA-based decentralized algorithms. To ATC-based algorithms, the communication of state information and local gradient information are involved. By this strategy, the stepsize conditions of ATC-based algorithms are often independent of the network topology, such as NIDS \cite{NDIS}, which exchanges the gradient adapted estimations, i.e., $2\M{x}^{k+1}-\M{x}^k-\tau(\nabla F(\M{x}^{k+1})-\nabla F(\M{x}^k))$, and the stepsize condition of NIDS is $0<\tau<2/{L_F}$, where $L_F=\max_{i\in\mathcal{V}}\{L_i\}$. To ATC-based algorithms, they only need to share state information; however, their stepsize conditions are often related to the network topology, such as PG-EXTRA \cite{PGEXTR}, which exchanges only the state information, and the stepsize condition PG-EXTRA is $0<\tau<(1+\sigma_m(\M{W}))/{L_F}$. To blend the respective superiority of these two types of algorithms, we propose the following decentralized algorithm to problem \eqref{Problem}.
\begin{subequations}\label{Ca8}
\begin{align}
\tilde{\mathbf{x}}^{k}&=\mathrm{prox}_G^{\Gamma^{-1}}(\mathbf{x}^k-\Gamma(\nabla F(\mathbf{x}^k)-\bm{\lambda}^k)),\label{Ca81}\\
\bm{\lambda}^{k+1}&=\bm{\lambda}^k-\beta (\mathbf{I-W})\tilde{\mathbf{x}}^{k},\\
\M{x}^{k+1}&=\tilde{\mathbf{x}}^{k}-\Gamma(\bm{\lambda}^{k}-\bm{\lambda}^{k+1})\label{correctionstep},
\end{align}
\end{subequations}
where $\Gamma=\mathrm{diga}\{\tau_1\M{I}_n,\cdots,\tau_N\M{I}_n\}$ is the primal stepsize matrix, $\beta>0$ is the dual stepsize, $\bm{\lambda}^k=\sqrt{\mathbf{I-W}}\bm{\alpha}^k$ is the ``scaled'' Lagrange multiplier, and $\bm{\alpha}^k$ is the Lagrange multiplier. To achieve a network independent stepsize, we use the simple primal correction step \eqref{correctionstep}. In Section \ref{Sec4}, we will prove that the algorithm \eqref{Ca8} is convergent when $0<\tau_i<{2}/{L_i},i\in\mathcal{V}$, $0<\beta\leq0.5/\max_i \tau_i$. By eliminating the dual variable $\bm{\lambda}$, it holds that
\begin{align*}
\tilde{\M{x}}^{k}&=\mathrm{prox}^{\Gamma^{-1}}_{G}(\M{s}^{k}),\\
\M{x}_c^k&=\beta\Gamma(\M{I-W})\tilde{\M{x}}^{k},\\
\M{x}^{k+1}&=\tilde{\M{x}}^{k}-\M{x}_c^k,\\
\M{s}^{k+1}&=\M{s}^{k}-\M{x}^{k}+\tilde{\M{x}}^{k}-2\M{x}_c^k+\Gamma\nabla F(\mathbf{x}^{k})-\Gamma\nabla F(\mathbf{x}^{k+1}).
\end{align*}
Apparently, the implementation of the proposed algorithm only requires neighboring variables $\tilde{\m{x}}_j^k$, and there is only one round of communication in each iteration. Compared with \eqref{ABC-1}, we know that the proposed algorithm is indeed a CTA-based decentralized algorithm and is independent of this primal-dual algorithmic framework.

On the other hand, when the proximal mapping of $g_i$ has no analytical solution, solving subproblem \eqref{Ca81} very accurately is a waste. It is more efficient to develop a
proper condition and stop the subproblem procedure, i.e., inner iterations, once the condition is satisfied. By introducing an absolutely summable error criteria for subproblem \eqref{Ca81}, we give D-iPGM (see Algorithm \ref{alg1}).
\begin{algorithm}[!t]
  \caption{D-iPGM} 
    \label{alg1}
  \begin{algorithmic}[1]
    \Require
      Mixing matrix $W$, $0<\tau_i<{2}/{L_i}$, and $0<\beta\leq0.5/\max_i \tau_i$. Let $\{\varepsilon_k\}$ be a summable sequence of nonnegative numbers.
    \State Initial $\mathbf{x}^0\in \mathbb{R}^{Nn}$, $\bm{\lambda}^0=\mathbf{0}$.
    \For{$k=1,2,\ldots,K$}
      \State Primal Update:
      \begin{align}\label{Subx}
\tilde{\mathrm{x}}_i^{k}\approx \tilde{\mathrm{y}}_i^k=\mathrm{prox}_{\tau_ig_i}(\mathrm{x}_i^k-\tau_i(\nabla f_i(\mathrm{x}_i^k)-{\lambda}_i^k)),\end{align}
such that there exists a vector
\begin{subequations}\label{Condition1}
\begin{align}
&\m{d}_i^k\in \partial g_i(\tilde{\mathrm{x}}_i^{k})+\nabla f_i(\mathrm{x}_i^k)-\lambda_i^k+\frac{1}{\tau_i}(\tilde{\mathrm{x}}_i^{k}-\mathrm{x}_i^{k}),\label{Condition1.1}\\
&\text{satisfying } \|\m{d}_i^k\|\leq  \varepsilon_k. \label{Condition1.2}
\end{align}
\end{subequations}
      \State Agent $i$ exchanges $\tilde{\mathrm{x}}_i^{k}$ with neighbors.
      \State Dual Update:
      \begin{align}\label{Suby}
      \lambda^{k+1}_i={\lambda}_i^{k}-\beta (\tilde{\mathrm{x}}_i^{k}-\sum\nolimits_{j=1}^{N}W_{ij}\tilde{\mathrm{x}}_j^{k}).
      \end{align}
      \State Correction Step:
      $\mathrm{x}_i^{k+1}=\tilde{x}^k_i-\tau_i(\lambda_i^k-\lambda_i^{k+1}).$
    \EndFor
    \Ensure
      $\mathbf{x}^K$.
  \end{algorithmic}
\end{algorithm}

In Algorithm \ref{alg1}, to save computational expenses, in the $k$-th iteration of D-iPGM, inexact proximal gradient descent is allowed, i.e., we can solve subproblem \eqref{Subx} inexactly by finding an approximate solution. It indicates that we have to find an appropriate inner iterative method to solve it inexactly. This is possible for a majority of the applications. For instance, to general LASSO regularizer \cite{Example3}, fused LASSO regularizer \cite{Tibshirani2005}, OSCAR regularizer \cite{Howard2008}, and the group LASSO regularizer \cite{Jacob2009}, in which $g_i$ has the form as $g_i(\m{x})=r_{1,i}(\m{x})+r_{2,i}(\m{D}_i\m{x}_i)$ with $r_{1,i},r_{2,i}$ being nonsmooth functions and $\m{D}_i$ being a given matrix, we can use the operator splitting methods \cite{PDS2021} or accelerated linearized ADMM \cite{Xuyangyang2017} to solve \eqref{Subx} with the inner stopping criterion \eqref{Condition1}.

\begin{rmk}
Since $g_i$ is convex, closed, and proper, by the definition of proximal mapping, each of the subproblems \eqref{Subx} is strongly convex so that each $\m{y}_i^k$ is uniquely determined by $(\m{x}_i^k,\lambda_i^k)$. From \cite{YaoW2017}, it holds that for any $\varepsilon_k\geq 0$, a certain $\tilde{\m{x}}^k$ can be found such that the absolutely summable error criteria \eqref{Condition1} holds. Therefore, the Algorithm \ref{alg1} is well-defined. Moreover, by \cite[Lemma 3.2]{Han2021}, the inner loops of Algorithm \ref{alg1} always terminate after a finite number of iterations.
\end{rmk}
\begin{rmk}
It follows from \eqref{Condition1.1} and the definition of proximal mapping that $\tilde{\m{x}}_i^k=\tilde{\m{y}}_i^k$ if $\m{d}^k_i=\M{0}$. Therefore, $\tilde{\m{x}}_i^k$ can be seen as an approximation of the exact solution $\tilde{\m{y}}_i^k$. To simplify convergence analysis, we give the inner stopping criterion \eqref{Condition1}. However, in practice, we can alternatively use $\|\m{x}_i^{k,l+1}-\m{x}_i^{k,l}\|\leq\varepsilon_0\varepsilon_k$ as the inner stopping criterion, where $\varepsilon_0>0$ is a constant, and $\m{x}_i^{k,l}$ is the estimate of $\tilde{\m{y}}_i^k$ in $l$-the inner iteration, because $\m{d}_i^{k,l+1}\in \partial g_i(\m{x}_i^{k,l+1})+\nabla f_i(\mathrm{x}_i^k)-\lambda_i^k+\frac{1}{\tau_i}(\m{x}_i^{k,l+1}-\mathrm{x}_i^{k})$ can be controlled by $\|\m{x}_i^{k,l+1}-\m{x}_i^{k,l}\|$, i.e., there exists a constant $\varepsilon_0$ such that $\|\m{d}_i^{k,l+1}\|\leq \varepsilon_0 \|\m{x}_i^{k,l+1}-\m{x}_i^{k,l}\|$. We will specify this in Section \ref{Sec5}.
\end{rmk}
\begin{rmk}
Let $\M{d}^k=\mathrm{col}\{\m{d}_1^k,\cdots,\m{d}_N^k\}$. Different from \cite{Han2021}, we require $\|\m{d}_i^k\|\leq  \varepsilon_k$ rather than $\|\mathbf{d}^k\|\leq  \frac{\varepsilon_k}{\max \{\beta,\|\tilde{\mathbf{x}}^k\|\}}$. Thus, by this setting, all agents have their own independent subproblems and independent error criteria for these subproblems, which implies that these subproblems can be solved in a decentralized manner. In addition, to control the residual $\m{d}_i^k$, we introduce a summable sequence $\{\varepsilon_k\}$, which implies that we can choose $\varepsilon_k$ by $\varepsilon_k=\varepsilon_0/(k+1)^r$ with any $\varepsilon_0>0$ and $r>1$, or by $\varepsilon_k=r^{k}$ with $0<r<1$. Apparently, when the proximal mapping of $g_i$ has a closed-form representation, we can set $\varepsilon_k\equiv0,k\geq0$. In this scenario, we use the exact proximal gradient descent to perform the primal update, thus it has the same complexity as the existing proximal gradient methods \cite{Sulaiman2021,Xu2021,Lorenzo2016,Sun2019,Chang2015,Aybat2017,PGEXTR,NDIS,Zhang2021}.
\end{rmk}

\begin{table*}[!h]
\renewcommand\arraystretch{1.5}
\begin{center}
\caption{The Comparison of D-iPGM ($\varepsilon_k=0$) with PG-EXTRA \cite{PGEXTR} and NIDS \cite{NDIS}.}
\scalebox{0.95}{
\begin{tabular}{|c|c|c|c|c|c|c|c|}
\hline
\textbf{Algorithm}& \textbf{Communication}&\textbf{Stepsize, uncoordinated?}&\textbf{Non-ergodic Rate}&\textbf{Ergodic Rate}&\textbf{Linear Rate}&\textbf{Inexact Iteration}\\
\hline
PG-EXTRA \cite{PGEXTR}&$2\M{x}^{k+1}-\M{x}^k$&$\tau<(1+\sigma_m(\M{W}))/{L_F}$, no&$O(1/k)$&No results&No results&\XSolidBrush\\
\hline
\multirow{2}{*}{NIDS \cite{NDIS}}&$\tau(\nabla F(\M{x}^{k})-\nabla F(\M{x}^{k+1}))$&\multirow{2}{*}{$\tau_i<2/L_i$, yes}&\multirow{2}{*}{$o(1/k)$}&\multirow{2}{*}{No results}&\multirow{2}{*}{No results}&\multirow{2}{*}{\XSolidBrush}\\
~&$+2\M{x}^{k+1}-\M{x}^k$~~~~~~~~~~~~~~~&~&~&~&~&~\\
\hline
D-iPGM, $\varepsilon_k=0$&$\tilde{\M{x}}^{k}$&$\tau_i<2/L_i$, yes&$o(1/k)$&$O(1/k)$&\Checkmark&\Checkmark\\
\hline
\end{tabular}}
\label{Comparison:2}
\end{center}
\end{table*}

\subsection{Comparison With PG-EXTRA \cite{PGEXTR} and NIDS \cite{NDIS}}
In this subsection, we compare D-iPGM ($\varepsilon_k=0$), PG-EXTRA \cite{PGEXTR}, and NIDS \cite{NDIS} from the perspective of operator splitting. Recall problem \eqref{Problem2}, which is equivalent to
\begin{align}\label{TSproblem1}
\min_{\M{x}} F(\M{x})+G(\M{x})+\delta_{\{\M{0}\}}(\M{Vx}),
\end{align}
where $\M{V}=\sqrt{\M{I-W}}$, and $\delta_{\{\M{0}\}}(\M{V}\M{x})$ is an indicator function defined as $\delta_{\{\M{0}\}}(\M{V}\M{x})=\M{0}$ if $\M{V}\M{x}=\M{0}$; otherwise, $\delta_{\{\M{0}\}}(\M{V}\M{x})=\infty$, which encodes the constraint $\M{V}\M{x}=\M{0}$. For compactness of exposition, we define the following operators
$$
\bm{A}=\left(
         \begin{array}{cc}
           \partial G & \M{0} \\
           \M{0} & \M{0} \\
         \end{array}
       \right), \bm{B}=\left(
         \begin{array}{cc}
           \M{0} & -\M{V} \\
           \M{V} & \M{0} \\
         \end{array}
       \right),\bm{C}=\left(
         \begin{array}{cc}
           \nabla F & \M{0} \\
           \M{0} & \M{0} \\
         \end{array}
       \right).
$$

According to \cite[Section 11.3.3]{BookYin2022}, PG-EXTRA is equivalent to the C-V splitting algorithm \cite{CV1,CV2} applied to problem \eqref{TSproblem1}. In addition, it also can be derived by forward-backward operator splitting \cite{PDS2021} under the metric defined by
\begin{align*}
\M{H}_{\mathrm{P}}=
\left(
  \begin{array}{cc}
    \frac{1}{\tau}\M{I} & \M{V} \\
    \M{V} & 2\tau\M{I} \\
  \end{array}
\right),
\end{align*}
i.e., it can be rewritten as
$$
\M{u}^{k+1}=(\M{H}_{\m{P}}+\bm{A+B})^{-1}(\M{H}_{\m{P}}-\bm{C})\M{u}^k,
$$
where $\M{u}^k=(\M{x}^k,\bm{\alpha}^k)$. Since $\bm{\lambda}^k=\sqrt{\M{I-W}}\bm{\alpha}^k$, the primal-dual update of PG-EXTRA is
\begin{subequations}\label{PD-PG-EXTRA}
\begin{align}
\M{x}^{k+1}&=\mathrm{prox}_{\tau G}(\mathbf{x}^k-\tau(\nabla F(\mathbf{x}^k)-\bm{\lambda}^k)),\label{PD-PG-EXTRA11}\\
\bm{\lambda}^{k+1}&=\bm{\lambda}^k-\frac{1}{2\tau}(\M{I-W})(2\M{x}^{k+1}-\M{x}^k).
\end{align}
\end{subequations}
By eliminating the dual variable $\bm{\lambda}$, it holds that
\begin{subequations}\label{IT:EXTRA2}
\begin{align}
\mathbf{x}^{k+1}&=\mathrm{prox}_{\tau G}(\mathbf{s}^{k}),\label{IT:EXTRA22}\\
\mathbf{s}^{k+1}&=\mathbf{s}^{k}-\M{x}^{k+1}+\tilde{\M{W}}(2\M{x}^{k+1}-\M{x}^{k})\nonumber\\
&\quad +\tau\nabla F(\M{x}^{k})-\tau\nabla F(\mathbf{x}^{k+1}).
\end{align}
\end{subequations}
In the implementation of PG-EXTRA, only the communication of the state variable $\m{x}_i$ is involved. Since the eigenvalues of $\M{W}$ lie in $(-1,1]$, and the multiplicity of eigenvalue $1$ is one, we obtain that the metric $\M{H}_{\m{P}}$ is positive definite. To ensure the convergence of PG-EXTRA, the positive definiteness of the following metric matrix is required
\begin{align*}
\widehat{\M{H}}_{\mathrm{P}}=\M{H}_{\mathrm{P}}-\frac{1}{2}\mathrm{diag}\{L_F\M{I},\M{0}\},
\end{align*}
where $L_F=\max_{i\in\mathcal{V}}\{L_i\}$. Thus, we have that the stepsize condition of PG-EXTRA is $0<\tau<(1+\sigma_m(\M{W}))/{L_F}$.

To NIDS, in terms of \cite[Section 11.3.3]{BookYin2022}, it is equivalent to the primal-dual three-operator splitting algorithm (PD3O) \cite{PD3O} applied to problem \eqref{TSproblem1} with the metric matrix being
$$
\M{H}_{\m{N}}=\left(
                \begin{array}{cc}
                  \frac{1}{\tau}\M{I} & \M{0} \\
                  \M{0} & \tau(2\M{I}-\M{V}^2) \\
                \end{array}
              \right).
$$
The primal-dual update of NIDS is
\begin{subequations}\label{IT:NIDS1}
\begin{align}
\M{x}^{k+1}&=\mathrm{prox}_{\tau G}(\mathbf{x}^k-\tau(\nabla F(\mathbf{x}^k)-\bm{\lambda}^k)),\label{IT:NIDS11}\\
\bm{\lambda}^{k+1}&=\bm{\lambda}^k-\frac{1}{2\tau}(\M{I-W})\big(2\M{x}^{k+1}-\M{x}^k\nonumber\\
&\quad+\tau\nabla F(\M{x}^{k})-\tau\nabla F(\M{x}^{k+1})\big).
\end{align}
\end{subequations}
By eliminating the dual variable $\bm{\lambda}$, it holds that
\begin{subequations}\label{IT:NIDS2}
\begin{align}
\mathbf{x}^{k+1}&=\mathrm{prox}_{\tau G}(\mathbf{s}^{k}),\label{IT:NIDS22}\\
\mathbf{s}^{k+1}&=\mathbf{s}^{k}-\M{x}^{k+1}+\tilde{\M{W}}\big(2\M{x}^{k+1}-\M{x}^{k}\nonumber\\
&\quad +\tau\nabla F(\M{x}^{k})-\tau\nabla F(\mathbf{x}^{k+1})\big).
\end{align}
\end{subequations}
To achieve a more relaxed convergence condition than PG-EXTRA, there are two additional terms $(\nabla F(\M{x}^{k})$ and $\nabla F(\M{x}^{k+1}))$ compared with the primal-dual version of PG-EXTRA \eqref{PD-PG-EXTRA} in the dual update. Moreover, in the implementation of NIDS, the communication of both the state variable $\m{x}_i$ and $\nabla f_i$ are involved. Since $2\M{I}-\M{V}^2=\M{I+W}\succ0$, the metric $\M{H}_{\m{N}}$ is positive definite. Similarly, to ensure the convergence of NIDS, the metric matrix
$$
\widehat{\M{H}}_{\mathrm{N}}=\M{H}_{\m{N}}-\frac{1}{2}\mathrm{diag}\{L_F\M{I},\M{0}\},
$$
is required to be positive definite. Thus, the stepsize condition of NIDS is $0<\tau<{2}/{L_F}$.

Inspired by asymmetric forward-backward-adjoint splitting \cite{AFBA}, which is a very general three operator splitting scheme, and by prediction-correction framework \cite{He2012}, the proposed D-iPGM can be seen as the triangularly preconditioned forward-backward operator splitting algorithm with a primal corrector. For a clearer view of this, we define the triangularly preconditioner $\M{Q}_{\m{D}}$ and the corrector $\M{M}_{\m{D}}$
$$
\M{Q}_{\m{D}}=\left(
        \begin{array}{cc}
          \frac{1}{\tau}\M{I} & \M{V} \\
           \mathbf{0} & 2\tau\M{I} \\
        \end{array}
      \right),~\M{M}_{\m{D}}=\left(
                       \begin{array}{cc}
                         \M{I} & \tau \M{V} \\
                         \M{0} & \M{I} \\
                       \end{array}
                     \right).
$$
Let $\varepsilon_k=0$, $\Gamma=\tau\M{I}, \tau>0$ and $\beta=\frac{1}{2\tau}$. By these definitions, D-iPGM is equivalent to
\begin{align*}
\text{[D-iPGM]}:~~~\tilde{\M{u}}^k&=(\M{Q}_{\m{D}}+\bm{A+B})^{-1}(\M{Q}_{\m{D}}-\bm{C})\M{u}^k,\\
\M{u}^{k+1}&=\M{u}^k-\M{M}_{\m{D}}(\M{u}^k-\tilde{\M{u}}^k).
\end{align*}
By eliminating the dual variable, it holds that
\begin{align*}
\tilde{\M{x}}^{k}&=\mathrm{prox}_{\tau G}(\M{s}^{k}),\\
\M{x}^{k+1}&=\tilde{\M{W}}\tilde{\M{x}}^{k},\\
\M{s}^{k+1}&=\M{s}^{k}-\M{x}^{k}-\M{W}\tilde{\M{x}}^{k}+\tau\nabla F(\mathbf{x}^{k})-\tau\nabla F(\mathbf{x}^{k+1}).
\end{align*}
Different from forward-backward operator splitting, the preconditioner $\M{Q}_{\m{D}}$ is not required to be positive definite. To achieve a network-independent convergence condition, the correction step
$\M{u}^{k+1}=\M{u}^k-\M{M}_{\m{D}}(\M{u}^k-\tilde{\M{u}}^k)$
is employed, which enables the metric matrices of the D-PGM,
\begin{align*}
\M{H}_{\m{D}}&=\M{Q}_{\m{D}}\M{M}_{\m{D}}^{-1}=\mathrm{diag}\{\frac{1}{\tau}\M{I},2\tau\M{I}\}\succ0,\\
\widehat{\M{H}}_{\mathrm{D}}&=\M{Q}_{\m{D}}\tr+\M{Q}_{\m{D}}-\M{M}_{\m{D}}\tr\M{H}_{\m{D}}\M{M}_{\m{D}}-\frac{1}{2}\mathrm{diag}\{L_F\M{I},\M{0}\}=\widehat{\M{H}}_{\mathrm{N}},
\end{align*}
to be diagonal. Therefore, similar as NIDS, the stepsize condition of D-iPGM is $0<\tau<{2}/{L_F}$. We compare it with PG-EXTRA and NIDS in Table \ref{Comparison:2}. In comparison, we observe that the primal sequence $\{\M{x}^k\}$ produced by D-iPGM differs from PG-EXTRA and NIDS. Therefore, we conclude that D-iPGM is a unique algorithm that distinguishes itself from PG-EXTRA and NIDS. To accomplish a more lenient convergence condition (or network-independent stepsize condition) that differs from NIDS, D-iPGM adopts an extra primal correction step. This primal correction step can be executed independently by each agent without exchanging any private information. Additionally, the extra computational cost is almost negligible when compared to PG-EXTRA. It is worth noting that in the absence of the proximal term, the sequences ${\M{x}^k}$ generated by NIDS and D-iPGM are identical. Furthermore, as shown in Section \ref{Scenario1}, D-iPGM ($\varepsilon_k=0$) and NIDS exhibit similar performance in solving composite optimization problems.

\section{Convergence Analysis}\label{Sec4}
This section analyzes the convergence and the convergence rate of the D-iPGM. The convergence analysis will be conducted in the variational inequality context \cite{He2012}, and the forthcoming analysis of the proposed D-iPGM is based on the following fundamental lemma.
\begin{lem}[\cite{He2012}]\label{LEMF}
Let $\theta_1(\m{x})$ and $\theta_2(\m{x})$ be proper closed convex functions. If $\theta_1$ is differentiable, and the solution set of the problem $\min\{\theta_1(\m{x})+\theta_2(\m{x}):\m{x}\in\mathbb{R}^n\}$ is nonempty, then it holds that $\m{x}^*\in \arg\min\{\theta_1(\m{x})+\theta_2(\m{x}):\m{x}\in \mathbb{R}^n\}$ if and only if
$
\theta_2(\m{x})-\theta_2(\m{x}^*)+\langle \m{x}-\m{x}^*,\nabla \theta_1(\m{x}^*)\rangle\geq 0, \forall \m{x}\in \mathbb{R}^n.
$
\end{lem}

Recall the Lagrangian of problem \eqref{Problem2}:
$L(\mathbf{x},\bm{\alpha})=F(\mathbf{x})+G(\mathbf{x})-\langle\bm{\alpha},\M{V}\mathbf{x}\rangle,$
where $\M{V}=\sqrt{\M{I-W}}$ and $\bm{\alpha}\in \mathbb{R}^{Nn}$ is the Lagrange multiplier. By strongly duality, to solve problem \eqref{Problem2}, we can
consider the saddle-point problem
\begin{align}\label{SD}
\min_{\mathbf{x}}\max_{\bm{\alpha}}F(\mathbf{x})+G(\mathbf{x})-\langle\bm{\alpha},\M{V}\mathbf{x}\rangle.
\end{align}
Define $\mathcal{M}^*$ as the set of saddle-points to the above saddle-point problem. It holds that $\mathcal{M}^*=\mathcal{X}^*\times \mathcal{Y}^*$, where $\mathcal{X}^*$ is the optimal solution set to problem \eqref{Problem2} and $\mathcal{Y}^*$ is the optimal solution set to its dual problem. Denote the KKT mapping as
$$
\mathcal{J}(\mathbf{u})=\left(\begin{array}{c}
                                             \partial G(\mathbf{x})+\nabla F(\mathbf{x})-\M{V}\bm{\alpha} \\
                                            \M{V}\mathbf{x}
                                          \end{array}
\right).
$$
We have $(\mathbf{x}^{*}, \bm{\alpha}^{*}) \in \mathcal{X}^* \times \mathcal{Y}^*$ if and only if $\mathbf{0} \in \mathcal{J}(\mathbf{\mathbf{x}^*,\bm{\alpha}^*})$. A point $\mathbf{u}^*=(\mathbf{x}^*,\bm{\alpha}^*)\in \mathcal{M}^* \subseteq \mathbb{R}^{Nn}\times \mathbb{R}^{Nn}:=\mathcal{M}$ is called a saddle point of the Lagrangian function $L(\mathbf{x},\bm{\alpha})$ if
$
L(\mathbf{x}^{*}, \bm{\alpha}) \leq L(\mathbf{x}^{*}, \bm{\alpha}^{*}) \leq L(\mathbf{x}, \bm{\alpha}^{*}), \forall (\mathbf{x},\bm{\alpha}) \in \mathcal{M},
$
which can be alternatively rewritten as the following variational inequalities (by Lemma \ref{LEMF})
\begin{align}
&G(\mathbf{x})-G(\mathbf{x}^{*}) +\langle\mathbf{u}-\mathbf{u}^*,\mathcal{K}_1(\mathbf{u}^*)\rangle\geq 0, \forall \mathbf{u}\in \mathcal{M},\label{VI1}\\
&\Longleftrightarrow \nonumber\\
&\phi(\mathbf{x})-\phi(\mathbf{x}^{*}) +\langle\mathbf{u}-\mathbf{u}^*,\mathcal{K}_2(\mathbf{u}^*)\rangle\geq 0, \forall \mathbf{u}\in \mathcal{M},\label{VI2}
\end{align}
where $\phi(\M{x})=G(\M{x})+F(\M{x})$, $\mathbf{u}=\mathrm{col}\{\mathbf{x},\bm{\alpha}\}$,
$$
\mathcal{K}_1(\mathbf{u})=\left(\begin{array}{c}
                   \nabla F(\M{x})-\M{V}\bm{\alpha} \\
                   \M{V}\mathbf{x}
                 \end{array}
\right),\text{ and }
 \mathcal{K}_2(\mathbf{u})=\left(\begin{array}{c}
                   -\M{V}\bm{\alpha} \\
                   \M{V}\mathbf{x}
                 \end{array}
\right).
$$
Note that the mapping $\mathcal{K}_2$ is affine with a skew-symmetric matrix and thus we have
\begin{align}\label{skew-symmetric}
\langle \M{u}_1-\M{u}_2,\mathcal{K}_1(\M{u}_1)-\mathcal{K}_1(\M{u}_2) \rangle\equiv0,\forall \M{u}_1,\M{u}_2\in \mathcal{M}.
\end{align}
From the above analysis, it holds that the set $\mathcal{M}^*$ is also the solution set of \eqref{VI1} and \eqref{VI2}. Therefore, the convergence of D-iPGM can be proved if the generated sequence can be expressed by a similar inequality to \eqref{VI1} or \eqref{VI2} with an extra term converging to zero. Based on such observation, we will show the convergence properties of the sequence generated by D-iPGM in Sections \ref{SEC3-A}, \ref{SEC3-B}, and \ref{SEC3-C}.

\subsection{Convergence Analysis Under General Convexity}\label{SEC3-A}
To simplify the notation for the convergence analysis, we first define two matrices as the following
$$
\M{Q}=\left(
        \begin{array}{cc}
          \Gamma^{-1} & \M{V} \\
           \mathbf{0} & \frac{1}{\beta}\M{I} \\
        \end{array}
      \right),~\M{M}=\left(
                       \begin{array}{cc}
                         \M{I} & \Gamma \M{V} \\
                         \M{0} & \M{I} \\
                       \end{array}
                     \right).
$$
Then, with the matrices $\M{Q}$ and $\M{M}$, we define the following two matrices which are important in the convergence and convergence rate establishment of D-iPGM.
\begin{align}\label{MAXCON}
\M{H}:=\M{QM}^{-1},~\M{G}:=\M{Q}\tr+\M{Q}-\M{M}\tr\M{H}\M{M}.
\end{align}
\begin{prop}
If $\tau_i>0,i\in\mathcal{V}$ and $0<\beta\leq0.5/\max_i \tau_i$, the matrices $\M{H}$ and $\M{G}$ are positive definite.
\end{prop}
\begin{IEEEproof}
It follows from the definitions of $\M{Q}$ and $\M{M}$ that
\begin{align*}
\M{H}=\left(
        \begin{array}{cc}
          \Gamma^{-1} & \M{V} \\
           \mathbf{0} & \frac{1}{\beta}\M{I} \\
        \end{array}
      \right)\left(
                       \begin{array}{cc}
                         \M{I} & -\Gamma \M{V} \\
                         \M{0} & \M{I} \\
                       \end{array}
                     \right)=\left(
        \begin{array}{cc}
          \Gamma^{-1} & \M{0} \\
           \mathbf{0} & \frac{1}{\beta}\M{I} \\
        \end{array}
      \right).
\end{align*}
In addition, it holds that
\begin{align*}
\M{G}&=\M{Q}\tr+\M{Q}-\M{M}\tr\M{H}\M{M}=\M{Q}\tr+\M{Q}-\M{M}\tr\M{Q}\\
&=\left(
        \begin{array}{cc}
          2\Gamma^{-1} & \M{V} \\
           \mathbf{V} & \frac{2}{\beta}\M{I} \\
        \end{array}
      \right)-\left(
                       \begin{array}{cc}
                         \M{I} & \M{0} \\
                         \M{V}\Gamma & \M{I} \\
                       \end{array}
                     \right)\left(
        \begin{array}{cc}
          \Gamma^{-1} & \M{V} \\
           \mathbf{0} & \frac{1}{\beta}\M{I} \\
        \end{array}
      \right)\\
&=\left(
        \begin{array}{cc}
          2\Gamma^{-1} & \M{V} \\
           \mathbf{V} & \frac{2}{\beta}\M{I} \\
        \end{array}
      \right)-\left(
        \begin{array}{cc}
          \Gamma^{-1} & \M{V} \\
           \mathbf{V} & \frac{1}{\beta}\M{I}+\M{V}\Gamma\M{V} \\
        \end{array}
      \right)\\
&=\left(
        \begin{array}{cc}
          \Gamma^{-1} & \M{0} \\
           \mathbf{0} & \frac{1}{\beta}\M{I}-\M{V}\Gamma\M{V} \\
        \end{array}
      \right).
\end{align*}
Since $\M{V}=\sqrt{\M{I-W}}$ and $\sigma_M(\M{I-W})\in[0,2)$, $0<\beta\leq0.5/\max_i \tau_i$ implies that $\frac{1}{\beta}\M{I}-\M{V}\Gamma\M{V}\succ\M{0}$. This completes the proof of this proposition.
\end{IEEEproof}
With the matrix $\M{G}$, we define two more matrices as
$$
\widehat{\M{H}}_1=\M{G}-\frac{1}{2}\mathrm{diag}\{\M{L}_F,\M{0}\},\text{ and }\widehat{\M{H}}_2=\M{G}-\mathrm{diag}\{\M{L}_F,\M{0}\}.
$$
The positive definiteness of $\widehat{\M{H}}_1$ is necessary for establishing the convergence and the non-ergodic convergence rate of D-iPGM, whereas the positive definiteness of $\widehat{\M{H}}_2$ is for establishing ergodic convergence rate for D-iPGM. If $\widehat{\M{H}}_1\succ0$, we have for any $\M{u}_1,\M{u}_2\in\mathcal{M}$
\begin{align*}
c_2\|\M{u}_1-\M{u}_2\|&\leq \|\M{u}_1-\M{u}_2\|_{\widehat{\M{H}}_1}\\
&\leq\|\M{u}_1-\M{u}_2\|_{\M{H}}\leq c_1\|\M{u}_1-\M{u}_2\|,
\end{align*}
where $c_1=\sigma_M^{1/2}(\M{H}) $ and
$c_2=\sigma_m^{1/2}(\widehat{\M{H}}_1) $. Let $\tilde{\mathbf{u}}^k=\mathrm{col}\{\tilde{\mathbf{x}}^k,\tilde{\bm{\alpha}}^k\}$. We can rewrite Algorithm \ref{alg1} as
\begin{subequations}\label{IDALM}
\begin{align}
\tilde{\mathbf{x}}^{k}&=\mathrm{prox}_G^{\Gamma^{-1}}(\mathbf{x}^k-\Gamma(\nabla F(\mathbf{x}^k)-\M{V}\bm{\alpha}^k-\M{d}^k)),\label{DALM2.1}\\
\tilde{\bm{\alpha}}^{k}&=\bm{\alpha}^k-\beta \M{V}\tilde{\mathbf{x}}^{k},\label{DALM2.2}\\
\M{u}^{k+1}&=\M{u}^k-\M{M}(\M{u}^k-\tilde{\M{u}}^k) \label{DALM2.3},
\end{align}
\end{subequations}
where $\M{d}^k=\mathrm{col}\{\m{d}_1^k,\cdots,\m{d}_N^k\}$.
Note that if $\bm{\alpha}^0=\bm{\lambda}^0=\mathbf{0}$, the sequence $\{\mathbf{x}^k\}$ generated by Algorithm \ref{alg1} and recursion \eqref{IDALM} are identical. In the analysis, we consider the sequence $\{(\mathbf{x}^k,\bm{\alpha}^k)\}$ to illustrate the convergence of Algorithm \ref{alg1}.

To establish the global convergence and derive the convergence rate of D-iPGM, we first give the following lemma.
\begin{lem}\label{Lemma:Prediction1}
Suppose that Assumptions \ref{ass2} and \ref{ass1} hold. For any $\forall (\M{x},\bm{\alpha})\in\mathcal{M}$, the sequence $\{\tilde{\mathbf{u}}^k\}$ generated by \eqref{IDALM} satisfies
\begin{align}\label{Lemma:P1}
&G(\mathbf{x})-G(\tilde{\mathbf{x}}^{k})+\langle\mathbf{u}-\tilde{\mathbf{u}}^{k},\mathcal{K}_1(\mathbf{u})\rangle \geq \langle\mathbf{u}-\tilde{\mathbf{u}}^k,\M{Q}(\mathbf{u}^k-\tilde{\mathbf{u}}^k)\rangle\nonumber\\
&\quad+\langle\mathbf{x}-\tilde{\mathbf{x}}^k,\mathbf{d}^k\rangle-\frac{1}{4}\|\mathbf{x}^k-\tilde{\mathbf{x}}^k\|_{\mathbf{L}_F}^2,
\end{align}
and
\begin{align}\label{Lemma:P2}
&\phi(\mathbf{x})-\phi(\tilde{\mathbf{x}}^{k})+\langle\mathbf{u}-\tilde{\mathbf{u}}^{k},\mathcal{K}_2(\mathbf{u})\rangle \geq \langle\mathbf{u}-\tilde{\mathbf{u}}^k,\M{Q}(\mathbf{u}^k-\tilde{\mathbf{u}}^k)\rangle\nonumber\\
&\quad+\langle\mathbf{x}-\tilde{\mathbf{x}}^k,\mathbf{d}^k\rangle-\frac{1}{2}\|\mathbf{x}^k-\tilde{\mathbf{x}}^k\|_{\mathbf{L}_F}^2.
\end{align}
\end{lem}
\begin{IEEEproof}
See Appendix \ref{Appendix:Lemma2}.
\end{IEEEproof}
To ensure the positive definiteness of $\widehat{\M{H}}_1$, we assume that the stepsize satisfies $0<\tau_i<{2}/{L_i},i\in\mathcal{V}$, $0<\beta\leq0.5/\max_i \tau_i$. Then, we give the following theorem to show that the sequence $\{(\mathbf{x}^k,\bm{\alpha}^k)\}$ generated by \eqref{IDALM} converges to a primal-dual optimal solution of problem \eqref{SD}.
\begin{thm}\label{TH1}
Suppose that Assumptions \ref{ass2} and \ref{ass1} hold. If $0<\tau_i<{2}/{L_i},i\in\mathcal{V}$, $0<\beta\leq0.5/\max_i \tau_i$, and $\sum_{k=1}^{\infty}\varepsilon_k<\infty$, the sequence $\{\mathbf{u}^k\}$ generated by \eqref{IDALM} satisfies,
\begin{enumerate}
  \item for any $\mathbf{u}^*\in \mathcal{M}^*$, it holds that
\begin{align}\label{Theorem1:P1}
&\|\mathbf{u}^{k+1}-\mathbf{u}^*\|^2_{\M{H}}\leq \|\mathbf{u}^k-\mathbf{u}^*\|^2_{\M{H}}-\|\mathbf{u}^k-\tilde{\mathbf{u}}^{k}\|^2_{\widehat{\M{H}}_1}\nonumber\\
&\quad+2\langle\tilde{\mathbf{x}}^k-\mathbf{x}^*,\mathbf{d}^k\rangle;
\end{align}
  \item the sequence $\{\mathbf{u}^k\}$ is bounded;
  \item these exists $\mathbf{u}^{\infty}\in \mathcal{M}^*$ such that $\lim_{k\rightarrow \infty}\mathbf{u}^k=\mathbf{u}^{\infty}$, i.e., $\M{x}^{\infty}=\mathbf{1}_N\otimes \M{x}^{\infty}$, and $\M{x}^{\infty}$ solves problem \eqref{Problem}.
\end{enumerate}
\end{thm}
\begin{IEEEproof}
See Appendix \ref{Appendix:Theorem1}.
\end{IEEEproof}
Although \eqref{Theorem1:P1} may not guarantee the monotonicity of sequence $\{\|\mathbf{u}^k-\mathbf{u}^{\infty}\|^2_{\M{H}}\}$, since $\{\M{u}^k\}$ is bounded and $\|\mathbf{d}^k\|$ is summable, the sequence $\|\mathbf{u}^k-\mathbf{u}^{\infty}\|^2_{\M{H}}$ is a quasi-Fej\'{e}r monotone sequence in $\M{H}$-norm, i.e.,
\begin{align*}
\|\M{u}^{k+1}-\M{u}^{\infty}\|_{\M{H}}\leq \|\mathbf{u}^k-\mathbf{u}^{\infty}\|_{\M{H}}+\bar{\mu}\varepsilon_k,
\end{align*}
which guarantees the convergence of the Di-PGM.

When the subproblem \eqref{Subx} is solved exactly or the proximal operator of $g_i$ has an analytical solution, we have $\|\M{d}^k\|=0$. In this case, the sequence $\{\mathbf{u}^k\}$ generated by \eqref{IDALM} satisfies
\begin{align*}
\|\mathbf{u}^{k+1}-\mathbf{u}^*\|^2_{\M{H}}\leq \|\mathbf{u}^k-\mathbf{u}^*\|^2_{\M{H}}-\|\mathbf{u}^k-\tilde{\mathbf{u}}^{k}\|^2_{\widehat{\M{H}}_1}, \forall \M{u}^*\in\mathcal{M}^*,
\end{align*}
which implies that $\{\M{u}^k\}$ is a Fej\'{e}r monotone sequence with respect to $\mathcal{M}^*$ in $\M{H}$-norm.

\subsection{Sublinear Convergence Rate Under General Convexity}\label{SEC3-B}
With general convexity, this subsection established the sublinear convergence rate of D-iPGM.

We give the following theorem to show the $O(1/k)$ non-ergodic and ergodic convergence rate of D-iPGM.
\begin{thm}\label{TH2}
Suppose that Assumptions \ref{ass2} and \ref{ass1} hold. If $0<\tau_i<{2}/{L_i},i\in\mathcal{V}$, $0<\beta\leq0.5/\max_i \tau_i$, and $\sum_{k=1}^{\infty}\varepsilon_k<\infty$, it holds that
\begin{enumerate}
  \item Running-average first-order optimality residual:
  \begin{align*}
  \frac{1}{K+1}\sum_{k=0}^{K}\mathrm{dist}^2(\mathbf{0},\mathcal{J}(\tilde{\M{u}}^{k}))=O\left(\frac{1}{K+1}\right).
  \end{align*}
  \item Running-best first-order optimality residual:
  \begin{align*}
  \min_{0\leq k\leq K} \{\mathrm{dist}^2(\mathbf{0},\mathcal{J}(\tilde{\M{u}}^{k}))\} =o\left(\frac{1}{K+1}\right).
  \end{align*}
\end{enumerate}
Let $\M{X}^K=\frac{1}{K}\sum_{k=0}^{K-1}\tilde{\M{x}}^k$ and $\M{\Lambda}^K=\frac{1}{K}\sum_{k=0}^{K-1}\tilde{\bm{\alpha}}^k$. Furthermore, if $0<\tau_i<{1}/{L_i},i\in\mathcal{V}$, it holds that
\begin{enumerate}
  \item Primal-dual gap: $L(\M{X}^K,\bm{\alpha})-L(\M{x},\M{\Lambda}^K)=O\left(\frac{1}{K}\right)$.
  \item Primal optimality gap: $|\phi(\M{X}^K)-\phi(\M{x}^*)|=O\left(\frac{1}{K}\right)$.
  \item Consensus error: $\|\sqrt{\M{I-W}}\M{X}^K\|=O\left(\frac{1}{K}\right)$.
\end{enumerate}
\end{thm}
\begin{IEEEproof}
See Appendix \ref{Appendix:Theorem2}.
\end{IEEEproof}
With the stepsize contion $\tau_i<2/L_i,i\in\mathcal{V}$, we establish the non-ergodic iteration complexity of D-iPGM. With a stronger convergence condition $\tau_i<1/L_i,i\in\mathcal{V}$, the ergodic $O(1/k)$ convergence rate of the primal-dual gap, the primal optimality gap, and the consensus error are proved.

The next theorem gives the $o(1/k)$ non-ergodic convergence rates of D-iPGM when $\varepsilon_k=0,k\geq0$.
\begin{thm}\label{Theorem:point-wise-rate}
Suppose that Assumptions \ref{ass2} and \ref{ass1} hold. If $0<\tau_i<{2}/{L_i},i\in\mathcal{V}$, $0<\beta\leq0.5/\max_i \tau_i$, and $\varepsilon_k=0,\forall k\geq0$, for any integer $k\geq0$, it holds that
\begin{align}\label{eq:point-wise-rate1}
\|\M{M}(\M{u}^{k+1}-\tilde{\M{u}}^{k+1})\|^2_{\M{H}}\leq\|\M{M}(\M{u}^{k}-\tilde{\M{u}}^{k})\|^2_{\M{H}}.
\end{align}
Moreover, for any $\M{u}^*\in\mathcal{M}^*$, it holds that
\begin{enumerate}
  \item Successive difference:
  \begin{align*}
  \|\M{u}^{k+1}-\M{u}^k\|^2_{\M{H}}=o\left(\frac{1}{k+1}\right).
  \end{align*}
  \item First-order optimality residual:
  \begin{align*}
  \mathrm{dist}^2(\M{0},\mathcal{J}(\tilde{\M{u}}^k))=o\left(\frac{1}{k+1}\right).
  \end{align*}
\end{enumerate}
\end{thm}
\begin{IEEEproof}
See Appendix \ref{Appendix:Theorem:point-wise-rate}.
\end{IEEEproof}
It follows from \eqref{Proof:lemma1:3}, i.e.,
\begin{align*}
&G(\mathbf{x})-G(\tilde{\mathbf{x}}^{k})+\langle\mathbf{x}-\tilde{\mathbf{x}}^{k},\nabla F(\mathbf{x}^{k})\rangle+\langle\mathbf{u}-\tilde{\mathbf{u}}^{k},\mathcal{K}_2(\tilde{\mathbf{u}}^k)\rangle\nonumber\\
&\geq\langle\mathbf{u}-\tilde{\mathbf{u}}^{k}, \M{Q}(\mathbf{u}^{k}-\tilde{\mathbf{u}}^{k})\rangle+\langle\mathbf{x}-\tilde{\mathbf{x}}^{k},\mathbf{d}^k\rangle, \forall \mathbf{u} \in \mathcal{M},
\end{align*}
the successive difference $\|\M{u}^{k+1}-\M{u}^k\|^2_{\M{H}}$ can be used to measure the accuracy of the iterate $\mathbf{u}^{k+1}$ to $\mathcal{M}^*$, when $\varepsilon_k=0,k\geq0$. More specifically, if $\|\M{u}^{k+1}-\M{u}^k\|^2_{\M{H}}=0$, by \eqref{DALM2.3}, we have $\M{u}^{k}=\tilde{\M{u}}^k=\M{u}^{k+1}$. Hence, it deduces that
\begin{align*}
G(\mathbf{x})-G(\tilde{\mathbf{x}}^{k})+\langle\mathbf{u}-\tilde{\mathbf{u}}^{k},\mathcal{K}_1(\tilde{\mathbf{u}}^k)\rangle\geq0, \forall \mathbf{u} \in \mathcal{M}.
\end{align*}
Compared to \eqref{VI1}, $\M{u}^{k+1}=\tilde{\mathbf{u}}^k\in \mathcal{M}^*$. Therefore, $\|\M{u}^{k+1}-\M{u}^k\|^2_{\M{H}}$ can be viewed as an error measurement after $k$ iterations of the D-iPGM. The same error measurement is considered in PG-EXTRA \cite{PGEXTR} and NIDS \cite{NDIS}.

\subsection{Linear Convergence Rate Under Metric Subregularity}\label{SEC3-C}
By some variational analysis techniques, this subsection proves the linear convergence of D-iPGM under metric subregularity. Recall the notion of metric subregularity \cite{va2004}, which is important in the establishment of the local linear convergence rate \cite{va2020,va2021}.
\begin{defn}
A set-valued mapping $\Psi:\mathbb{R}^{n}\rightrightarrows \mathbb{R}^{m}$ is metrically subregular at $(\bar{u},\bar{v})\in \mathrm{gph}(\Psi)$ if for some $\epsilon>0$ there exists $\kappa\geq0$ such that
$$
\mathrm{dist}(u,\Psi^{-1}(\bar{v}))\leq \kappa~ \mathrm{dist}(\bar{v},\Psi(u)),~ \forall u\in \mathcal{B}_{\epsilon}(\bar{u}),
$$
where $\mathrm{gph}(\Psi):=\{(u,v):v=\Phi(u)\}$, $\Psi^{-1}(v):=\{u\in\mathbb{R}^n:(u,v)\in\mathrm{gph}(\Psi)\}$.
\end{defn}
To ensure the linear rate convergence of D-iPGM, we give the following assumption.
\begin{ass}\label{ASS3}
There exists an integer $\bar{k}>0$ and a sequence $\{\varrho_k\}$ such that
\begin{equation}\label{Cr1}
\sup_{k\geq \bar{k}}\{\varrho_k\}<{1}/{\bar{\mu}},\text{ and }\|\mathbf{d}^k\|\leq \varrho_k\|\mathbf{x}^{k}-\mathbf{x}^{k+1}\|, \forall k\geq \bar{k},
\end{equation}
where $\bar{\mu}=\sigma_M(\M{H}^{\frac{1}{2}}\M{M})\sigma_M(\Gamma)(\beta\sigma_M(\sqrt{\M{I-W}})+1)$.
\end{ass}
\begin{rmk}
From Theorem \ref{TH1}, we have $\sum_{k=0}^{\infty}\|\M{x}^{k+1}-\M{x}^k\|<\infty$. It implies that $\lim_{k\rightarrow\infty}(\|\M{x}^{k}-\M{x}^{k-1}\|-\|\M{x}^{k+1}-\M{x}^k\|)=0$. Therefore, to facilitate implementation, we can use $\|\M{x}^{k}-\M{x}^{k-1}\|$ instead of $\|\M{x}^{k+1}-\M{x}^{k}\|$ in practice.
\end{rmk}

Then, the linear convergence rate of D-iPGM is presented.
\begin{thm}\label{TH3}
Suppose that Assumptions \ref{ass2}, \ref{ass1}, and \ref{ASS3} hold, $0<\tau_i<{2}/{L_i},i\in\mathcal{V}$, $0<\beta\leq0.5/\max_i \tau_i$, and the sequence $\{\mathbf{u}^k\}$ generated by \eqref{IDALM} converges to $\mathbf{u}^{\infty}$. If the KKT mapping
 $\mathcal{J}(\mathbf{u})$
is metrically subregular at $(\mathbf{u}^{\infty},\M{0})$, then there exist an integer $\bar{k}>0$ and a constant $\kappa>0$ such that for all $k\geq\bar{k}$
\begin{align}\label{linear-rate}
\mathrm{dist}_{\M{H}}(\mathbf{u}^{k+1},\mathcal{M}^*)\leq\vartheta_k~ \mathrm{dist}_{\M{H}}(\mathbf{u}^{k},\mathcal{M}^*),
\end{align}
with
\begin{align*}
\vartheta_k&=\frac{{\mu}(1+\omega)\varrho_k+\bar{\vartheta}}{(1-{\mu}\varrho_k)}, ~\mu={\bar{\mu}}/{c_2},~\omega=\sqrt{c_1/c_2},\\
\bar{\vartheta}^2&={1-\frac{c_2^2}{c^2_1(\kappa\kappa_2c_1+1)^2}}<1,\\
\kappa^2_2&={\frac{\max\{3\sigma_M(\M{V}^2)+\frac{1}{\beta^2},3\max_i\{L^2_i\}+3\sigma_m(\Gamma^{-2})\}}{\sigma_m(\widehat{\M{H}}_1)}}.
\end{align*}
Moreover, if
\begin{equation}\label{Cr2}
\sup_{k\geq\bar{k}}\{\varrho_k\}<\frac{1-\bar{\vartheta}}{\mu(2+\omega)},
\end{equation}
then one has $\vartheta_k<1$, and thus the convergence rate of $\mathrm{dist}_{\M{H}}(\mathbf{u}^k,\mathcal{M}^*)$ is Q-linear.
\end{thm}
\begin{IEEEproof}
See Appendix \ref{Appendix:Theorem3}.
\end{IEEEproof}
Note that if $\{\varrho_k\}$ converges to $0$ as $k\rightarrow\infty$, condition \eqref{Cr2} hold eventually for $\bar{k}$ sufficiently large, since $\bar{\mu}$ and $\bar{\vartheta}$ are constant. In addition, if $\varrho_k=0,k\geq0$, it holds that $\vartheta_k\equiv\overline{\vartheta}$, which implies that for all $k\geq\bar{k}$
\begin{align*}
\mathrm{dist}^2_{\M{H}}(\mathbf{u}^{k+1},\mathcal{M}^*)\leq\left(1-\frac{c_2^2}{c^2_1(\kappa\kappa_2c_1+1)^2}\right) \mathrm{dist}^2_{\M{H}}(\mathbf{u}^{k},\mathcal{M}^*).
\end{align*}
Thus, the Q-linear convergence of D-iPGM holds for exact proximal iteration.

\begin{rmk}
In \cite{Sulaiman2021}, a dimension independent lower bound of problem \eqref{Problem} is given, and it is reported that if each agent owns a different local nonsmooth term, then exact linear convergence cannot be attained in the worst case ($f_i$ is a strongly convex smooth function and $g_i$ is nonsmooth convex with closed form proximal mappings). However, in Theorem \ref{TH3}, the dimension dependent linear rate of D-iPGM can be established, which does not contradict with the exiting results.
\end{rmk}

Next, using the analysis techniques in \cite{va2020} and \cite{va2021}, we give the sufficient condition such that $\mathcal{J}(\mathbf{u})$ is metrically subregular at $(\mathbf{u}^{\infty},\mathbf{0})$. Similar as \cite{va2020} and \cite{va2021}, a convex function $F(\mathbf{x})$ is said to satisfy the structured assumption if $F(\mathbf{x})=\sum_{i=i}^{N}f_i(\mathrm{x}_i)=\sum_{i=1}^{N}h_i(A_i\mathrm{x}_i)+\langle q_i,\mathrm{x}_i\rangle$,
where $A_i$ is a $m_i\times n$ matrix, $q_i$ is a vector in $\mathbb{R}^n$, and $h_i$ is smooth and essentially
locally strongly convex, i.e., for any compact and convex subset $\mathbb{K}$, $h_i$ is strongly convex on $\mathbb{K}$. Some commonly used loss functions in machine learning automatically satisfy the structured assumption. We give Table \ref{Table2} to summarize cases satisfying the structured assumption, where $\mathrm{b}^1_i\in \mathbb{R}^n$,  $\mathrm{b}^2_i\in \{0,1\}^{n}$, $\mathrm{b}^3_i\in \mathbb{R}^n_+$, $\mathrm{b}^4_i\in \{-1,1\}^{n}$, and $\mathrm{b}_i^5\in\{0,1,2,\cdots\}$.
\begin{table}[!t]
\renewcommand\arraystretch{1.5}
\begin{center}
\caption{Some Commonly Used Convex Loss Functions Satisfying The Structured Assumption.}
\scalebox{0.9}{
\begin{tabular}{ccc}
\hline
~&\textbf{Local loss function $h_i$}&\textbf{Convexity of $f_i$ }\\
\hline
Linear regression&$\frac{1}{2}\|\mathrm{x}_i-\mathrm{b}^1_i\|^2$&Convex\\
Logistic regression&$\sum_{j=1}^{n}\log(1+e^{\mathrm{x}_{ij}})-(\mathrm{b}_i^2)^{\top}\mathrm{x}_i$&Convex\\
Logistic regression&$-\sum_{j=1}^{n}\log(\mathrm{x}_{ij})+(\mathrm{b}_i^3)^{\top}\mathrm{x}_i$&Convex\\
Likelihood estimation&$\sum_{j=1}^{n}\log(1+e^{(\mathrm{b}_i^4)^{\top}\mathrm{x}_i})$&Convex\\
Poisson regression&$\sum_{j=1}^{n}e^{\mathrm{x}_{ij}}-(\mathrm{b}_i^5)^{\top}\mathrm{x}_i$&Convex\\
\hline
\end{tabular}}
\label{Table2}
\end{center}
\end{table}

By the structure of $F$, the following proposition is provided, which presents an alternative characterization of $\mathcal{M}^*$.
\begin{prop}\label{PRO1}
Suppose that $F(\mathbf{x})$ satisfies the structured assumption. The set $\mathcal{M}^*$ can be characterized as
\begin{align*}
\mathcal{M}^*=\{\mathbf{u}:&\widetilde{\mathbf{A}} \mathbf{x}=\tilde{\mathbf{t}}, \sqrt{\mathbf{I-W}}\mathbf{x}=\mathbf{0}, \\
&\mathbf{0}\in \partial G(\mathbf{x})+\tilde{\bm{\xi}}-\sqrt{\mathbf{I-W}}\bm{\alpha}
\},
\end{align*}
with some $\tilde{\mathbf{t}}$ satisfying $\widetilde{\mathbf{A}}\mathbf{x}=\tilde{\mathbf{t}}$ for all $\mathbf{x}\in \mathcal{X}^*$, and
$\tilde{\bm{\xi}}=\mathrm{col}\{A_1^{\top}\nabla h_1(\tilde{t}_1) +q_1,\cdots,A_N^{\top}\nabla h_N(\tilde{t}_N) +q_N\}$,
where $\tilde{\mathbf{t}}=\mathrm{col}\{\tilde{t}_1,\cdots,\tilde{t}_N\}$, and $\widetilde{\mathbf{A}}$ is a block diagonal matrix with its $(i,i)$-th block being $A_i$ and other blocks being $\mathbf{0}$.
\end{prop}
The proof of Proposition \ref{PRO1} is rather standard by \cite[Lemma 2.1]{luo1992} and thus is omitted here. To facilitate the analysis, an auxiliary perturbed set-valued mapping with perturbation $\mathbf{P}=(P_1,P_2,P_3)$ is introduced associated with $\mathcal{M}^*$:
$
\mathbf{\Gamma}(\mathbf{P})=\{\mathbf{u}:P_1=\widetilde{\mathbf{A}}\mathbf{x}-\tilde{\mathbf{t}},P_2=\M{V}\mathbf{x},P_3\in\partial G(\mathbf{x})+\tilde{\bm{\xi}}-\M{V}\bm{\alpha}\}.
$
It is obvious that $\Gamma(\mathbf{0})=\mathcal{M}^*$. Similar to \cite[Proposition 40]{va2020} and \cite[Proposition 4.1]{va2021}, one has the following equivalence.

\begin{prop}\label{PRO2}
Assume that structured assumption of $F$ is satisfied. The metric subregularity conditions of $\mathbf{\Gamma}^{-1}$ and $\mathcal{J}$ are equivalent. Precisely, given $\bar{\mathbf{u}}\in\mathcal{M}^*$, the following two statements are equivalent:

(i) There exist $\tilde{\kappa}_{1}, \tilde{\epsilon}_{1}>0$ such that
$$
\mathrm{dist}(\mathbf{u}, \M{\Gamma}(\mathbf{0})) \leq \tilde{\kappa}_{1}~ \mathrm{dist}(\mathbf{0}, \mathbf{\Gamma}^{-1}(\mathbf{u})),  \forall\mathbf{u} \in \mathcal{B}_{\tilde{\epsilon}_{1}}(\bar{\mathbf{u}}).
$$

(ii) There exist $\tilde{\kappa}_{2}, \tilde{\epsilon}_{2}>0$ such that
$$
\mathrm{dist}(\mathbf{u}, \mathcal{J}^{-1}(\mathbf{0})) \leq \tilde{\kappa}_{2}~ \mathrm{dist}(\mathbf{0}, \mathcal{J}(\mathbf{u})),  \forall\mathbf{u} \in \mathcal{B}_{\tilde{\epsilon}_{2}}(\bar{\mathbf{u}}).
$$
\end{prop}

Suppose that $F(\M{x})$ satisfies the structured assumption and $\partial G$ is a polyhedral multifunction, i.e., those whose graphs are unions of finite number of polyhedral convex sets. In this scenario, it holds that $\M{\Gamma}$ is a polyhedral multifunction and hence $\mathbf{\Gamma}^{-1}$ is also a polyhedral multifunction. By \cite[Proposition 1]{Robinson}, for any point $(\mathbf{u},\mathbf{P})\in\mathrm{gph}(\mathbf{\Gamma}^{-1})$, $\mathbf{\Gamma}^{-1}$ is metrically subregular at $(\mathbf{u},\mathbf{P})$. Thus, by Proposition \ref{PRO2}, $\mathcal{J}$ is metrically subregular at $(\mathbf{u}^*,\mathbf{0})$ for any given $\mathbf{u}^*\in\mathcal{M}^*$. It follows from \cite[Proposition 7]{va2021} that if $g_i$ is a polyhedral convex function, which includes the indicator function of a polyhedral set and the polyhedral convex regularizer, or $g_i$ is convex piecewise linear-quadratic function, then $\partial g_i$ is a polyhedral multifunction. This case covers scenarios where $g_i$ is the $l_{\infty}$-norm regularizer, the $l_{1}$-norm regularizer, the elastic net regularizer \cite{Hui2005}, the generalized LASSO regularizer \cite{Example3}, and the fused LASSO regularizer \cite{Tibshirani2005}.

\begin{figure*}[!t]
\centering
\setlength{\abovecaptionskip}{-2pt}
\subfigure{
\includegraphics[width=0.32\linewidth]{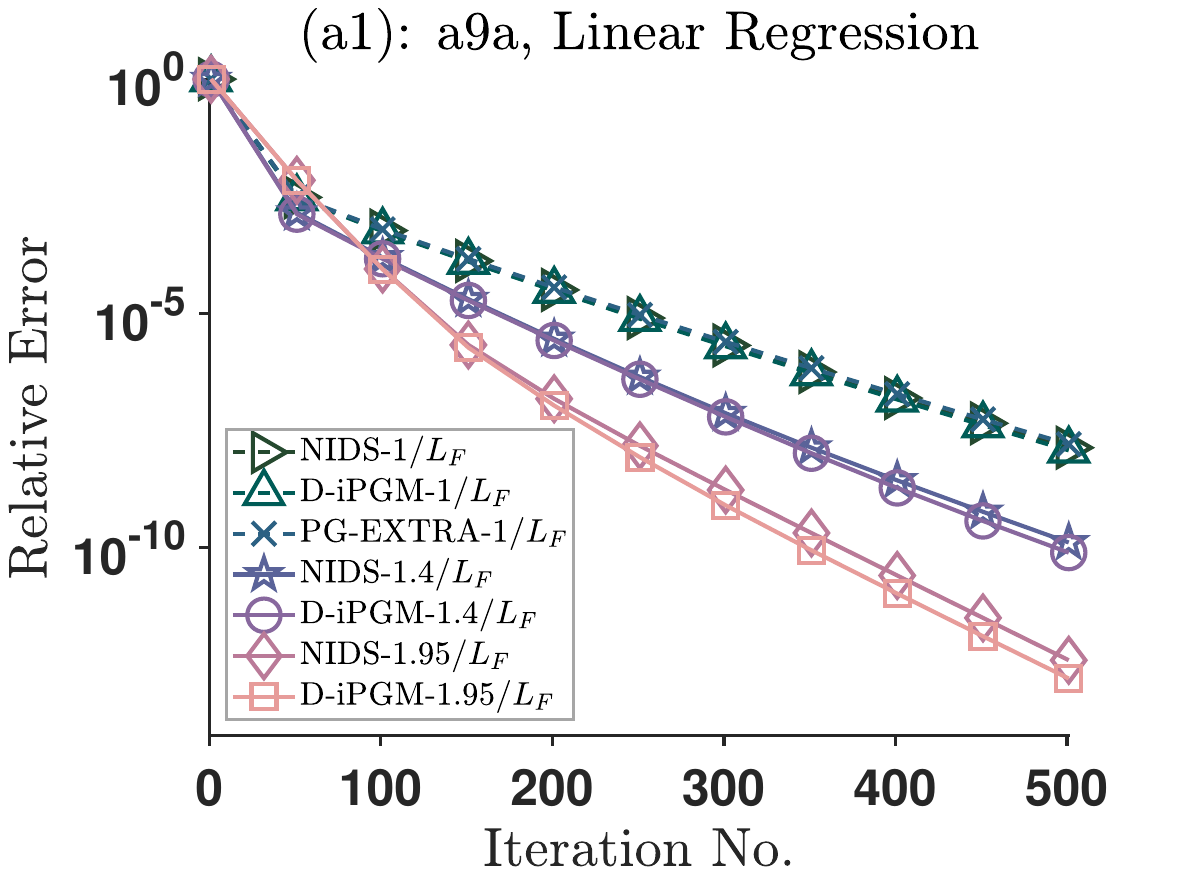}}
\subfigure{
\includegraphics[width=0.32\linewidth]{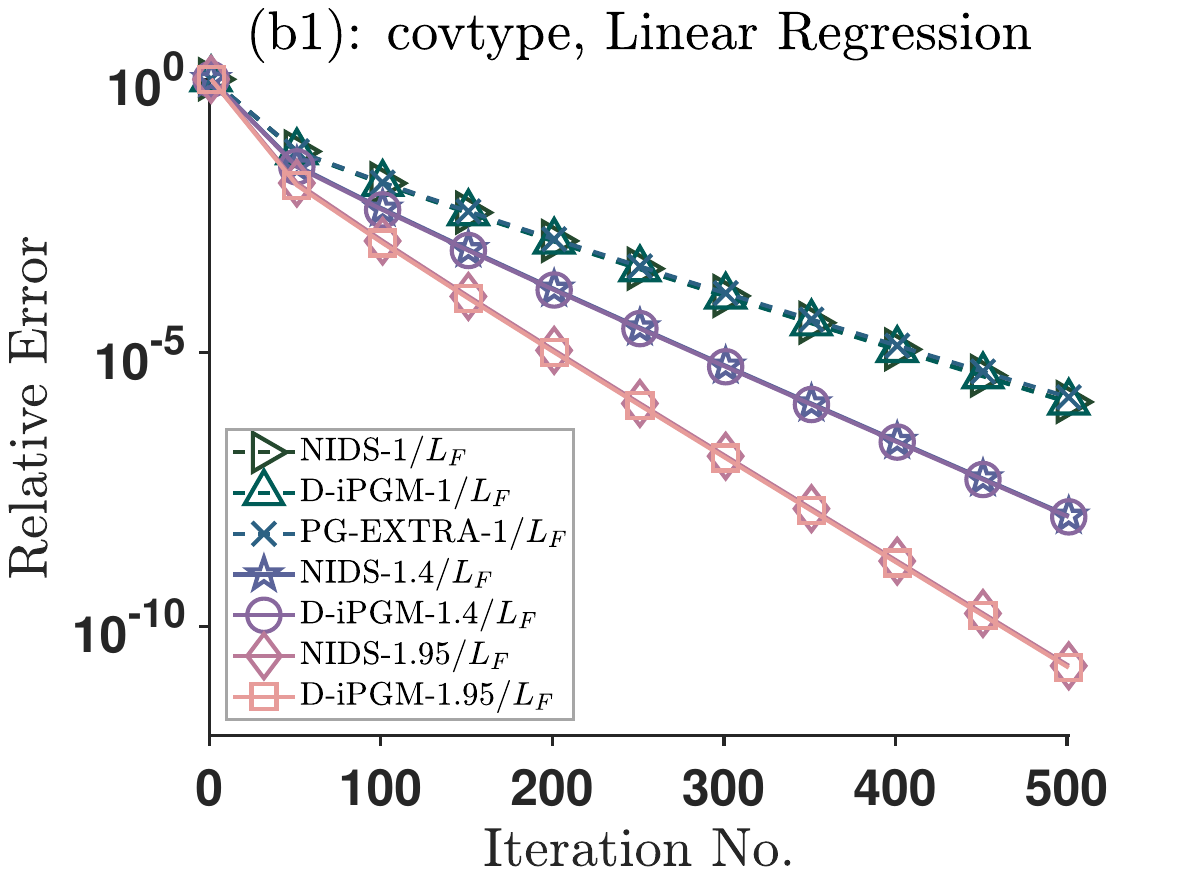}}
\subfigure{
\includegraphics[width=0.32\linewidth]{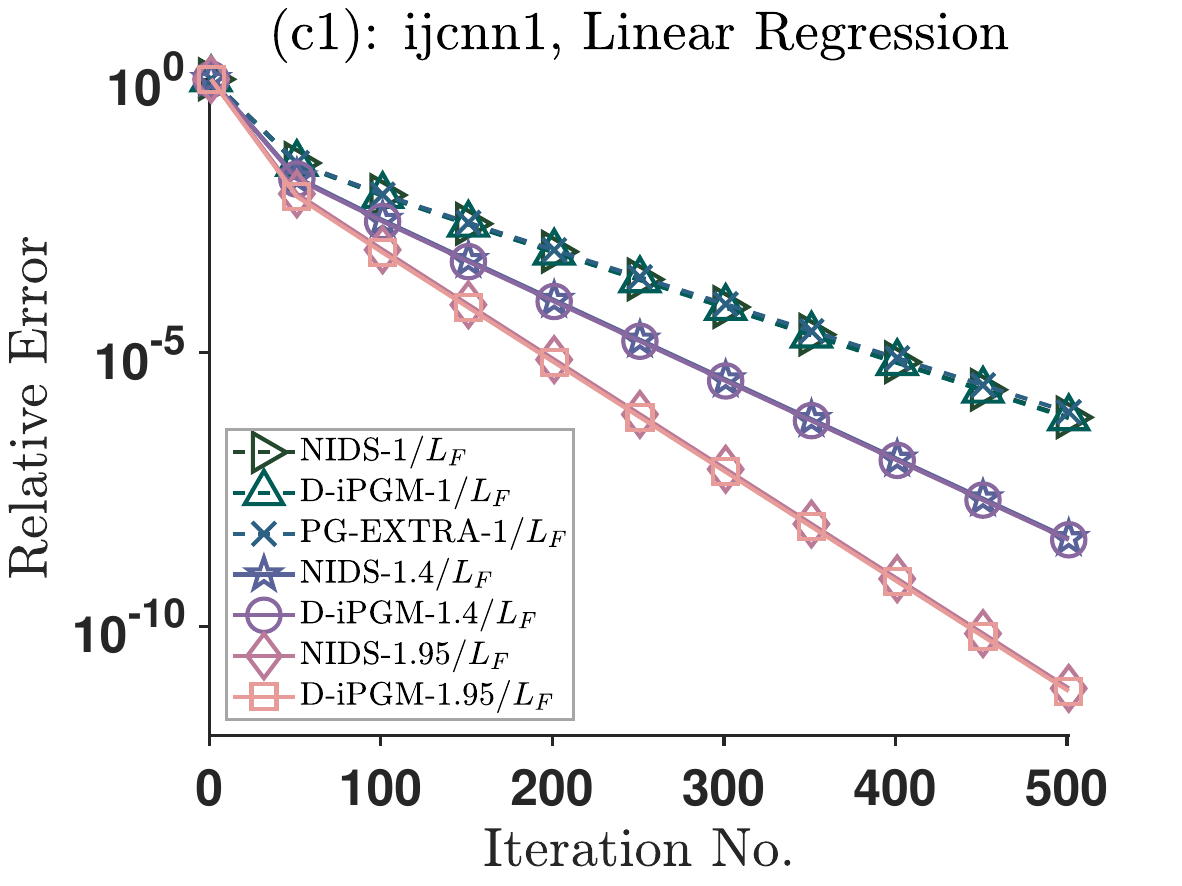}}
\subfigure{
\includegraphics[width=0.32\linewidth]{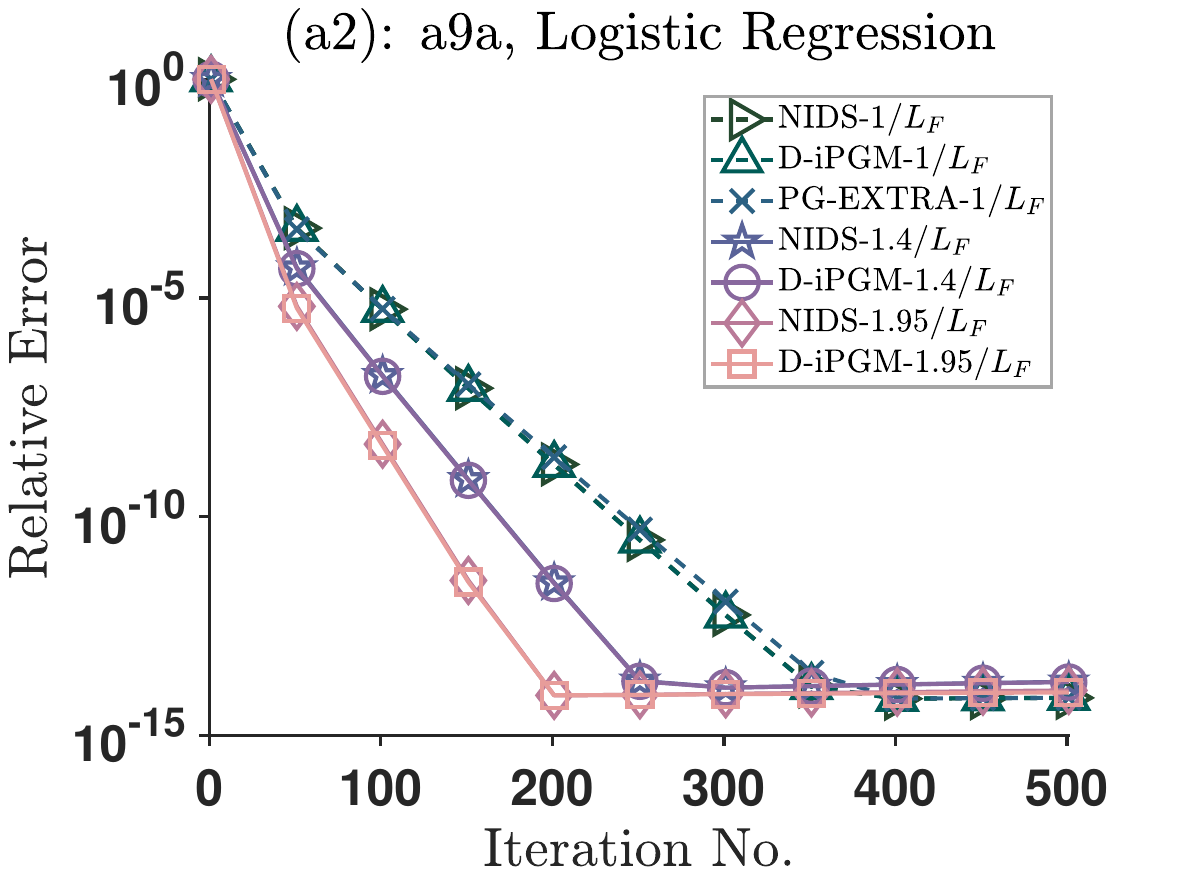}}
\subfigure{
\includegraphics[width=0.32\linewidth]{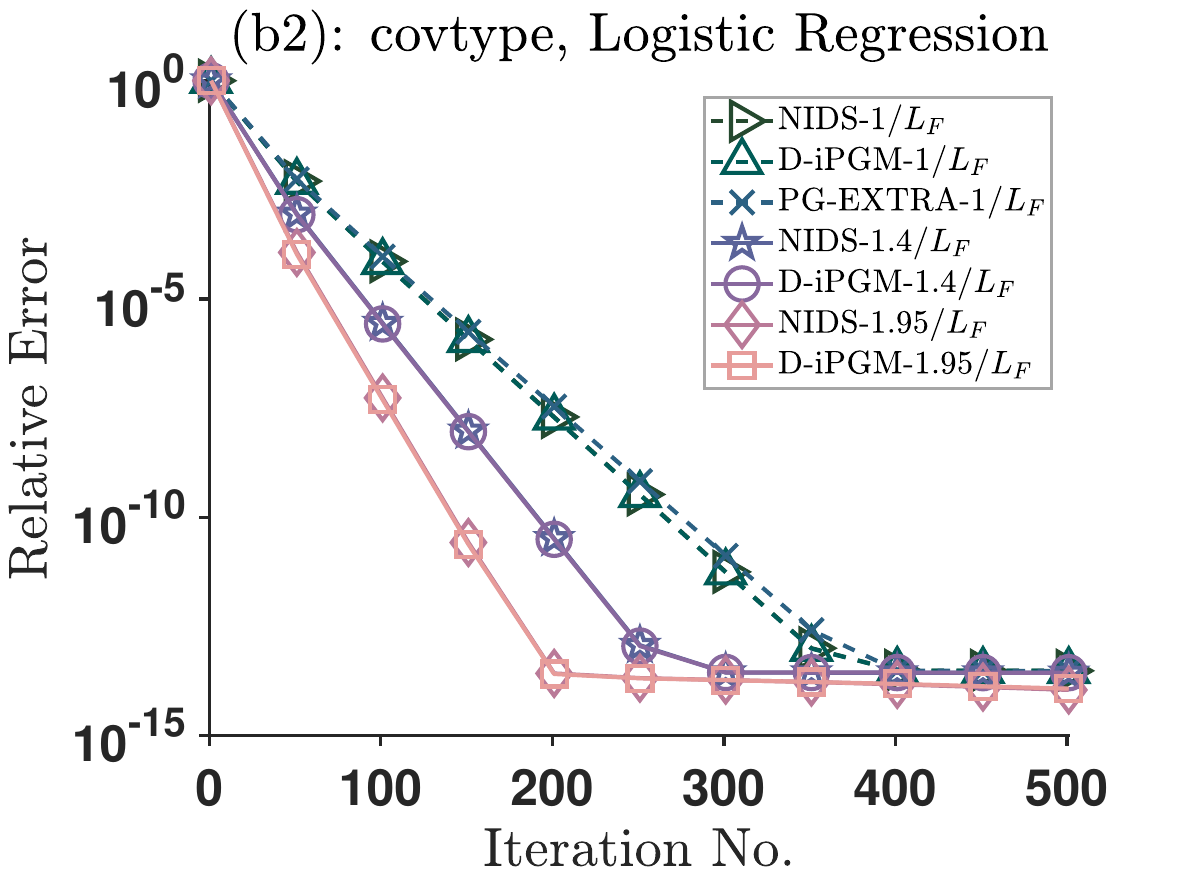}}
\subfigure{
\includegraphics[width=0.32\linewidth]{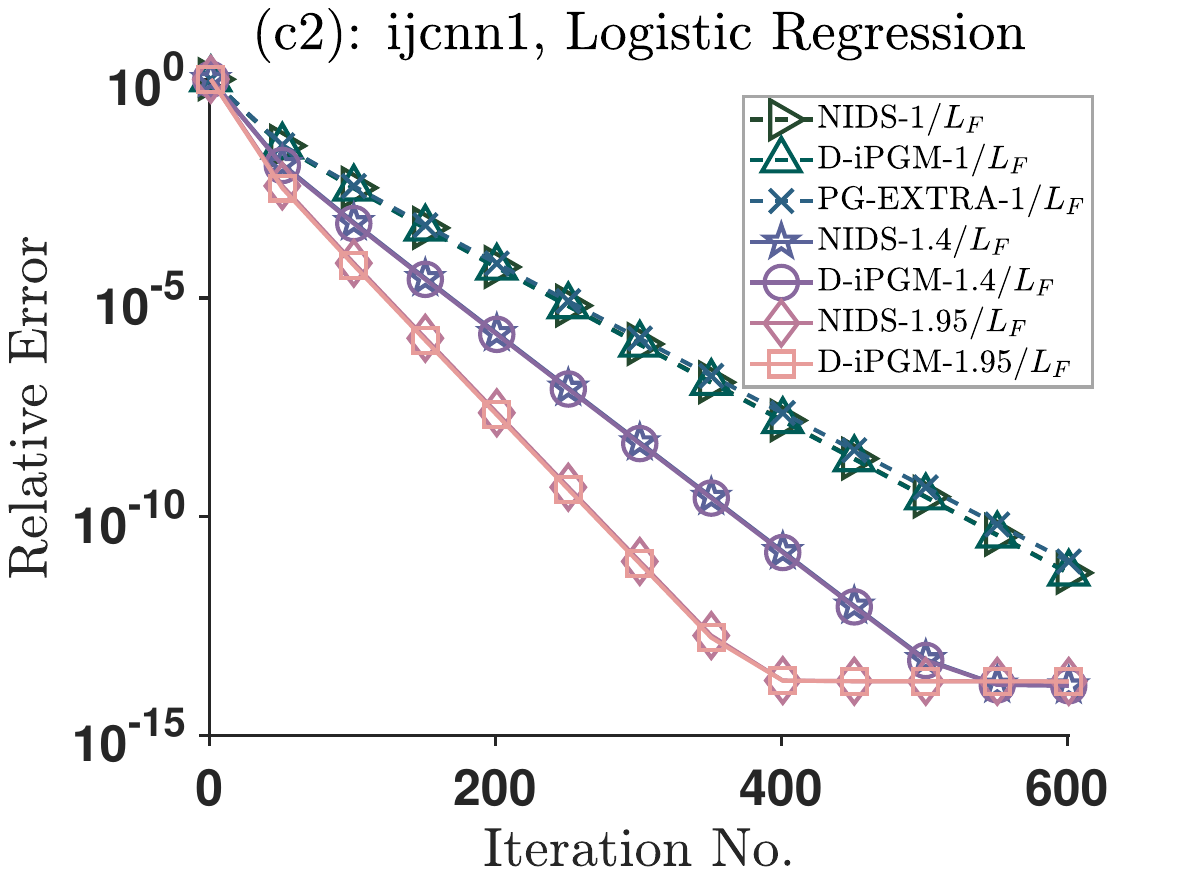}}
\caption{The relative error $\|\M{x}^k-\M{x}^*\|/\|\M{x}^*\|$ with respect to iteration numbers when the regularizer is $\nu_{1,i}\|\m{x}\|+\frac{\nu_{2,i}}{2}\|\m{x}\|^2$, where $\nu_{1,i}=0.01$ and $\nu_{2,i}=1$. Here, different stepsizes for D-iPGM, PG-EXTRA, and NIDS are considered, and the stepsize across the network of agents are the same.
}
\label{Sim:case1.1}
\end{figure*}

\begin{figure*}[!t]
\centering
\setlength{\abovecaptionskip}{-2pt}
\subfigure{
\includegraphics[width=0.32\linewidth]{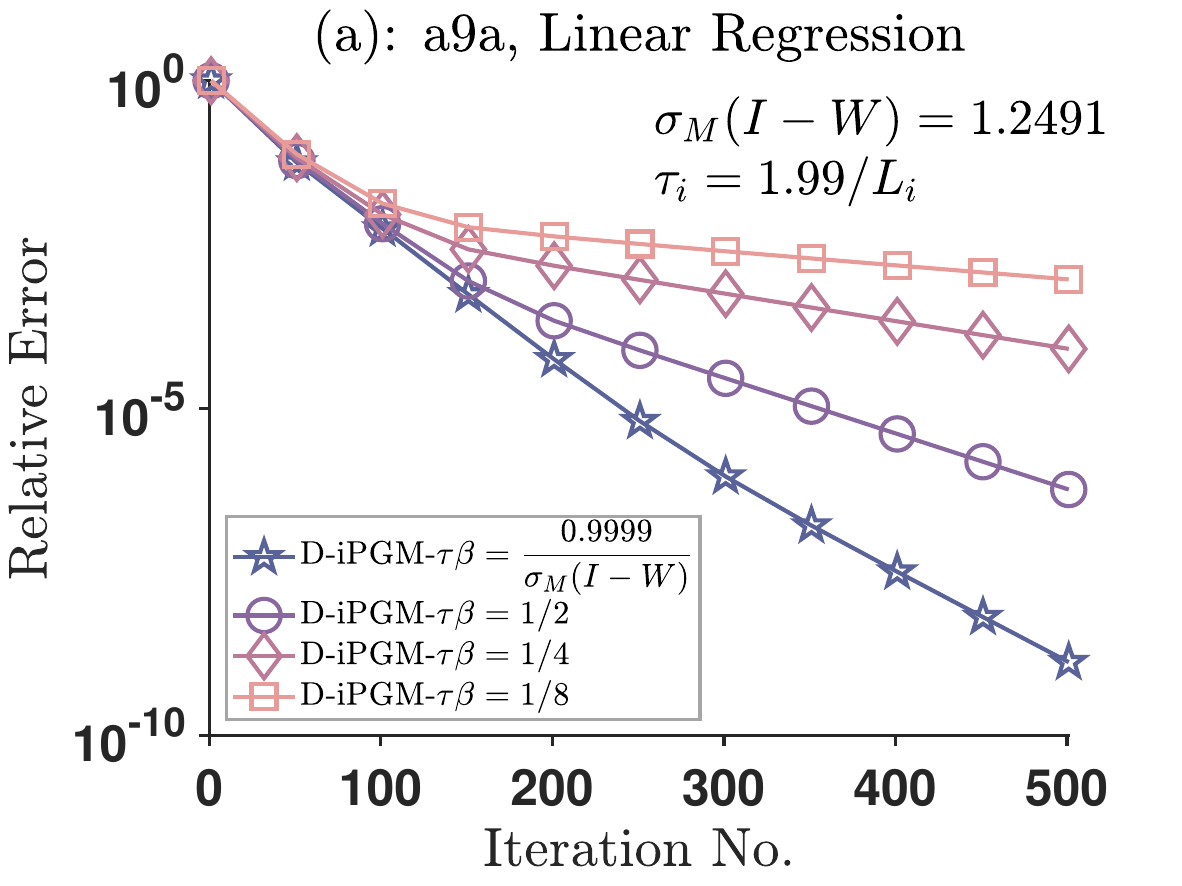}}
\subfigure{
\includegraphics[width=0.32\linewidth]{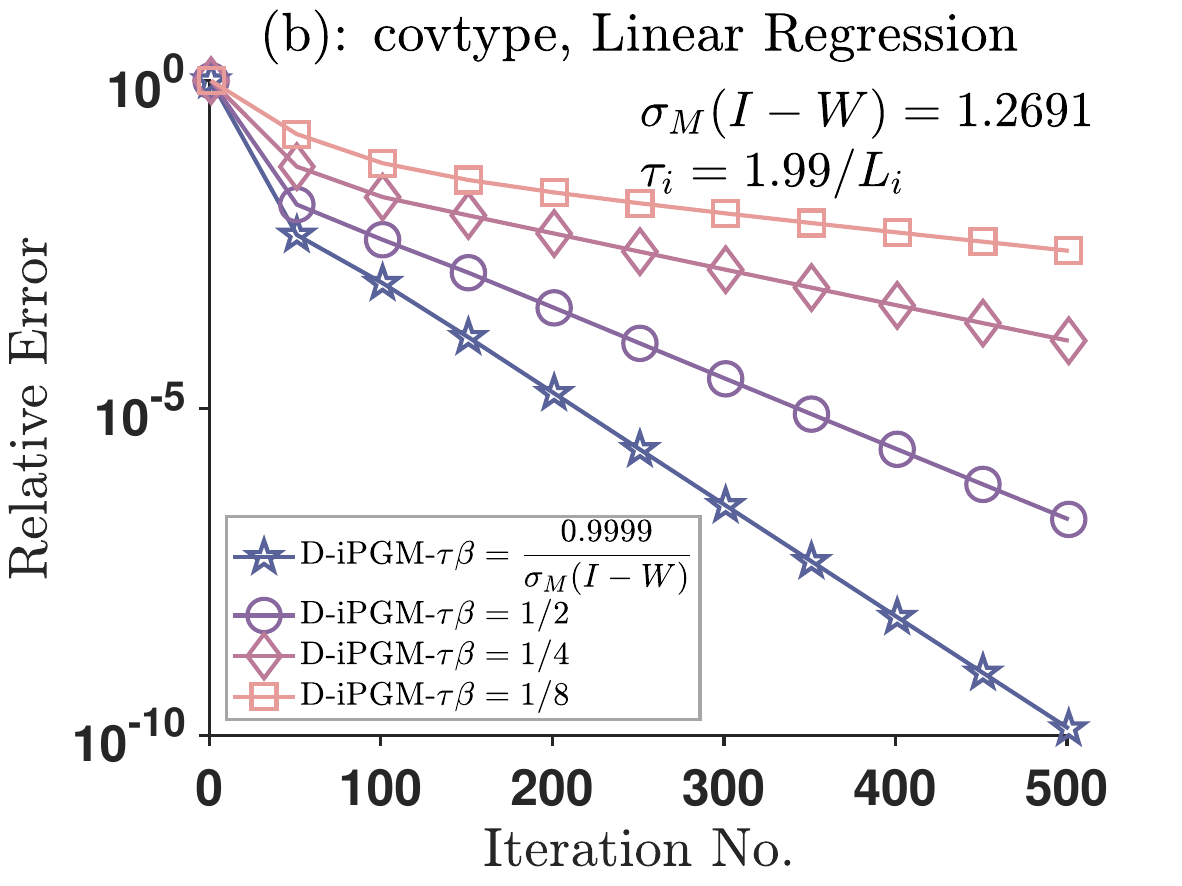}}
\subfigure{
\includegraphics[width=0.32\linewidth]{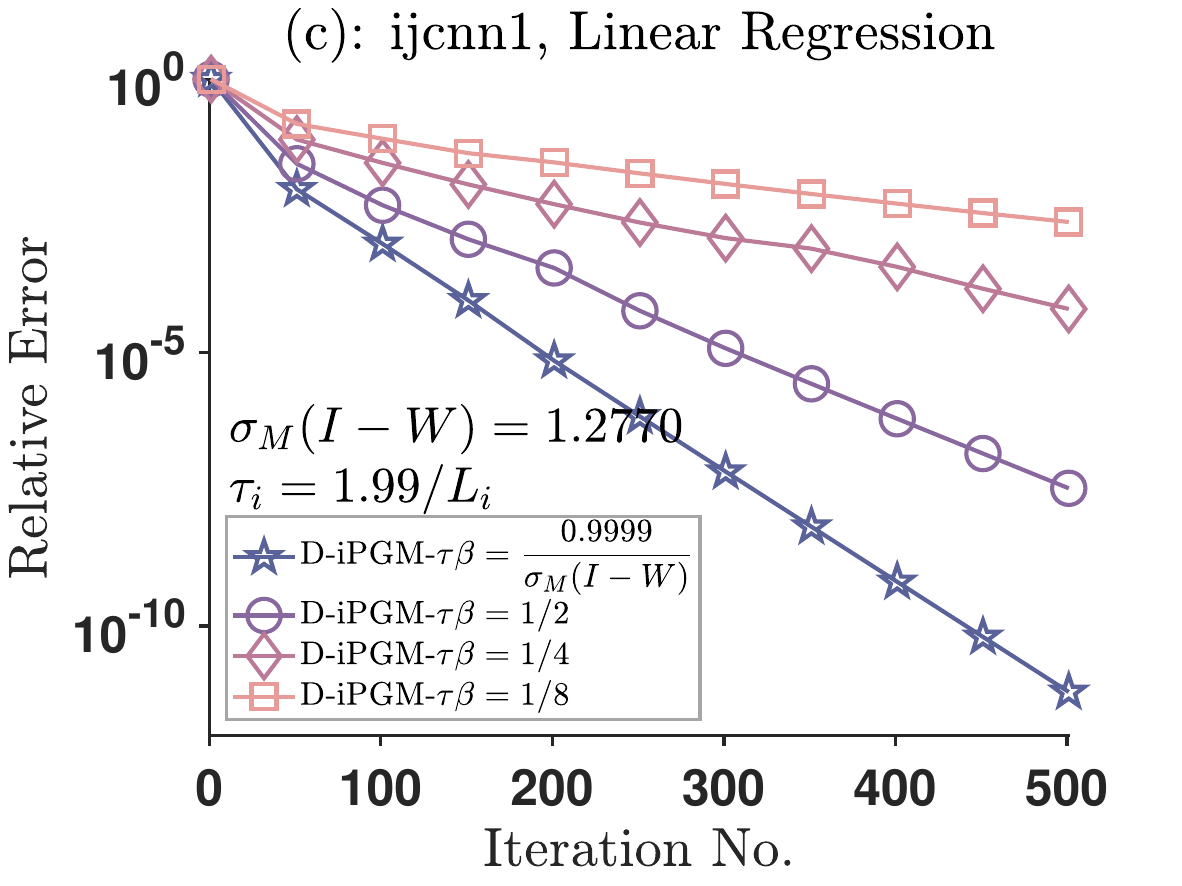}}
\caption{The relative error $\|\M{x}^k-\M{x}^*\|/\|\M{x}^*\|$ with respect to iteration numbers when the regularizer is $\nu_{1,i}\|\m{x}\|_1+\frac{\nu_{2,i}}{2}\|\m{x}\|^2$, where $\nu_{1,i}=0.01$ and $\nu_{2,i}=1$. Here the each agent has its own local stepsize $\tau_i=1.99/L_i$. We fix $\tau_i$ and let $\tau\beta=0.9999/\sigma_M(\M{I-W})$, $1/2$, $1/4$, and $1/8$, where $\tau=\max_{i\in\mathcal{V}}\{\tau_i\}$.
}
\label{Sim:case1.2}
\end{figure*}

\section{Numerical Experiments}\label{Sec5}
In this section, we compare the performance of D-iPGM with PG-EXTRA \cite{PGEXTR} (see \eqref{PD-PG-EXTRA} or \eqref{IT:EXTRA2}) and NIDS \cite{NDIS} (see \eqref{IT:NIDS1} or \eqref{IT:NIDS2}) for decentralized linear and logistic regression problems with different regularizers on real datasets.  All the algorithms are implemented in Matlab R2020b in a computer with 3.30 GHz AMD Ryzen 9 5900HS with Radeon Graphics and 16 GB memory.

For all experiments, we first compute the solution $\m{x}^*$ to \eqref{Problem} by centralized methods, and then run over a randomly generated connected network with $N$ agents and $\frac{\iota N(N-1)}{2}$ undirected edges, where $N=50$ and $\iota=0.1$ is the connectivity ratio. The mixing matrix $W$ is generated with the Metropolis-Hastings rule (see \cite[Section 2.4]{EXTRA}). The performance is evaluated by the relative error $\|\M{x}^k-\M{x}^*\|/\|\M{x}^*\|$. For decentralized linear regression, the loss function is
$$
\sum_{i=1}^{N}\big\{\frac{1}{m_i}\sum_{j=1}^{m_{i}} \frac{1}{2}\|\mathcal{A}_{ij}\tr{\m{x}}-\mathcal{B}_{ij}\|^2\big\}.
$$
For decentralized logistic regression, the loss function is
$$
\sum_{i=1}^{N}\big\{\frac{1}{m_i}\sum_{j=1}^{m_{i}} \ln (1+e^{-(\mathcal{A}_{i j}\tr {\m{x}})\mathcal{B}_{i j}})\big\}.
$$
Here, any agent $i$ holds its own training date $\left(\mathcal{A}_{i j}, \mathcal{B}_{i j}\right) \in$ $\mathbb{R}^{n} \times\{-1,1\}, j=1, \cdots, m_{i}$, including sample vectors $\mathcal{A}_{i j}$ and corresponding classes $\mathcal{B}_{i j}$. We use three real datasets including a9a, covtype, and ijcnn1\footnote[1]{\href{https://www.csie.ntu.edu.tw/~cjlin/libsvmtools/datasets/}{https://www.csie.ntu.edu.tw/~cjlin/libsvmtools/datasets/}}, whose attributes are $n=123$ and $\sum_{i=1}^{N}m_i=32550$, $n=54$ and $\sum_{i=1}^{N}m_i=55500$, and $n=22$ and $\sum_{i=1}^{N}m_i=49950$, respectively. Moreover, the training samples are
randomly and evenly distributed over all the $N$ agents.

\begin{table*}[!t]
\renewcommand\arraystretch{1.3}
\begin{center}
\caption{Comparison of Decentralized Linear Regression With Generalized LASSO Solved by PG-EXTRA, NIDS, and D-iPGM}
\scalebox{1}{
\begin{tabular}{cccccccccc}
\hline
\multirow{2}{*}{\textbf{Algorithm}} & \multicolumn{3}{c}{\textbf{a9a: Linear Regression}} & \multicolumn{3}{c}{\textbf{covtype: Linear Regression}} & \multicolumn{3}{c}{\textbf{ijcnn1: Linear Regression}}\\
\cline{2-10}
                    &Inner Iter. & Outer Iter.&Time&Inner Iter. & Outer Iter.&Time&Inner Iter. & Outer Iter.&Time\\
\hline
PG-EXTRA&32687 &146 &14.1916 & 63843& 186& 6.8362&42901&139&1.7367\\
NIDS&31678 & 141& 13.7203&61940&  180& 6.7148 &41027&133&1.6229\\
D-iPGM-$\varepsilon_k=10^{-10}$&31460&140&13.1901&61809&180&6.4375&40528&132&1.6263\\
D-iPGM-$\varepsilon_k=O(1/k^2)$&\textbf{18254} & 140&\textbf{11.4739}& \textbf{35067}& 185& \textbf{5.1945}&\textbf{27595}& 132&\textbf{1.2733}\\
D-iPGM-$\varepsilon_k=O(\ln k/k^2)$&\textbf{16640}& 146&\textbf{11.2607}&\textbf{33887}&189&\textbf{4.9323}&\textbf{24375}& 132&\textbf{1.1896}\\
D-iPGM-$\varepsilon_k=O(\|\m{x}_i^k-\m{x}_i^{k-1}\|/k)$& \textbf{11836}&142 &\textbf{10.6125} &\textbf{23520} &184 & \textbf{4.4608}&\textbf{13388}&132&\textbf{0.8935}\\
D-iPGM-$\varepsilon_k=O(\|\m{x}_i^k-\m{x}_i^{k-1}\|/\ln k)$&\textbf{8568} &141& \textbf{9.9773}&\textbf{17539} &185 & \textbf{3.9734}&\textbf{11529}&132&\textbf{0.8552}\\
\hline
\end{tabular}}
\label{Table-Comparison-LR}
\end{center}
\end{table*}

\begin{table*}[!t]
\renewcommand\arraystretch{1.3}
\begin{center}
\caption{Comparison of Decentralized Logistic Regression With Generalized LASSO Solved by PG-EXTRA, NIDS, and D-iPGM}
\scalebox{1}{
\begin{tabular}{cccccccccc}
\hline
\multirow{2}{*}{\textbf{Algorithm}} & \multicolumn{3}{c}{\textbf{a9a: Logistic Regression}} & \multicolumn{3}{c}{\textbf{covtype: Logistic Regression}} & \multicolumn{3}{c}{\textbf{ijcnn1: Logistic Regression}}\\
\cline{2-10}
                    &Inner Iter. & Outer Iter.&Time&Inner Iter. & Outer Iter.&Time&Inner Iter. & Outer Iter.&Time\\
\hline
PG-EXTRA&32912&169&15.3620& 52761&199 & 6.8032&54041&135&1.9368\\
NIDS&32910&169&15.3008&55084& 208&7.0974 &53794&130&1.8648\\
D-iPGM-$\varepsilon_k=10^{-10}$&33794&169&15.7742&55281&209&7.0255&54050&130&1.9370\\
D-iPGM-$\varepsilon_k=O(1/k^2)$&\textbf{17701}&165&\textbf{13.6317}&\textbf{31920} & 210&\textbf{5.8638} &\textbf{37284}&130&\textbf{1.4901}\\
D-iPGM-$\varepsilon_k=O(\ln k/k^2)$&\textbf{16055}&166 &\textbf{13.0957}&\textbf{30864}&211&\textbf{5.6870}&\textbf{32520}&130&\textbf{1.4303}\\
D-iPGM-$\varepsilon_k=O(\|\m{x}_i^k-\m{x}_i^{k-1}\|/k)$&\textbf{13458}&166&\textbf{12.9017}&\textbf{19383} & 160& \textbf{4.1614}&\textbf{17887}&123&\textbf{1.0219}\\
D-iPGM-$\varepsilon_k=O(\|\m{x}_i^k-\m{x}_i^{k-1}\|/\ln k)$&\textbf{8521} & 156& \textbf{11.5524}&\textbf{18795} & 173& \textbf{4.2117}&\textbf{13759}&117&\textbf{0.8875}\\
\hline
\end{tabular}}
\label{Table-Comparison-Log}
\end{center}
\end{table*}

\subsection{Scenario 1: $\mathrm{prox}_{\tau_i g_i}(\cdot)$ Has a Closed-form Solution}\label{Scenario1}
We first consider the decentralized linear and logistic regression problems with $\nu_{1,i}\|\m{x}\|+\frac{\nu_{2,i}}{2}\|\m{x}\|^2$ regularizer. To these problems the nonsmooth term is $g_i(\m{x})=\nu_{1,i}\|x\|$, and we have
$
\mathrm{prox}_{\tau_i g_i}(\m{x})=(1-\frac{\tau_i\nu_{1,i}}{\max\{\|\m{x}\|,\tau_i\nu_{1,i}\}})\m{x}.
$
It is not difficult to verify that the smooth coefficient of these two problems are $L_i=\|\mathcal{A}_i\tr\mathcal{A}_i\|+1$, where $\mathcal{A}_i=[\mathcal{A}_{i1}\tr;\cdots;\mathcal{A}_{im_i}\tr]$. In this case, since the proximal mapping of $g_i$ has a closed-form solution, we use the exact version of D-iPGM, i.e., $\varepsilon_k=0,k\geq0$. Thus, the per-iteration complexity of these three algorithms are the same. Note that there is only one round of communication in each iteration of these algorithms. The amount of information exchange over the network is proportional to their numbers of iterations. Consequently, only the number of iterations is recorded.

The comparison results of relative error for these algorithms after each iteration are shown as Fig. \ref{Sim:case1.1}, where we fix $\tau\beta=\frac{1}{2}$ for D-iPGM. For PG-EXTRA, we only show $\tau=1/L_F$, because when $\tau=1.4/L_F$ and $\tau=1.95/L_F$, PG-EXTRA is divergent. When $\tau=1/L_F$, these three algorithms have similar convergence performance. When $\tau=1.4/L_F$ and $\tau=1.95/L_F$, D-iPGM and NIDS converge at almost the same speed. In addition, D-iPGM and NIDS converge faster with a larger stepsize.

Then, we consider the decentralized linear regression problem with $\nu_{1,i}\|\m{x}\|_1+\frac{\nu_{2,i}}{2}\|\m{x}\|^2$ regularizer. Thus, the nonsmooth term is $g_i(\m{x})=\nu_{1,i}\|\m{x}\|_1$, it also has a closed-form solution. We use the uncoordinated steps $\tau_i=1.99/L_i$ and let $\tau\beta=0.9999/\sigma_M(\M{I-W})$, $1/2$, $1/4$, and $1/8$, where $\tau=\max_{i\in\mathcal{V}}\{\tau_i\}$. The relative error after each iteration are shown in Fig. \ref{Sim:case1.2}. It shows that D-iPGM converge faster with a larger $\tau\beta$.

\subsection{Scenario 2: $\mathrm{prox}_{\tau_i g_i}(\cdot)$ Has no Closed-form Solution}
In this subsection, we consider the decentralized linear and logistic regression problems with $l_2$+generalized LASSO \cite{Example3} regularizer, i.e., $\frac{\nu_{2,i}}{2}\|\m{x}\|^2+\nu_{1,i}\|\m{D}_i\m{x}\|_1$, where $\m{D}_i\in\mathbb{R}^{10\times n}$ is drawn from the normal distribution $N(0,1)$, $\nu_{1,i}=0.01$, and $\nu_{2,i}=1$. To this scenario, the nonsmooth term $g_i(\m{x})=\|\m{D}_i\m{x}\|_1$, whose proximal mapping has no closed-form solution for general $\m{D}_i$, and we will use the C-V algorithm \cite{CV1,CV2} to compute the numerical solution of $\mathrm{prox}_{\tau_i g_i}(\cdot)$.

Recall the primal update of D-iPGM,
$$
\tilde{\mathrm{x}}_i^{k}\approx \tilde{\mathrm{y}}_i^k:=\mathrm{prox}_{\tau_ig_i}(\mathrm{x}_i^k-\tau_i(\nabla f_i(\mathrm{x}_i^k)-{\lambda}_i^k)).
$$
Let $s_i^k=\mathrm{x}_i^k-\tau_i(\nabla f_i(\mathrm{x}_i^k)-{\lambda}_i^k)$. By the definition of proximal mapping of $g_i$, it holds that
\begin{align}\label{Sim:Scenario2:1}
\tilde{\mathrm{x}}_i^{k}\approx \tilde{\mathrm{y}}_i^k=\mathop{\arg\min}\limits_{\m{x}\in \mathbb{R}^n}\{\tau_i\nu_{1,i}\|\m{D}_i\m{x}\|_1+\frac{1}{2}\|\m{x}-\m{s}^k_i\|^2\}.
\end{align}
Applying the C-V algorithm to problem \eqref{Sim:Scenario2:1}, we get
\begin{align*}
\m{x}_i^{k,l+1}&=\m{x}_i^{k,l}-t_{1,i}(\m{x}_i^{k,l}-s_i^k+\m{D}_i\tr\m{w}_i^l),\\
\m{v}_i^{l+1}&=\m{w}_i^l+t_{2,i}\m{D}_i(2\m{x}_i^{k,l+1}-\m{x}_i^{k,l}),\\
\m{w}_i^{l+1}&=\m{v}_i^{l+1}-t_{2,i}\mathrm{prox}_{\frac{\tau_i\nu_{1,i}}{t_{2,i}}\|\cdot\|_1}(\m{v}_i^{l+1}/t_{2,i}),
\end{align*}
where $\m{x}_i^{k,0}=\m{x}_i^k$, $\m{w}_i^0=\M{0}$, and $t_{1,i}$ and $t_{2,i}$ are the stepsizes satisfying that $\frac{t_{1,i}}{2}+t_{1,i}t_{2,i}\|\m{D}_i\tr\m{D}_i\|<1$. Let
$\m{d}_i^{k,l+1}\in \partial g_i(\m{x}_i^{k,l+1})+\nabla f_i(\mathrm{x}_i^k)-\lambda_i^k+\frac{1}{\tau_i}(\m{x}_i^{k,l+1}-\mathrm{x}_i^{k})$.
It follows from \cite[Section 7]{PDS2021} and \cite[Lemma 3]{Jiangfan2022} that there exists a constant $\varepsilon_0$ such that $\|\m{d}_i^{k,l+1}\|\leq \varepsilon_0 \|\m{x}_i^{k,l+1}-\m{x}_i^{k,l}\|$. Therefore, if $\|\m{x}_i^{k,l+1}-\m{x}_i^{k,l}\|<{\varepsilon_k}/{\varepsilon_0}$, we stop the inner procedure, and set $\tilde{\m{x}}_i^{k}=\m{x}_i^{k,l+1}$, i.e., we can choose the inner stopping criterion of D-iPGM as
$$
\text{[D-iPGM]: }\|\m{x}_i^{k,l+1}-\m{x}_i^{k,l}\|\leq\varepsilon_k.
$$
To this scenario, we use the primal-dual forms of PG-EXTRA \eqref{PD-PG-EXTRA} and NIDS \eqref{IT:NIDS1} for solving the considered problems. Note that the primal update of these three algorithms are the same, and PG-EXTRA and NIDS are exact iterative algorithms. We apply the C-V algorithm to compute \eqref{PD-PG-EXTRA11} and \eqref{IT:NIDS11} with the inner stopping criterion
$$
\text{[PG-EXTRA/NIDS]: }\|\m{x}_i^{k,l+1}-\m{x}_i^{k,l}\|\leq10^{-10}.
$$
In order to compare the effect of inexact iteration and to ensure that the performance of outer iteration is similar for all the considered algorithms, according to the experimental results of Scenario 1, we set $\tau={1}/{L_F}$, $\tau\beta={1}/{2}$ for the outer iteration. In addition, we use the stopping criterion ${\|\M{x}^k-\M{x}^*\|}/{\|\M{x}^*\|}\leq10^{-5}$ for all considered methods.
\begin{figure}[!h]
  \centering
  \includegraphics[width=1\linewidth]{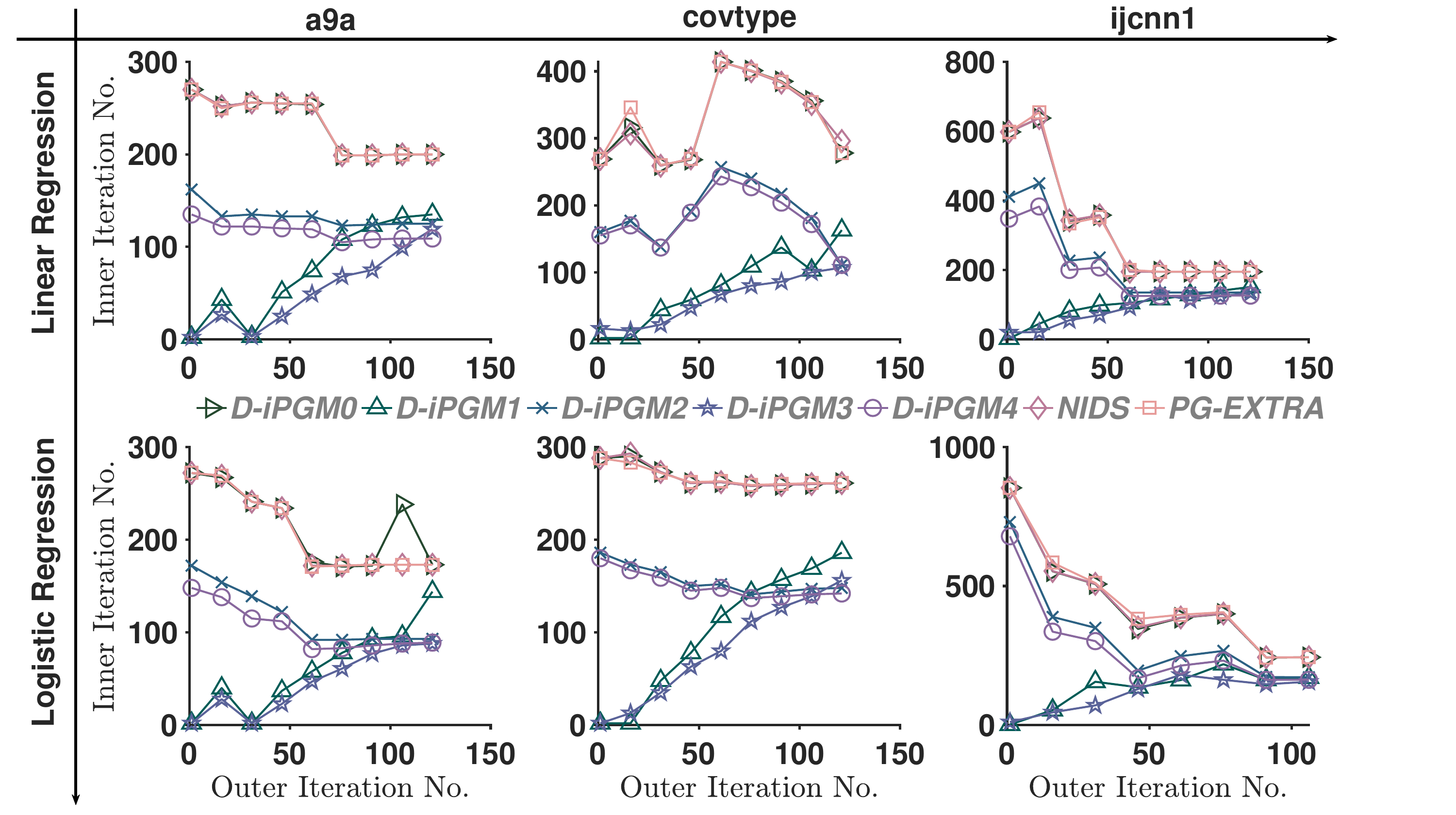}
  \caption{The number of inner iterations against the outer iterations. D-iPGM0, D-iPGM1, D-iPGM2, D-iPGM3, and D-iPGM4 represent D-iPGM with $\varepsilon_k=10^{-10},O(\|\m{x}_i^k-\m{x}_i^{k-1}\|/k),O(1/k^2),O(\|\m{x}_i^k-\m{x}_i^{k-1}\|/\ln k),O(\ln k/k^2)$.}
  \label{Fig-2-inner}
\end{figure}

Note that different algorithms have different complexities for each iteration. In Tables \ref{Table-Comparison-LR} and \ref{Table-Comparison-Log}, we report the number of inner iterations (Inner Iter.) and outer iterations (Outer Iter.), as well as the computational time in seconds (Time), when each method achieves the outer iteration stopping criterion. From Tables \ref{Table-Comparison-LR} and \ref{Table-Comparison-Log}, we can observe that, the computational time and the inner iteration number of the D-iPGM are significantly less than exact algorithms including its exact version, and the outer iteration number of D-iPGM is very close to exact algorithms. Therefore, it is obvious that D-iPGM performs better than the exact version of D-iPGM, PG-EXTRA, and NIDS in terms of the computational time. Moreover, by comparing the performance of D-iPGM with various inner stopping criterions, the result suggests that it is better to choose the $\varepsilon_k=O(\|\m{x}_i^k-\m{x}_i^{k-1}\|/\ln k)$ as the inner stopping criterion. To further visualize the numerical results, we present Fig. \ref{Fig-2-inner} to show the required number of inner iterations for each iteration, which indicates that D-iPGM requires fewer inner iterations per iteration than exact algorithms. If we choose $\varepsilon_k=O(\|\m{x}_i^k-\m{x}_i^{k-1}\|/\ln k)$ or $\varepsilon_k=O(\|\m{x}_i^k-\m{x}_i^{k-1}\|/ k)$ as the inner stopping criterion, the number of inner iterations required for each iteration is roughly on the rise. If we choose $\varepsilon_k=O(1/k^2)$ or $\varepsilon_k=O(\ln k/k^2)$, the trend of the number of internal iterations required for each iteration is similar to that of the exact algorithms.


\section{Conclusion}\label{Sec6}
In this paper, we proposed a CTA-based decentralized algorithm called D-iPGM for convex composite optimizations over networks, utilizing uncoordinated network-independent constant stepsizes. D-iPGM does not require the exact solution of the proximal mapping of $g_i$ in each iteration, making it suitable for decentralized composite optimizations where the proximal mapping may be difficult to compute. We established the convergence rate of D-iPGM to be $O(1/k)$ with general convexity, which can be improved to $o(1/k)$ if the proximal mapping is solved exactly. With metric subregularity, we further established its linear convergence rate. The numerical experiments confirmed the effectiveness of D-iPGM. Acceleration, asynchronous, and stochastic settings of D-iPGM were left as future work.

\appendices
\section{Proof of Lemma \ref{Lemma:Prediction1}}\label{Appendix:Lemma2}
\begin{IEEEproof}
It follows from \eqref{DALM2.1} that
$G(\mathbf{x})-G(\tilde{\mathbf{x}}^{k})+\langle\mathbf{x}-\tilde{\mathbf{x}}^{k},\nabla F(\mathbf{x}^{k})-\M{V} \bm{\alpha}^{k}+\Gamma^{-1} (\tilde{\mathbf{x}}^{k}-\mathbf{x}^{k})-\mathbf{d}^k\rangle \geq 0, \forall \mathbf{x}\in\mathbb{R}^{Nn}$,
which implies that for $\forall \mathbf{x} \in \mathbb{R}^{Nn}$,
\begin{align}\label{Proof:lemma1:1}
&G(\mathbf{x})-G(\tilde{\mathbf{x}}^{k})+\langle\mathbf{x}-\tilde{\mathbf{x}}^{k},\nabla F(\mathbf{x}^{k})-\M{V} \tilde{\bm{\alpha}}^{k}\nonumber\\
&+\Gamma^{-1} (\tilde{\mathbf{x}}^{k}-\mathbf{x}^{k})+\M{V}(\tilde{\bm{\alpha}}^{k}-\bm{\alpha}^{k})-\mathbf{d}^k\rangle \geq 0.
\end{align}
It follows from \eqref{DALM2.2} that for $\forall \bm{\alpha}\in \mathbb{R}^{Nn}$
\begin{align}\label{Proof:lemma1:2}
&\langle\bm{\alpha}-\tilde{\bm{\alpha}}^{k},\M{V}\tilde{\mathbf{x}}^k\rangle
=\langle\bm{\alpha}-\tilde{\bm{\alpha}}^{k},\frac{1}{\beta}(\bm{\alpha}^k-\tilde{\bm{\alpha}}^{k})\rangle.
\end{align}
Combining \eqref{Proof:lemma1:1} and \eqref{Proof:lemma1:2}, it gives that
\begin{align}\label{Proof:lemma1:3}
&G(\mathbf{x})-G(\tilde{\mathbf{x}}^{k})+\langle\mathbf{x}-\tilde{\mathbf{x}}^{k},\nabla F(\mathbf{x}^{k})\rangle+\langle\mathbf{u}-\tilde{\mathbf{u}}^{k},\mathcal{K}_2(\tilde{\mathbf{u}}^k)\rangle\nonumber\\
&\geq\langle\mathbf{u}-\tilde{\mathbf{u}}^{k}, \M{Q}(\mathbf{u}^{k}-\tilde{\mathbf{u}}^{k})\rangle+\langle\mathbf{x}-\tilde{\mathbf{x}}^{k},\mathbf{d}^k\rangle, \forall \mathbf{u} \in \mathcal{M}.
\end{align}
Using \eqref{skew-symmetric}, \eqref{Proof:lemma1:3} can be rewritten as
\begin{align}\label{Proof:lemma1:4}
&G(\mathbf{x})-G(\tilde{\mathbf{x}}^{k})+\langle\mathbf{x}-\tilde{\mathbf{x}}^{k},\nabla F(\mathbf{x}^{k})\rangle+\langle\mathbf{u}-\tilde{\mathbf{u}}^{k},\mathcal{K}_2(\M{u})\rangle\nonumber\\
&\geq\langle\mathbf{u}-\tilde{\mathbf{u}}^{k}, \M{Q}(\mathbf{u}^{k}-\tilde{\mathbf{u}}^{k})\rangle+\langle\mathbf{x}-\tilde{\mathbf{x}}^{k},\mathbf{d}^k\rangle, \forall \mathbf{u} \in \mathcal{M}.
\end{align}
Note that
$$
\langle\mathbf{u}-\tilde{\mathbf{u}}^{k},\mathcal{K}_1(\M{u})\rangle=\langle\mathbf{x}-\tilde{\mathbf{x}}^{k},\nabla F(\mathbf{x})\rangle+\langle\mathbf{u}-\tilde{\mathbf{u}}^{k},\mathcal{K}_2(\M{u})\rangle.
$$
By rearranging some terms of \eqref{Proof:lemma1:4}, it holds that for $\forall \M{u}\in \mathcal{M}$
\begin{align*}
&G(\mathbf{x})-G(\tilde{\mathbf{x}}^{k})+\langle\mathbf{u}-\tilde{\mathbf{u}}^{k},\mathcal{K}_1(\M{u})\rangle\\
&\geq -\langle \M{x}-\tilde{\M{x}}^k, \nabla F(\M{x}^k)-\nabla F(\M{x}) \rangle\\
&\quad+\langle\mathbf{u}-\tilde{\mathbf{u}}^{k}, \M{Q}(\mathbf{u}^{k}-\tilde{\mathbf{u}}^{k})\rangle+\langle\mathbf{x}-\tilde{\mathbf{x}}^{k},\mathbf{d}^k\rangle.
\end{align*}
By \eqref{SM1}, the assertion \eqref{Lemma:P1} holds. Combining \eqref{SM2} and \eqref{Proof:lemma1:4}, the assertion \eqref{Lemma:P2} holds.
\end{IEEEproof}
\section{Proof of Theorem \ref{TH1}}\label{Appendix:Theorem1}
\begin{IEEEproof}
1) According to \eqref{DALM2.3} and $\M{Q}=\M{HM}$, one has that
\begin{align*}
\langle\mathbf{u}-\tilde{\mathbf{u}}^{k}, \M{Q}(\mathbf{u}^{k}-\tilde{\mathbf{u}}^{k})\rangle=\langle\mathbf{u}-\tilde{\mathbf{u}}^{k}, \M{H}(\mathbf{u}^{k}-{\mathbf{u}}^{k+1})\rangle.
\end{align*}
Applying the identity $\langle a-b,\M{H}(c-d)\rangle=\frac{1}{2}\{\|a-d\|^2_{\M{H}}-\|a-c\|^2_{\M{H}}+\|b-c\|^2_{\M{H}}-\|b-d\|^2_{\M{H}}\}$ to the $\langle\tilde{\mathbf{u}}^{k}-\mathbf{u}, \mathrm{H}(\mathbf{u}^{k+1}-\mathbf{u}^{k})\rangle$ with specifications $a=\tilde{\mathbf{u}}^k$, $b=\mathbf{u}$, $c=\mathbf{u}^{k+1}$ and $d=\mathbf{u}^k$, it holds that
\begin{align}\label{Proof:theorem1:1}
&\langle\tilde{\mathbf{u}}^{k}-\mathbf{u}, \M{H}(\mathbf{u}^{k+1}-\mathbf{u}^{k})\rangle=\frac{1}{2}\{\|\tilde{\mathbf{u}}^k-\mathbf{u}^k\|^2_{\M{H}}\nonumber\\
&\quad-\|\mathbf{u}^{k+1}-\tilde{\mathbf{u}}^k\|^2_{\M{H}}+\|\mathbf{u}^{k+1}-\mathbf{u}\|^2_{\M{H}}-\|\mathbf{u}^k-\mathbf{u}\|^2_{\M{H}}\}.
\end{align}
Note that $\M{Q}=\M{HM}$ and $\M{G}=\M{Q}\tr+\M{Q}-\M{M}\tr\M{H}\M{M}$, one has
\begin{align}\label{Proof:theorem1:2}
&\|\mathbf{u}^k-\tilde{\mathbf{u}}^k\|^2_{\M{H}}-\|\mathbf{u}^{k+1}-\tilde{\mathbf{u}}^k\|^2_{\M{H}}\nonumber\\
&=\|\mathbf{u}^k-\tilde{\mathbf{u}}^k\|^2_{\M{H}}-\|(\M{u}^k-\tilde{\M{u}}^k)-\M{M}(\M{u}^k-\tilde{\M{u}}^k)\|^2_{\M{H}}\nonumber\\
&=2\langle \mathbf{u}^k-\tilde{\mathbf{u}}^k,\M{HM}(\mathbf{u}^k-\tilde{\mathbf{u}}^k) \rangle\nonumber\\
&\quad-\langle\mathbf{u}^k-\tilde{\mathbf{u}}^k,\M{M}\tr\M{H}\M{M}(\mathbf{u}^k-\tilde{\mathbf{u}}^k)\rangle\nonumber\\
&=\langle\mathbf{u}^k-\tilde{\mathbf{u}}^k,(\M{Q}\tr+\M{Q}-\M{M}\tr\M{HM})(\mathbf{u}^k-\tilde{\mathbf{u}}^k)\rangle\nonumber\\
&=\|\mathbf{u}^k-\tilde{\mathbf{u}}^k\|^2_{\M{G}}.
\end{align}
Thus, substituting \eqref{Proof:theorem1:1} and \eqref{Proof:theorem1:2} into \eqref{Lemma:P1}, it holds that
\begin{align}\label{Proof:theorem1:3}
&G(\mathbf{x})-G(\tilde{\mathbf{x}}^{k})+\langle\mathbf{u}-\tilde{\mathbf{u}}^{k},\mathcal{K}_1(\mathbf{u})\rangle \geq \langle\mathbf{x}-\tilde{\mathbf{x}}^k,\mathbf{d}^k\rangle\nonumber\\
&+\frac{1}{2}(\|\M{u}-\M{u}^{k+1}\|^2_{\M{H}}-\|\M{u}-\M{u}^{k}\|^2_{\M{H}})+\frac{1}{2}\|\mathbf{u}^k-\tilde{\mathbf{u}}^k\|^2_{\widehat{\M{H}}_1}.
\end{align}
Letting $\M{u}$ in \eqref{Proof:theorem1:3} as arbitrary $\M{u}^*\in\mathcal{M}^*$, we have
\begin{align*}
&\frac{1}{2}\|\mathbf{u}^{k+1}-\mathbf{u}^*\|^2_{\M{H}}\leq \frac{1}{2}\|\mathbf{u}^k-\mathbf{u}^*\|^2_{\M{H}}-\frac{1}{2}\|\mathbf{u}^k-\tilde{\mathbf{u}}^{k}\|^2_{\widehat{\M{H}}_1}\\
&+\langle\tilde{\mathbf{x}}^k-\mathbf{x}^*,\mathbf{d}^k\rangle
+G(\mathbf{x}^*)-G(\tilde{\mathbf{x}}^{k})+\langle\mathbf{u}^*-\tilde{\mathbf{u}}^{k},\mathcal{K}_1(\mathbf{u}^*)\rangle.
\end{align*}
Then, by \eqref{VI1}, the assertion \eqref{Theorem1:P1} holds.

2) We introduce the following notations.
\begin{subequations}\label{EXDiPGM}
\begin{align}
\tilde{\mathbf{y}}^{k}&=\mathrm{prox}_G^{\Gamma^{-1}}(\mathbf{x}^k-\Gamma(\nabla F(\mathbf{x}^k)-\M{V}\bm{\alpha}^k)),\\
\tilde{\M{z}}^{k}&=\bm{\alpha}^k-\beta \M{V}\tilde{\mathbf{y}}^{k},\\
\M{v}^{k+1}&=\M{u}^k-\M{M}(\M{u}^k-\tilde{\M{v}}^k),
\end{align}
\end{subequations}
where $\tilde{\M{v}}^k=\mathrm{col}\{\tilde{\mathbf{y}}^{k},\tilde{\M{z}}^{k}\}$. It follows from \eqref{Theorem1:P1} that
\begin{align*}
\|\mathbf{v}^{k+1}-\mathbf{u}^*\|^2_{\M{H}}\leq \|\mathbf{u}^k-\mathbf{u}^*\|^2_{\M{H}}-\|\mathbf{u}^k-\tilde{\mathbf{v}}^{k}\|^2_{\widehat{\M{H}}_1}, \forall k\geq0,
\end{align*}
which implies that $\|\mathbf{v}^{k+1}-\mathbf{u}^*\|_{\M{H}}\leq \|\mathbf{u}^k-\mathbf{u}^*\|_{\M{H}},k\geq0$. Hence, for $\forall \M{u}^*\in\mathcal{M}^*$ one has
\begin{align}\label{Proof:theorem1:4}
\|\M{u}^{k+1}-\M{u}^*\|_{\M{H}}&\leq\|\M{v}^{k+1}-\M{u}^*\|_{\M{H}}+\|\M{u}^{k+1}-\M{v}^{k+1}\|_{\M{H}}\nonumber\\
&\leq \|\mathbf{u}^k-\mathbf{u}^*\|_{\M{H}}+\|\M{u}^{k+1}-\M{v}^{k+1}\|_{\M{H}}.
\end{align}
By the nonexpansivity of proximal mapping, it holds that
$$
\|\tilde{\mathbf{x}}^{k}-\tilde{\mathbf{y}}^{k}\|\leq \sigma_M(\Gamma)\|\M{d}^k\|\leq\sqrt{N}\sigma_M(\Gamma)\varepsilon_k, \forall k\geq0.
$$
Therefore, it follows that for $\forall k\geq0$
\begin{align}\label{Proof:theorem1:5}
&\|\M{u}^{k+1}-\M{v}^{k+1}\|_{\M{H}}=\|\M{M}(\tilde{\M{u}}^k-\tilde{\M{v}}^k)\|_{\M{H}}\nonumber\\
&\leq\sigma_M(\M{H}^{\frac{1}{2}}\M{M})(\|\tilde{\mathbf{x}}^{k}-\tilde{\mathbf{y}}^{k}\|+\|\beta\mathbf{V}(\tilde{\mathbf{x}}^{k}-\tilde{\mathbf{y}}^{k})\|)\nonumber\\
&\leq\underbrace{\sqrt{N}\sigma_M(\M{H}^{\frac{1}{2}}\M{M})\sigma_M(\Gamma)(\beta\sigma_M(\M{V})+1)}_{:=\bar{\mu}>0}\varepsilon_k.
\end{align}
Combining \eqref{Proof:theorem1:4} and \eqref{Proof:theorem1:5}, it holds that
\begin{align}\label{Proof:theorem1:6}
\|\M{u}^{k+1}-\M{u}^*\|_{\M{H}}\leq \|\mathbf{u}^k-\mathbf{u}^*\|_{\M{H}}+\bar{\mu}\varepsilon_k,\forall \M{u}^*\in\mathcal{M}^*.
\end{align}
Summing the inequality \eqref{Proof:theorem1:6} over $k=0,1,\cdots,\mathrm{K}-1$, one has for any $\mathrm{K}\geq 1$,
\begin{align*}
\sum_{k=0}^{\mathrm{K}-1}\left\{\|\mathbf{u}^{k+1}-\mathbf{u}^{*}\|_{\M{H}}-\|\mathbf{u}^k-\mathbf{u}^{*}\|_{\M{H}}\right\}\leq \sum_{k=0}^{\mathrm{K}-1}\bar{\mu}\varepsilon_k,
\end{align*}
which implies
$$\|\mathbf{u}^{\mathrm{K}}-\mathbf{u}^*\|_{\M{H}}\leq \|\mathbf{u}^0-\mathbf{u}^*\|_{\M{H}}+\sum_{k=0}^{\mathrm{K}-1}\bar{\mu}\varepsilon_k, \forall \mathrm{K}\geq 1.
$$
Since $\M{H}\succ 0$ and $\{\varepsilon_k\}$ is summable, it deduces that for any $\mathbf{u}^*\in \mathcal{M}^*$, $\{\|\mathbf{u}^k-\mathbf{u}^*\|\}$ is bounded.

3) Summing \eqref{Theorem1:P1} over $k=0,1,\cdots,\infty$, one obtains
$$
\sum_{k=0}^{\infty}\|\M{u}^k-\tilde{\M{u}}^{k}\|^2_{\widehat{\M{H}}_1}\leq \|\M{u}^{0}-\M{u}^*\|_{\M{H}}^2 + \sum_{k=0}^{\infty}2\sqrt{N}\|\tilde{\M{x}}^{k}-\M{x}^*\|\varepsilon_k.
$$
Since $\{\varepsilon_k\}$ is summable and $\{\mathbf{u}^k\}$ is bounded, one has
\begin{align}\label{profit-1}
\sum_{k=0}^{\infty}\|\M{u}^k-\tilde{\M{u}}^{k}\|^2_{\widehat{\M{H}}_1}<\infty,
\end{align}
which implies that $\|\tilde{\mathbf{u}}^{k}-\mathbf{u}^k\|^2\rightarrow 0,k\rightarrow\infty$. Let $\mathbf{u}^{\infty}$ be a cluster point of $\{\tilde{\mathbf{u}}^k\}$ and $\{\tilde{\mathbf{u}}^{k_j}\}$ be a subsequence converging to $\mathbf{u}^{\infty}$. Then, it follows from \eqref{Proof:lemma1:3} that
\begin{align*}
&G(\mathbf{x})-G(\tilde{\mathbf{x}}^{k_j})+\langle\mathbf{x}-\tilde{\mathbf{x}}^{k_j},\nabla F(\mathbf{x}^{k_j})\rangle+\langle\mathbf{u}-\tilde{\mathbf{u}}^{k_j},\mathcal{K}_2(\tilde{\mathbf{u}}^{k_j})\rangle\\
&\geq\langle\mathbf{u}-\tilde{\mathbf{u}}^{k_j}, \M{Q}(\mathbf{u}^{k_j}-\tilde{\mathbf{u}}^{k_j})\rangle+\langle\mathbf{x}-\tilde{\mathbf{x}}^{k_j},\mathbf{d}^{k_j}\rangle, \forall \mathbf{u} \in \mathcal{M}.
\end{align*}
Since the matrix $\M{Q}$ is nonsingular, and $G$, $\nabla F$, and $\mathcal{K}_2$ are continuous, taking $k_j\rightarrow\infty$ in the above inequality, one gets
\begin{align*}
G(\mathbf{x})-G(\mathbf{x}^{\infty}) +\langle\mathbf{u}-\mathbf{u}^{\infty},\mathcal{K}_1(\mathbf{u}^{\infty})\rangle\geq 0, \forall \mathbf{u}\in \mathcal{M}.
\end{align*}
Compared with \eqref{VI1}, we can conclude that $\M{u}^{\infty}\in\mathcal{M}^*$. Therefore, by \eqref{Proof:theorem1:6}, it holds that
\begin{align*}
\|\M{u}^{k+1}-\M{u}^{\infty}\|_{\M{H}}\leq \|\mathbf{u}^k-\mathbf{u}^{\infty}\|_{\M{H}}+\bar{\mu}\varepsilon_k.
\end{align*}
Since $\sum_{k=0}^{\infty}\varepsilon_k<\infty$, by \cite[Lemma 3.2]{MP2021} the quasi-Fej\'{e}r monotone sequence $\{\|\M{u}^k-\M{u}^{\infty}\|_{\mathbf{H}}\}$ converges to a unique limit point. Then, with $\M{u}^{\infty}$ being an accumulation point of $\{\M{u}^k\}$, one has $\lim_{k\rightarrow \infty}\M{u}^k=\M{u}^{\infty}$.
\end{IEEEproof}

\section{Proof of Theorem \ref{TH2}}\label{Appendix:Theorem2}
\begin{IEEEproof}
By \eqref{DALM2.1} and \eqref{DALM2.2}, we have
\begin{align*}
&\M{0}\in \partial G(\tilde{\M{x}}^{k})+\nabla F(\M{x}^k)-\M{V}\bm{\alpha}^k+\Gamma^{-1}(\tilde{\M{x}}^{k}-\M{x}^k)-\M{d}^k,\\
&\M{0}=\beta{\M{V}}\tilde{\M{x}}^{k}+(\tilde{\bm{\alpha}}^k-\bm{\alpha}^k),
\end{align*}
which implies that
\begin{align*}
&\nabla F(\tilde{\M{x}}^k)-\nabla F(\M{x}^k)-\M{V}(\tilde{\bm{\alpha}}^k-\bm{\alpha}^k)-\Gamma^{-1}(\tilde{\M{x}}^{k}-\M{x}^k)+\M{d}^k\\
&\in \partial G(\tilde{\M{x}}^{k}) + \nabla F(\tilde{\M{x}}^k)-\M{V}\tilde{\bm{\alpha}}^k,
\end{align*}
and $-\frac{1}{\beta}(\tilde{\bm{\alpha}}^k-\bm{\alpha}^k)={\M{V}}\tilde{\M{x}}^{k}$.
Therefore, it holds that
\begin{align}\label{ADDproof}
&\mathrm{dist}^2(\M{0},\mathcal{J}(\tilde{\M{u}}^k))\leq\|\nabla F(\tilde{\M{x}}^k)-\nabla F(\M{x}^k)-\M{V}(\tilde{\bm{\alpha}}^k-\bm{\alpha}^k)\nonumber\\
&\quad-\Gamma^{-1}(\tilde{\M{x}}^{k}-\M{x}^k)+\M{d}^k\|^2+\|\frac{1}{\beta}(\tilde{\bm{\alpha}}^k-\bm{\alpha}^k)\|^2\nonumber\\
&\leq(4\max_i\{L^2_i\}+4\sigma_m(\Gamma^2))\|\tilde{\M{x}}^k-\M{x}^k\|^2\nonumber\\
&\quad+(4\sigma_M(\M{V}^2)+\frac{1}{\beta^2})\|\tilde{\bm{\alpha}}^k-\bm{\alpha}^k\|^2+4\|\M{d}^k\|^2\nonumber\\
&\leq \kappa_1^2\|\tilde{\M{u}}^k-\M{u}^k\|^2_{\widehat{\M{H}}_1}+4N(\varepsilon_k)^2,
\end{align}
where $\kappa_1^2=\frac{\max\{4\sigma_M(\M{V}^2)+\frac{1}{\beta^2},4\max_i\{L^2_i\}+4\sigma_m(\Gamma^{-2})\}}{\sigma_m(\widehat{\M{H}}_1)}$. Thus, by \eqref{profit-1} and \cite[Proposition 1]{PGEXTR}, the $O(1/(K+1))$ running-average optimality residual and the $o(1/(K+1))$ running-best optimality residual hold.

Substituting \eqref{Proof:theorem1:1} and \eqref{Proof:theorem1:2} into \eqref{Lemma:P2}, it holds that
\begin{align}\label{Proof:theorem2:1}
&\phi(\mathbf{x})-\phi(\tilde{\mathbf{x}}^{k})+\langle\mathbf{u}-\tilde{\mathbf{u}}^{k},\mathcal{K}_2(\mathbf{u})\rangle \geq \langle\mathbf{x}-\tilde{\mathbf{x}}^k,\mathbf{d}^k\rangle\nonumber\\
&+\frac{1}{2}(\|\M{u}-\M{u}^{k+1}\|^2_{\M{H}}-\|\M{u}-\M{u}^{k}\|^2_{\M{H}})+\frac{1}{2}\|\mathbf{u}^k-\tilde{\mathbf{u}}^k\|^2_{\widehat{\M{H}}_2}.
\end{align}
Note that
$$
L(\tilde{\M{x}}^k,\bm{\alpha})-L(\M{x},\tilde{\bm{\alpha}}^k)
=\phi(\tilde{\mathbf{x}}^{k})-\phi(\mathbf{x})+\langle\tilde{\mathbf{u}}^{k}-\mathbf{u},\mathcal{K}_2(\mathbf{u})\rangle.
$$
Summing \eqref{Proof:theorem2:1} over $k=0,1,\cdots,K-1$, we obtain
\begin{align*}
&2\sum_{k=0}^{K-1}({L}(\tilde{\M{x}}^{k},\bm{\alpha})-{L}(\M{x},\tilde{\bm{\alpha}}^{k}))\\
&\leq \|\M{u}^0-\M{u}\|_{\M{H}}^2 + 2\sqrt{N}\sum_{k=0}^{K-1}\|\tilde{\M{x}}^{k}-\M{x}\| \varepsilon_k.
\end{align*}
It follows from the convexity of $F$, $G$ and the definition of $(\M{X}^{K},\M{\Lambda}^K)$ that
$$
K(L(\M{X}^{K},\bm{\alpha})-L(\M{x},\M{\Lambda}^{K}))\leq\sum_{k=0}^{K-1}(L(\tilde{\M{x}}^{k},\bm{\alpha})-L(\M{x},\tilde{\bm{\alpha}}^{k})).
$$
Therefore, the primal–dual gap
\begin{align}\label{primal-dual-gap}
&L(\M{X}^{K},\bm{\alpha})-L(\M{x},\M{\Lambda}^{K})\nonumber\\
&\leq \frac{\|\M{u}^0-\M{u}\|_{\M{H}}^2+ 2\sqrt{N}\sum_{k=0}^{K-1}\|\tilde{\M{x}}^{k}-\M{x}\| \varepsilon_k}{2K}
\end{align}
holds. Since $\{\varepsilon_k\}$ is summable and $\{\mathbf{u}^k\}$ is bounded, the primal-dual gap is $O(1/K)$.

Note the inequality \eqref{primal-dual-gap} holds for all $(\M{x},\bm{\alpha})\in \mathcal{M}$, hence it is also true for any optimal solution $\M{x}^*$ and $\mathcal{B}_{\rho}=\{\bm{\alpha}:\|\bm{\alpha}\|\leq\rho\}$, where $\rho>0$ is a any given positive number. Letting $\M{u}=\mathrm{col}\{\M{x}^*,\bm{\alpha}\}$ and $\M{U}^K=\mathrm{col}\{\M{X}^K,\M{\Lambda}^K\}$, it gives
\begin{align}\label{Proof:theorem2:2}
&\sup_{\bm{\alpha}\in \mathcal{B}_{\rho}}\{\mathcal{L}(\M{X}^{K},\bm{\alpha})-\mathcal{L}(\M{x}^*,\M{\Lambda}^{K})\}\nonumber\\
&=\sup_{\bm{\alpha}\in \mathcal{B}_{\rho}}\{\phi(\M{X}^{K})-\phi(\mathbf{x}^*)+\langle{\mathbf{U}}^{K}-\mathbf{u},\mathcal{K}_2(\mathbf{u})\rangle\}\nonumber\\
&\overset{\eqref{skew-symmetric}}{=}\sup_{\bm{\alpha}\in \mathcal{B}_{\rho}}\{\phi(\M{X}^{K})-\phi(\mathbf{x}^*)+\langle{\mathbf{U}}^{K}-\mathbf{u},\mathcal{K}_2({\mathbf{U}}^{K})\rangle\}\nonumber\\
&=\sup_{\bm{\alpha}\in \mathcal{B}_{\rho}}\{\phi(\M{X}^{K})-\phi(\mathbf{x}^*)+\langle \M{X}^K-\M{x}^*,-\M{V}\M{\Lambda}^K \rangle\nonumber\\
&\quad\quad\quad\quad+\langle \M{\Lambda}^k-\bm{\alpha},\M{V}\M{X}^K \rangle\}\nonumber\\
&=\sup_{\bm{\alpha}\in \mathcal{B}_{\rho}}\{\phi(\M{X}^{K})-\phi(\mathbf{x}^*)-\langle\bm{\alpha},\M{V}\M{X}^K\rangle\}\nonumber\\
&=\phi(\M{X}^{K})-\phi(\mathbf{x}^*)+\rho\|\M{VX}^K\|.
\end{align}
Note that $\bm{\alpha}^0=\M{0}$. Combining \eqref{primal-dual-gap} and \eqref{Proof:theorem2:2}, one has
\begin{align*}
&\phi(\M{X}^{K})-\phi(\mathbf{x}^*)+\rho\|\M{VX}^K\|\\
&\leq \sup_{\bm{\alpha}\in \mathcal{B}_{\rho}}\{\frac{\|\M{u}^0-\M{u}\|_{\M{H}}^2+ 2\sqrt{N}\sum_{k=0}^{K-1}\|\tilde{\M{x}}^{k}-\M{x}^*\| {\varepsilon_k}}{2K}\}\\
&=\frac{\|\M{x}^0-\M{x}^*\|^2_{\Gamma^{-1}}+\frac{1}{\beta}\rho^2+2\sqrt{N}\sum_{k=0}^{K-1}\|\tilde{\M{x}}^{k}-\M{x}^*\| {\varepsilon_k}}{2K}.
\end{align*}
Since $\{\M{u}^k\}$ is bounded and $\sum_{k=0}^{\infty}{\varepsilon_k}<\infty$, there exists $0<\psi<\infty$ such that
\begin{align}\label{Proof:theorem2:3}
\phi(\M{X}^{K})-\phi(\mathbf{x}^*)+\rho\|\M{VX}^K\|\leq\frac{\psi}{K}.
\end{align}
For any $\rho>\|\bm{\alpha}^*\|$, applying \cite[Lemma 2.3]{Xuyangyang2017} in \eqref{Proof:theorem2:3}, one has
\begin{align*}
\|\M{VX}^K\|&\leq\frac{\psi}{K(\rho-\|\bm{\alpha}^*\|)},\\
-\frac{\psi\|\bm{\alpha}^*\|}{K(\rho-\|\bm{\alpha}^*\|)}&\leq\phi(\M{X}^K)-\phi(\M{x}^*)\leq\frac{\psi}{K}.
\end{align*}
Letting $\rho=\max\{1+\|\bm{\alpha}^*\|,2\|\bm{\alpha}^*\|\}$, we have
\begin{align*}
|\phi(\M{X}^K)-\phi(\M{x}^*)|\leq\frac{\psi}{K},
\|\sqrt{\M{I-W}}\M{X}^K\|\leq\frac{\psi}{K}.
\end{align*}
Thus, the primal gap and the consensus error are $O(1/K)$.
\end{IEEEproof}

\section{Proof of Theorem \ref{Theorem:point-wise-rate}}\label{Appendix:Theorem:point-wise-rate}
\begin{IEEEproof}
Since ${\varepsilon_k}=0$ for all $k\geq0$, by \eqref{Proof:lemma1:3}, one has
\begin{align}\label{Proof:Theorem:point-wise-rate:1}
&G(\mathbf{x})-G(\tilde{\mathbf{x}}^{k})+\langle\mathbf{x}-\tilde{\mathbf{x}}^{k},\nabla F(\mathbf{x}^{k})\rangle\nonumber\\
&\quad+\langle\mathbf{u}-\tilde{\mathbf{u}}^{k},\mathcal{K}_2(\tilde{\mathbf{u}}^k)\rangle\nonumber\\
&\geq\langle\mathbf{u}-\tilde{\mathbf{u}}^{k}, \M{Q}(\mathbf{u}^{k}-\tilde{\mathbf{u}}^{k})\rangle, \forall \mathbf{u} \in \mathcal{M}.
\end{align}
Letting $\M{u}=\tilde{\M{u}}^{k+1}$ in \eqref{Proof:Theorem:point-wise-rate:1}, one has
\begin{align}\label{Proof:Theorem:point-wise-rate:2}
&G(\tilde{\M{x}}^{k+1})-G(\tilde{\mathbf{x}}^{k})+\langle\tilde{\M{x}}^{k+1}-\tilde{\mathbf{x}}^{k},\nabla F(\mathbf{x}^{k})\rangle\nonumber\\
&\quad+\langle\tilde{\M{u}}^{k+1}-\tilde{\mathbf{u}}^{k},\mathcal{K}_2(\tilde{\mathbf{u}}^k)\rangle\nonumber\\
&\geq\langle\tilde{\M{u}}^{k+1}-\tilde{\mathbf{u}}^{k}, \M{Q}(\mathbf{u}^{k}-\tilde{\mathbf{u}}^{k})\rangle.
\end{align}
Rewriting the inequality \eqref{Proof:Theorem:point-wise-rate:1} for the $(k+1)$-th iteration, it gives that
\begin{align}\label{Proof:Theorem:point-wise-rate:3}
&G(\mathbf{x})-G(\tilde{\mathbf{x}}^{k+1})+\langle\mathbf{x}-\tilde{\mathbf{x}}^{k+1},\nabla F(\mathbf{x}^{k+1})\rangle\nonumber\\
&\quad+\langle\mathbf{u}-\tilde{\mathbf{u}}^{k+1},\mathcal{K}_2(\tilde{\mathbf{u}}^{k+1})\rangle\nonumber\\
&\geq\langle\mathbf{u}-\tilde{\mathbf{u}}^{k+1}, \M{Q}(\mathbf{u}^{k+1}-\tilde{\mathbf{u}}^{k+1})\rangle, \forall \mathbf{u} \in \mathcal{M}.
\end{align}
Setting $\M{u}=\tilde{\M{u}}^k$ in \eqref{Proof:Theorem:point-wise-rate:3}, one has
\begin{align}\label{Proof:Theorem:point-wise-rate:4}
&G(\tilde{\M{x}}^k)-G(\tilde{\mathbf{x}}^{k+1})+\langle\tilde{\M{x}}^k-\tilde{\mathbf{x}}^{k+1},\nabla F(\mathbf{x}^{k+1})\rangle\nonumber\\
&\quad+\langle\tilde{\M{u}}^k-\tilde{\mathbf{u}}^{k+1},\mathcal{K}_2(\tilde{\mathbf{u}}^{k+1})\rangle\nonumber\\
&\geq\langle\tilde{\M{u}}^k-\tilde{\mathbf{u}}^{k+1}, \M{Q}(\mathbf{u}^{k+1}-\tilde{\mathbf{u}}^{k+1})\rangle, \forall \mathbf{u} \in \mathcal{M}.
\end{align}
Note that $\langle \tilde{\M{u}}^k- \tilde{\M{u}}^{k+1},\mathcal{K}_2(\tilde{\M{u}}^k)-\mathcal{K}_2(\tilde{\M{u}}^{k+1})\rangle\equiv0$. Adding \eqref{Proof:Theorem:point-wise-rate:2} and \eqref{Proof:Theorem:point-wise-rate:4}, we obtain
\begin{align}\label{Proof:Theorem:point-wise-rate:5}
&\langle \tilde{\M{u}}^k-\tilde{\M{u}}^{k+1},\M{Q}((\M{u}^k-\tilde{\M{u}}^k)-(\M{u}^{k+1}-\tilde{\M{u}}^{k+1})) \rangle\nonumber\\
&\geq\langle\tilde{\M{x}}^k-\tilde{\mathbf{x}}^{k+1},\nabla F(\mathbf{x}^{k})-\nabla F(\mathbf{x}^{k+1})\rangle.
\end{align}
Then, adding the term
$$
\langle(\M{u}^k-\tilde{\M{u}}^k)-(\M{u}^{k+1}-\tilde{\M{u}}^{k+1}),\M{Q}((\M{u}^k-\tilde{\M{u}}^k)-(\M{u}^{k+1}-\tilde{\M{u}}^{k+1}))\rangle
$$
to both sides \eqref{Proof:Theorem:point-wise-rate:5} and using the identity
$$
\M{u}\tr\M{Q}\M{u}=\frac{1}{2}\M{u}\tr(\M{Q}\tr+\M{Q})\M{u}, \forall \M{u}\in \mathcal{M},
$$
one obtains that
\begin{align}\label{Proof:Theorem:point-wise-rate:6}
&\underbrace{\langle {\M{u}}^k-{\M{u}}^{k+1},\M{Q}((\M{u}^k-\tilde{\M{u}}^k)-(\M{u}^{k+1}-\tilde{\M{u}}^{k+1})) \rangle}_{(\m{I})}\nonumber\\
&\geq\underbrace{\langle\tilde{\M{x}}^k-\tilde{\mathbf{x}}^{k+1},\nabla F(\mathbf{x}^{k})-\nabla F(\mathbf{x}^{k+1})\rangle}_{(\m{II})}\nonumber\\
&\quad+\frac{1}{2}\|(\M{u}^k-\tilde{\M{u}}^k)-(\M{u}^{k+1}-\tilde{\M{u}}^{k+1})\|^2_{\M{Q}\tr+\M{Q}}.
\end{align}
On one hand, we have
\begin{align*}
(\m{I})=&\langle {\M{u}}^k-{\M{u}}^{k+1},\M{Q}((\M{u}^k-\tilde{\M{u}}^k)-(\M{u}^{k+1}-\tilde{\M{u}}^{k+1})) \rangle\\
\overset{\eqref{DALM2.3}}{=}&\langle \M{M}(\M{u}^k-\tilde{\M{u}}^{k}),\M{Q}((\M{u}^k-\tilde{\M{u}}^k)-(\M{u}^{k+1}-\tilde{\M{u}}^{k+1})) \rangle\\
\overset{\eqref{MAXCON}}{=}&\langle \M{u}^k- \tilde{\M{u}}^{k},\M{M}\tr\M{HM}((\M{u}^k-\tilde{\M{u}}^k)-(\M{u}^{k+1}-\tilde{\M{u}}^{k+1}))\rangle.
\end{align*}
On the other hand, we have
\begin{align*}
-(\m{II})=&\langle \nabla F(\M{x}^{k+1})-\nabla F(\M{x}^k),(\tilde{\M{x}}^k-\M{x}^k)-(\tilde{\M{x}}^{k+1}-\M{x}^{k+1}) \rangle\\
&+\langle \nabla F(\M{x}^{k+1})-\nabla F(\M{x}^k),\M{x}^k-\M{x}^{k+1}\rangle\\
\leq&\frac{\epsilon}{2}\|(\tilde{\M{x}}^k-\M{x}^k)-(\tilde{\M{x}}^{k+1}-\M{x}^{k+1})\|^2_{\M{L}_F}\\
&+\frac{1}{2\epsilon}\| \nabla F(\M{x}^{k+1})-\nabla F(\M{x}^k)\|^2_{\M{L}_F^{-1}}\\
&+\langle \nabla F(\M{x}^{k+1})-\nabla F(\M{x}^k),\M{x}^k-\M{x}^{k+1}\rangle\\
\overset{\eqref{SM3}}{\leq}&\frac{\epsilon}{2}\|(\tilde{\M{x}}^k-\M{x}^k)-(\tilde{\M{x}}^{k+1}-\M{x}^{k+1})\|^2_{\M{L}_F}\\
&+(1-\frac{1}{2\epsilon})\langle \nabla F(\M{x}^{k+1})-\nabla F(\M{x}^k),\M{x}^k-\M{x}^{k+1}\rangle.
\end{align*}
Setting $\epsilon=\frac{1}{2}$, one has that
$$
(\m{II})\geq-\frac{1}{4}\|(\tilde{\M{x}}^k-\M{x}^k)-(\tilde{\M{x}}^{k+1}-\M{x}^{k+1})\|^2_{\M{L}_F}.
$$
Thus, by \eqref{Proof:Theorem:point-wise-rate:6}, it holds that
\begin{align}\label{Proof:Theorem:point-wise-rate:7}
&2\langle \M{u}^k- \tilde{\M{u}}^{k},\M{M}\tr\M{HM}((\M{u}^k-\tilde{\M{u}}^k)-(\M{u}^{k+1}-\tilde{\M{u}}^{k+1}))\rangle\nonumber\\
&\geq \|(\M{u}^k-\tilde{\M{u}}^k)-(\M{u}^{k+1}-\tilde{\M{u}}^{k+1})\|^2_{\M{Q}\tr+\M{Q}} \nonumber\\
&\quad-\frac{1}{2}\|(\tilde{\M{x}}^k-\M{x}^k)-(\tilde{\M{x}}^{k+1}-\M{x}^{k+1})\|^2_{\M{L}_F}.
\end{align}
By the identity $\|a\|^2_{\M{H}}-\|b\|^2_{\M{H}}=2\langle a,\M{H}(a-b)\rangle-\|a-b\|^2_{\M{H}}$ with $a=\M{M}(\M{u}^k-\tilde{\M{u}}^k)$ and $b=\M{M}(\M{u}^{k+1}-\tilde{\M{u}}^{k+1})$ and $\widehat{\M{H}}_1=\M{Q}\tr+\M{Q}-\M{M}\tr\M{HM}-\frac{1}{2}\mathrm{diag}\{\M{L}_F,\M{0}\}$,
we have
\begin{align*}
&\|\M{M}(\M{u}^k-\tilde{\M{u}}^k)\|^2_{\M{H}}-\|\M{M}(\M{u}^{k+1}-\tilde{\M{u}}^{k+1})\|^2_{\M{H}}\nonumber\\
&=2\langle \M{u}^k- \tilde{\M{u}}^{k},\M{M}\tr\M{HM}((\M{u}^k-\tilde{\M{u}}^k)-(\M{u}^{k+1}-\tilde{\M{u}}^{k+1}))\rangle\nonumber\\
&\quad+\|\M{M}((\M{u}^k-\tilde{\M{u}}^k)-(\M{u}^{k+1}-\tilde{\M{u}}^{k+1}))\|^2_{\M{H}}\nonumber\\
&\overset{\eqref{Proof:Theorem:point-wise-rate:7}}{\geq}\|(\M{u}^k-\tilde{\M{u}}^k)-(\M{u}^{k+1}-\tilde{\M{u}}^{k+1})\|^2_{\M{Q}\tr+\M{Q}} \nonumber\\
&\quad+\|\M{M}((\M{u}^k-\tilde{\M{u}}^k)-(\M{u}^{k+1}-\tilde{\M{u}}^{k+1}))\|^2_{\M{H}}\nonumber\\
&\quad-\frac{1}{2}\|(\tilde{\M{x}}^k-\M{x}^k)-(\tilde{\M{x}}^{k+1}-\M{x}^{k+1})\|^2_{\M{L}_F}\nonumber\\
&=\|(\M{u}^k-\tilde{\M{u}}^k)-(\M{u}^{k+1}-\tilde{\M{u}}^{k+1})\|^2_{\widehat{\M{H}}_1}\geq0.
\end{align*}
Thus, \eqref{eq:point-wise-rate1} holds.

By \eqref{Theorem1:P1} and the equivalence of different norms, there exists a constant $c_0$ such that for any $\M{u}^*\in\mathcal{M}^*$
\begin{align}\label{Proof:Theorem:point-wise-rate:8}
\|\mathbf{u}^{k+1}-\mathbf{u}^*\|^2_{\M{H}}\leq \|\mathbf{u}^k-\mathbf{u}^*\|^2_{\M{H}}-c_0\|\M{M}(\mathbf{u}^k-\tilde{\mathbf{u}}^{k})\|^2_{\M{H}}.
\end{align}
Summarizing \eqref{Proof:Theorem:point-wise-rate:8} over $k=0,\cdots,K$, it gives that $\forall \M{u}^*\in\mathcal{M}^*$
\begin{align*}
\sum_{k=0}^{K}c_0\|\M{M}(\mathbf{u}^k-\tilde{\mathbf{u}}^{k})\|^2_{\M{H}}\leq \|\M{u}^0-\M{u}^*\|^2_{\M{H}}<\infty.
\end{align*}
It follows from \eqref{eq:point-wise-rate1} that the sequence $\{\|\M{M}(\M{u}^k-\tilde{\M{u}}^k)\|^2_{\M{H}}\}$ is monotonically non-increasing. Thus, we have
$$
\|\M{M}(\M{u}^k-\tilde{\M{u}}^k)\|^2_{\M{H}}\leq\frac{\|\M{u}^0-\M{u}^*\|^2_{\M{H}}}{c_0(k+1)}.
$$
Since $\M{u}^{k+1}-\M{u}^k=-\M{M}(\M{u}^k-\tilde{\M{u}}^k) $, the $o({1}/{k})$ rate of $\|\M{u}^{k+1}-\M{u}^k\|^2_{\M{H}}$ follows from \cite[Proposition 1]{PGEXTR}. By \eqref{ADDproof}, the $o({1}/{k})$ rate of the first-order optimality residual hold.
\end{IEEEproof}

\section{Proof of Theorem \ref{TH3}}\label{Appendix:Theorem3}
\begin{IEEEproof}
Since $\{\M{u}^k\}$ and $\{\tilde{\M{u}}^k\}$ converge to $\M{u}^{\infty}$ and $\{\|\M{d}^k\|\}$ converges to $0$ as $k\rightarrow \infty$, it from \eqref{Proof:theorem1:5} that $\{\M{v}^k\}$ and $\{\tilde{\M{v}}^k\}$ also converges to $\M{u}^{\infty}$ as $k\rightarrow \infty$. By \eqref{EXDiPGM}, one has
\begin{align*}
&\M{0}\in \partial G(\tilde{\M{y}}^{k})+\nabla F(\M{x}^k)-\M{V}\bm{\alpha}^k+\Gamma^{-1}(\tilde{\M{y}}^{k}-\M{x}^k),\\
&\M{0}=\beta{\M{V}}\tilde{\M{y}}^{k}+(\tilde{\M{z}}^k-\bm{\alpha}^k),
\end{align*}
which implies that
\begin{align*}
&\nabla F(\tilde{\M{y}}^k)-\nabla F(\M{x}^k)-\M{V}(\tilde{\M{z}}^k-\bm{\alpha}^k)-\Gamma^{-1}(\tilde{\M{y}}^{k}-\M{x}^k)\\
&\in \partial G(\tilde{\M{y}}^{k}) + \nabla F(\tilde{\M{x}}^k)-\M{V}{\M{z}}^k,
\end{align*}
and $-\frac{1}{\beta}(\tilde{\M{z}}^k-\bm{\alpha}^k)={\M{V}}\tilde{\M{y}}^{k}$.
Therefore, it holds that
\begin{align*}
&\mathrm{dist}^2(\M{0},\mathcal{J}(\tilde{\M{v}}^k))\leq\|\nabla F(\tilde{\M{y}}^k)-\nabla F(\M{x}^k)-\M{V}(\tilde{\M{z}}^k-\bm{\alpha}^k)\\
&\quad-\Gamma^{-1}(\tilde{\M{y}}^{k}-\M{x}^k)\|^2+\|\frac{1}{\beta}(\tilde{\M{z}}^k-\bm{\alpha}^k)\|^2\\
&\leq(3\max_i\{L^2_i\}+3\sigma_m(\Gamma^2))\|\tilde{\M{y}}^k-\M{x}^k\|^2\\
&\quad+(3\sigma_M(\M{V}^2)+\frac{1}{\beta^2})\|\tilde{\M{z}}^k-\bm{\alpha}^k\|^2\\
&\leq \kappa_2^2\|\tilde{\M{v}}^k-\M{u}^k\|^2_{\widehat{\M{H}}_1},
\end{align*}
where $\kappa_2^2=\frac{\max\{3\sigma_M(\M{V}^2)+\frac{1}{\beta^2},3\max_i\{L^2_i\}+3\sigma_m(\Gamma^{-2})\}}{\sigma_m(\widehat{\M{H}}_1)}$. Since the KKT mapping $\mathcal{J}$ is metrically subregular at $(\mathbf{u}^{\infty},\mathbf{0})$, there exists $\kappa>0$ and $\epsilon>0$ such that when $\tilde{\M{v}}^k\in\mathcal{B}_{\epsilon}(\mathbf{u}^{\infty})$
\begin{align*}
&\mathrm{dist}(\tilde{\M{v}}^{k},\mathcal{M}^*)=\mathrm{dist}(\tilde{\M{v}}^{k},\mathcal{T}^{-1}(0))\\
&\leq \kappa \mathrm{dist}(0,\mathcal{T}(\tilde{\M{v}}^{k}))\leq \kappa\kappa_2\|\tilde{\M{v}}^{k}-\M{u}^k\|_{\widehat{\M{H}}_1}.
\end{align*}
Since $\{\M{u}^k\}$ and $\{\tilde{\M{v}}^k\}$ converge to $\M{u}^{\infty}$ as $k\rightarrow \infty$, there exists $\bar{k}\geq0$ such that
\begin{align}\label{Proof:theorem3:1}
\mathrm{dist}(\tilde{\M{v}}^{k},\mathcal{M}^*)\leq \kappa\kappa_2\|\tilde{\M{v}}^{k}-\M{u}^k\|_{\widehat{\M{H}}_1}, \forall k\geq \bar{k}.
\end{align}
Since
Note that $0\prec\widehat{\M{H}}_1\prec\M{H}$, one has
\begin{align}\label{Proof:theorem3:2}
&\mathrm{dist}(\tilde{\M{v}}^{k},\mathcal{M}^*)\geq \frac{1}{c_1} \mathrm{dist}_{\widehat{\M{H}}_1}(\tilde{\M{v}}^{k},\mathcal{M}^*)\nonumber\\
&\geq \frac{1}{c_1} (\mathrm{dist}_{\widehat{\M{H}}_1}(\M{u}^{k},\mathcal{M}^*)-\|\tilde{\M{v}}^{k}-\M{u}^k\|_{\widehat{\M{H}}_1}),
\end{align}
where $c_1=\sigma_M^{1/2}(\M{H}) $. Combining \eqref{Proof:theorem3:1} and \eqref{Proof:theorem3:2}, it gives that
\begin{align}\label{Proof:theorem3:3}
&\mathrm{dist}_{\M{H}}(\M{u}^k,\mathcal{M}^*)\leq\frac{c_1}{c_2}\mathrm{dist}_{\widehat{\M{H}}_1}(\M{u}^k,\mathcal{M}^*)\nonumber\\
&\leq\frac{c_1(\kappa\kappa_2c_1+1)}{c_2}\|\tilde{\M{v}}^{k}-\M{u}^k\|_{\widehat{\M{H}}_1},\forall k\geq \bar{k},
\end{align}
where $c_2=\sigma_m^{1/2}(\widehat{\M{H}}_1) $. It follows from \eqref{Theorem1:P1} that for $\forall k\geq0$
\begin{align}\label{Proof:theorem3:V1}
\|\mathbf{v}^{k+1}-\mathbf{u}^{\infty}\|^2_{\M{H}}\leq \|\mathbf{u}^k-\mathbf{u}^{\infty}\|^2_{\M{H}}-\|\tilde{\M{v}}^{k}-\M{u}^k\|^2_{\widehat{\M{H}}_1}.
\end{align}
Together with \eqref{Proof:theorem3:3}, one has
\begin{align*}
&\mathrm{dist}_{\M{H}}^2(\M{u}^k,\mathcal{M}^*)-\mathrm{dist}_{\M{H}}^2(\M{v}^{k+1},\mathcal{M}^*)\geq \|\tilde{\M{v}}^{k}-\M{u}^k\|^2_{\widehat{\M{H}}_1}\\
& \geq \frac{c_2^2}{c^2_1(\kappa\kappa_2c_1+1)^2}  \mathrm{dist}^2_{\M{H}}(\M{u}^k,\mathcal{M}^*), \forall k\geq \bar{k},
\end{align*}
which implies that
\begin{align}\label{Proof:theorem3:4}
\mathrm{dist}_{\M{H}}(\M{v}^{k+1},\mathcal{M}^*) \leq  \bar{\vartheta}~ \mathrm{dist}_{\M{H}}(\M{u}^k,\mathcal{M}^*),\forall k\geq \bar{k},
\end{align}
where $\bar{\vartheta}=\sqrt{1-\frac{c_2^2}{c^2_1(\kappa\kappa_2c_1+1)^2}}<1$. By \eqref{Proof:theorem3:V1} and the fact $0\prec\widehat{\M{H}}_1\prec\M{H}$, one has
\begin{align*}
\frac{c_2}{c_1}\|\tilde{\M{v}}^{k}-\M{u}^k\|^2_{\M{H}}\leq\|\tilde{\M{v}}^{k}-\M{u}^k\|^2_{\widehat{\M{H}}_1}\leq\mathrm{dist}^2_{\M{H}}(\M{u}^k,\mathcal{M}^*),\forall k\geq \bar{k}.
\end{align*}
It deduces that
\begin{align}\label{Proof:theorem3:5}
\|\tilde{\M{v}}^{k}-\M{u}^k\|_{\M{H}}\leq \omega~ \mathrm{dist}_{\M{H}}(\M{u}^k,\mathcal{M}^*),\forall k\geq \bar{k},
\end{align}
where $\omega=\sqrt{c_1/c_2}$. By the triangle inequality, one has
\begin{align}\label{Proof:theorem3:6}
&\|\mathbf{u}^k-\mathcal{P}_{\mathcal{M}^*}^{\M{H}}(\M{v}^{k+1})\|_{\M{H}}\leq \mathrm{dist}_{\M{H}}(\mathbf{u}^k,\mathcal{M}^*)\nonumber\\
&\quad+\|\mathcal{P}_{\mathcal{M}^*}^{\M{H}}(\mathbf{u}^k)-\mathcal{P}_{\mathcal{M}^*}^{\M{H}}({\mathbf{v}}^{k+1})\|_{\mathrm{H}},\forall k\geq \bar{k}.
\end{align}
Then, by the Lipschitz continuity of $\mathcal{P}_{\mathcal{M}^*}^{\M{H}}$, it holds that
\begin{align*}
\|\mathcal{P}_{\mathcal{M}^*}^{\M{H}}(\mathbf{u}^k)-\mathcal{P}_{\mathcal{M}^*}^{\M{H}}({\mathbf{v}}^{k+1})\|_{\M{H}}\leq \|\mathbf{u}^k-{\mathbf{v}}^{k+1}\|_{\M{H}},
\end{align*}
which together with \eqref{Proof:theorem3:5} and \eqref{Proof:theorem3:6} deduces that
\begin{align}\label{Proof:theorem3:7}
\|\mathbf{u}^k-\mathcal{P}_{\mathcal{M}^*}^{\M{H}}({\mathbf{v}}^{k+1})\|_{\M{H}}\leq (1+\omega) \mathrm{dist}_{\M{H}}(\mathbf{u}^k,\mathcal{M}^*).
\end{align}
By \eqref{Cr1}, \eqref{Proof:theorem1:5} and the triangle inequality, it holds that
\begin{align*}
&\|\mathbf{u}^{k+1}-\mathcal{P}_{\mathcal{M}^*}^{\M{H}}({\mathbf{v}}^{k+1})\|_{\M{H}}\\
&\leq \|{\mathbf{v}}^{k+1}-\mathbf{u}^{k+1}\|_{\M{H}}+\mathrm{dist}_{\M{H}}({\mathbf{v}}^{k+1},\mathcal{M}^*)\\
&\leq \bar{\mu}\|\mathbf{d}^k\|+\mathrm{dist}_{\M{H}}({\mathbf{v}}^{k+1},\mathcal{M}^*)\\
&\leq \mu\varrho_k\|\mathbf{u}^{k}-\mathbf{u}^{k+1}\|_{\M{H}}+\mathrm{dist}_{\M{H}}({\mathbf{v}}^{k+1},\mathcal{M}^*)\\
&\leq  \mu\varrho_k(\|\mathbf{u}^{k+1}-\mathcal{P}_{\mathcal{M}^*}^{\M{H}}({\mathbf{v}}^{k+1})\|\\
&\quad+\|\mathbf{u}^{k}-\mathcal{P}_{\mathcal{M}^*}^{\M{H}}({\mathbf{v}}^{k+1})\|)
+\mathrm{dist}_{\M{H}}({\mathbf{v}}^{k+1},\mathcal{M}^*),
\end{align*}
where $\mu={\bar{\mu}}/{c_2}$. Then, by the fact that
$$
\|\mathbf{u}^{k+1}-\mathcal{P}_{\mathcal{M}^*}^{\M{H}}({\mathbf{v}}^{k+1})\|_{\M{H}}\geq \mathrm{dist}_{\M{H}}(\mathbf{u}^{k+1},\mathcal{M}^*)
$$
and \eqref{Proof:theorem3:7}, it can be deduced that when $k\geq\bar{k}$
\begin{align*}
&(1-\mu\varrho_k)\mathrm{dist}_{\M{H}}(\mathbf{u}^{k+1},\mathcal{M}^*)\\
&\leq {\mu}\varrho_k(1+\omega) \mathrm{dist}_{\M{H}}(\mathbf{u}^{k},\mathcal{M}^*)+
\mathrm{dist}_{\M{H}}({\mathbf{v}}^{k+1},\mathcal{M}^*),
\end{align*}
which together with \eqref{Proof:theorem3:4}, implies that when $k\geq\bar{k}$
$$
\mathrm{dist}_{\M{H}}(\mathbf{u}^{k+1},\mathcal{M}^*)\leq \underbrace{\frac{{\mu}(1+\omega)\varrho_k+\bar{\vartheta} }{(1-{\mu}\varrho_k)} }_{:=\vartheta_k} \mathrm{dist}_{\M{H}}(\mathbf{u}^{k},\mathcal{M}^*).
$$
It follows that if
$$
\sup_{k\geq\bar{k}}\{\varrho_k\}<\frac{1-\bar{\vartheta}}{\mu(2+\omega)},
$$
then one has $\sup_{k\geq\bar{k}}\{\varrho_k\}<1$. This completes the proof.
\end{IEEEproof}

\newpage
\begin{IEEEbiography}[{\includegraphics[width=1.3in,height=1.2in,clip,keepaspectratio]{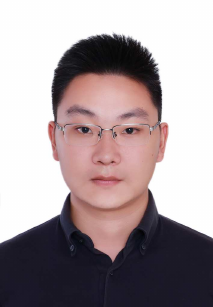}}]{Luyao Guo}
received the B.S. degree in information and computing science
from Shanxi University, Taiyuan, China, in 2020. He is currently pursuing the Ph.D.
degree in applied mathematics with the Jiangsu Provincial Key Laboratory of Networked Collective
Intelligence, School of Mathematics, Southeast University, Nanjing, China.

His current research focuses on distributed optimization and learning.
\end{IEEEbiography}

\vspace{-50 mm}

\begin{IEEEbiography}[{\includegraphics[width=1.3in,height=1.2in,clip,keepaspectratio]{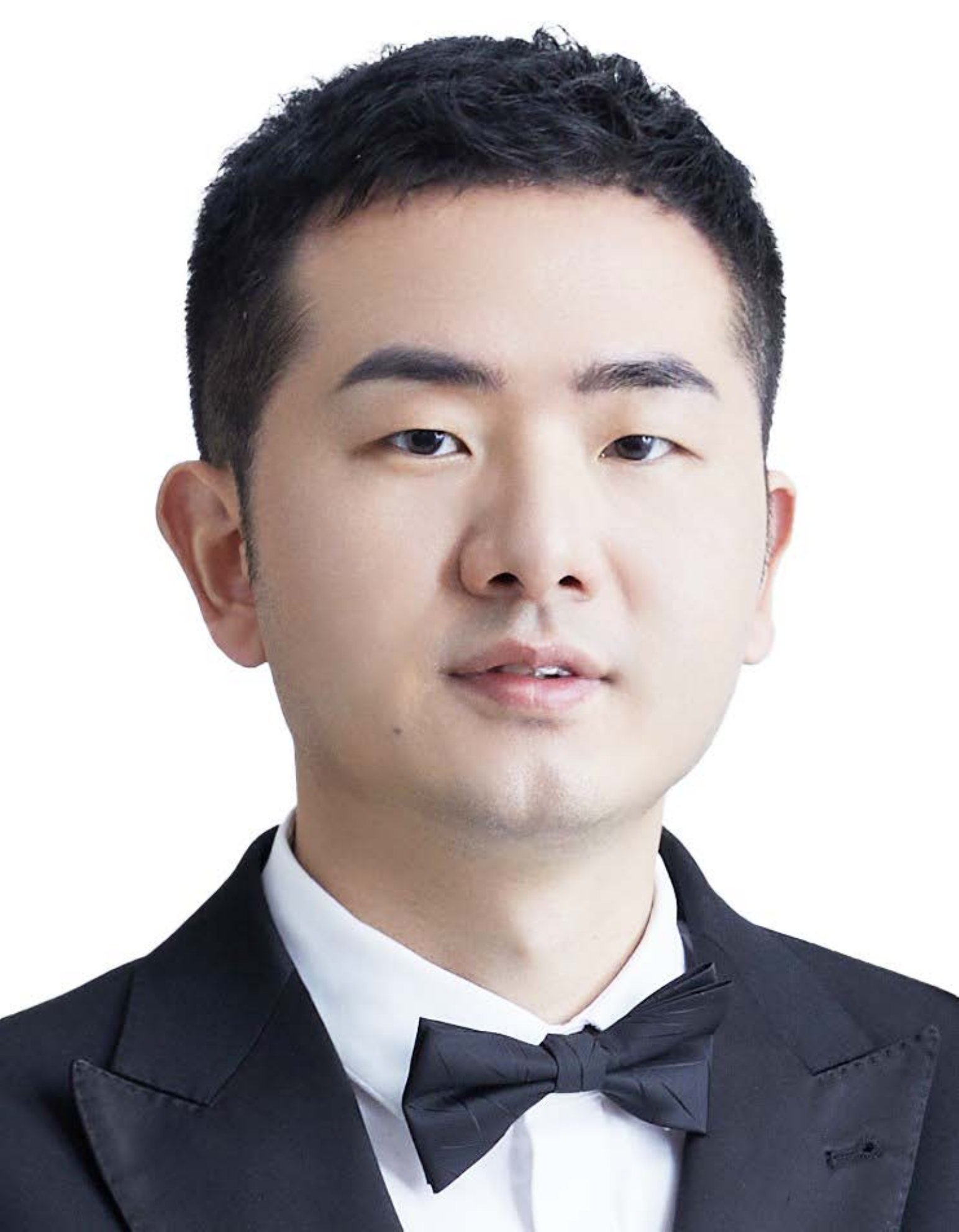}}]{Xinli Shi (Senior Member, IEEE)}
received the B.S. degree in software engineering, the M.S. degree in applied mathematics, and the Ph.D. degree in control science and engineering from Southeast University, Nanjing, China, in 2013, 2016, and 2019, respectively.

He held a China Scholarship Council Studentship for one-year study with the University of Royal Melbourne Institute of Technology, Melbourne, VIC, Australia, in 2018. He is currently an Associate Professor with the School of Cyber Science and Engineering, Southeast University. His current research interests include distributed optimization, nonsmooth analysis, and network control systems.

Dr. Shi was the recipient of the Outstanding Ph.D. Degree Thesis Award from Jiangsu Province, China.
\end{IEEEbiography}
\vspace{-50 mm}
\begin{IEEEbiography}[{\includegraphics[width=1.3in,height=1.2in,clip,keepaspectratio]{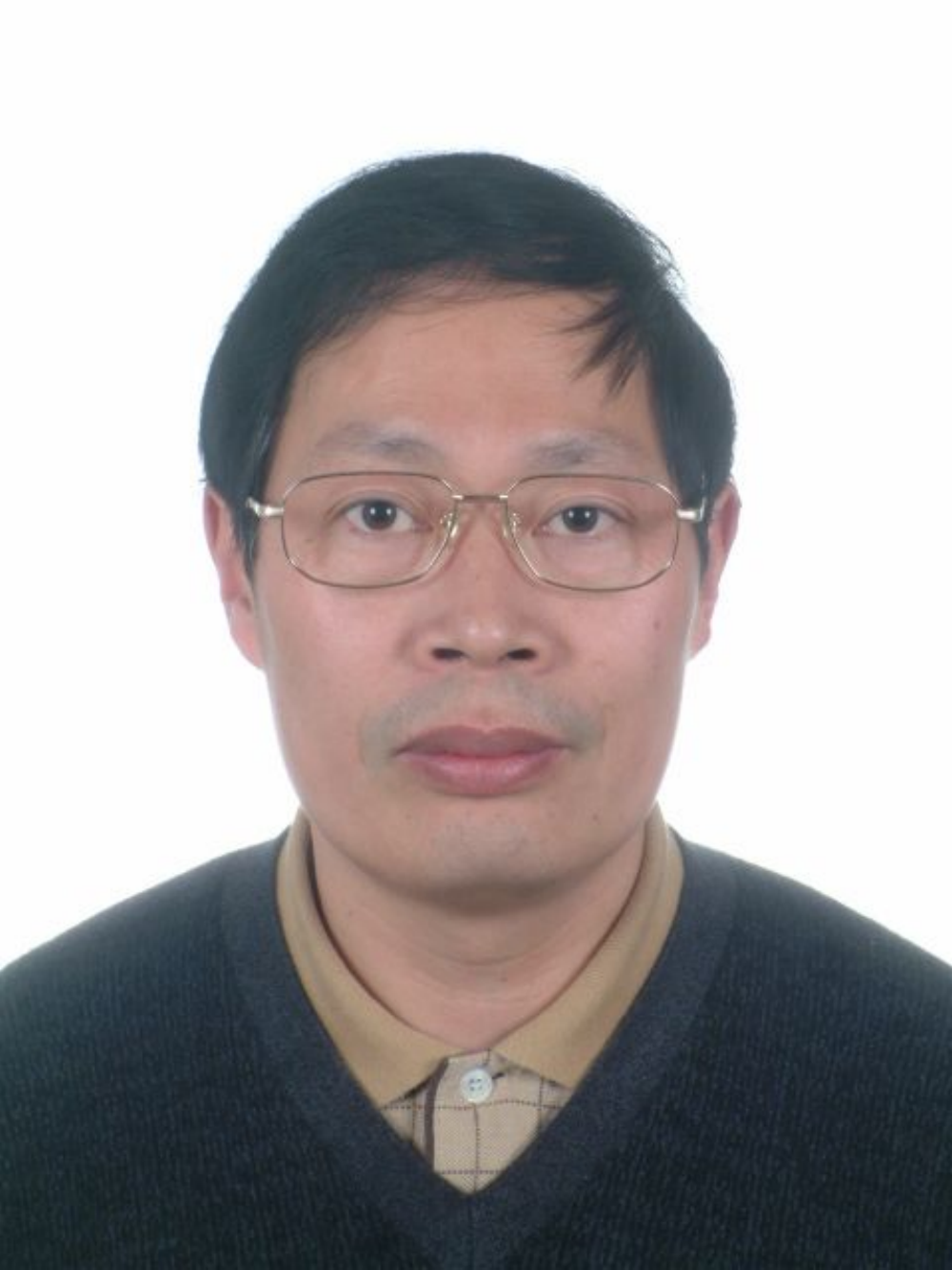}}]{Jinde Cao (Fellow, IEEE)}

received the B.S. degree in mathematics from Anhui Normal University, Wuhu, China, in 1986, the M.S. degree in mathematics/applied mathematics from Yunnan University, Kunming, China, in 1989, and the Ph.D. degree in mathematics/applied mathematics from Sichuan University, Chengdu, China, in 1998.

He is an Endowed Chair Professor, the Dean of the School of Mathematics, and the Director of the Jiangsu Provincial Key Laboratory of Networked Collective Intelligence of China and the Research Center for Complex Systems and Network Sciences with Southeast University, Nanjing, China.

Prof. Cao was a recipient of the National Innovation Award of China, the Gold medal of Russian Academy of Natural Sciences, the Obada Prize, and the Highly Cited Researcher Award in Engineering, Computer Science, and Mathematics by Thomson Reuters/Clarivate Analytics. He is elected as a foreign member of Russian Academy of Sciences, a member of the Academy of Europe, a foreign member of Russian Academy of Engineering, a member of the European Academy of Sciences and Arts, a foreign member of Russian Academy of Natural Sciences, a foreign fellow of Pakistan Academy of Sciences, a fellow of African Academy of Sciences, a foreign Member of the Lithuanian Academy of Sciences, and an IASCYS academician.
\end{IEEEbiography}

\vspace{-50 mm}

\begin{IEEEbiography}[{\includegraphics[width=1.3in,height=1.2in,clip,keepaspectratio]{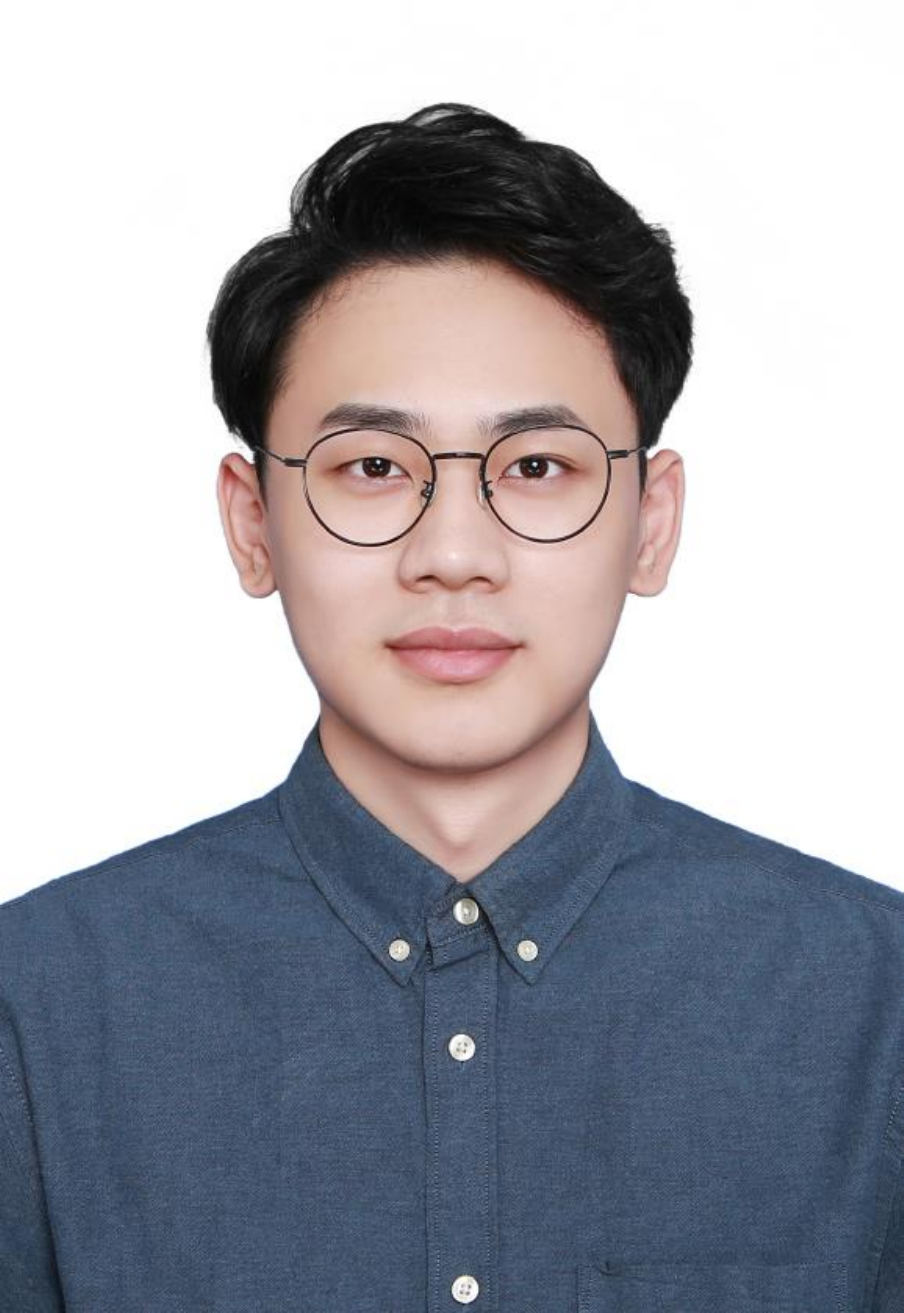}}]{Zihao Wang}
received the B.E. degree in network engineering from Nanjing Xiaozhuang University, Nanjing, China, in 2020. He is currently pursuing the M.E. degree in cyberspace security with the
School of Cyber Science and Engineering, Southeast University, Nanjing.

His current research focuses on privacy-preserving distributed optimization of networked multi-agent systems.
\end{IEEEbiography}

\end{document}